\documentclass[12pt]{amsart}
\usepackage{amssymb}
\usepackage{amsmath}
\usepackage{longtable}
\newcommand\g{{\mathfrak g}}
\newcommand\h{{\mathfrak h}}
\newcommand{\f}{\mathfrak{f}}
\renewcommand{\b}{\mathfrak{b}}
\newcommand\gl{\mathfrak{gl}}
\renewcommand{\a}{\mathfrak{a}}

\newcommand\m{\mathfrak m}
\renewcommand\l{\mathfrak l}
\newcommand\n{\mathfrak n}
\newcommand\q{\mathfrak q}
\newcommand\z{\mathfrak z}
\newcommand\s{\mathfrak s}

\renewcommand\u{{\mathfrak u}}
\newcommand\skewperp{\angle}
\renewcommand{\t}{\mathfrak{t}}
\newcommand\codim{\operatorname{codim}}

\newcommand\tr{\operatorname{tr}}
\newcommand\td{\operatorname{tr.deg}}
\newcommand\im{\operatorname{im}}

\newcommand\Quot{\operatorname{Quot}}

\newcommand\Aut{\operatorname{Aut}}

\newcommand\C{\mathbb C}
\newcommand\X{\mathfrak X}
\newcommand\Q{\mathbb Q}
\newcommand\R{\mathbb R}

\newcommand\Z{\mathbb Z}
\newcommand\A{\mathfrak A}
\newcommand\V{\mathcal V}
\newcommand\Rad{\operatorname{R}}
\newcommand\Gr{\operatorname{Gr}}
\renewcommand\O{\operatorname{O}}
\renewcommand\sl{\mathfrak{sl}}
\newcommand\so{\mathfrak{so}}
\renewcommand\sp{\mathfrak{sp}}
\newcommand\spin{\mathfrak{spin}}
\newcommand\GL{\mathop{\rm GL}\nolimits}

\newcommand\SL{\mathop{\rm SL}\nolimits}
\newcommand\Sp{\mathop{\rm Sp}\nolimits}
\newcommand\SO{\mathop{\rm SO}\nolimits}

\newcommand\Span{\operatorname{Span}}
\newcommand\Supp{\operatorname{Supp}}

\newcommand{\Ad}{\mathop{\rm Ad}\nolimits}

\newcommand{\rank}{\mathop{\rm rk}\nolimits}

\newcommand\Int{\mathop{\rm Int}\nolimits}
\renewcommand{\Ad}{\mathop{\rm Ad}\nolimits}

\newcommand\quo{/\!/}
\newtheorem{Thm}{Theorem}[subsection]
\newtheorem{Prop}[Thm]{Proposition}
\newtheorem{Cor}[Thm]{Corollary}
\newtheorem{Lem}[Thm]{Lemma}
\theoremstyle{definition}

\newtheorem{defi}[Thm]{Definition}
\newtheorem{Rem}[Thm]{Remark}
\newtheorem{Alg}[Thm]{Algorithm}
\unitlength=1mm   \numberwithin{equation}{section}
\numberwithin{table}{section} \oddsidemargin=0cm
\author{Ivan V. Losev}
\title{Computation of  Weyl groups of  $G$-varieties}
\thanks{{\it Key words and phrases}: reductive groups,  homogeneous spaces, Weyl groups,
spherical varieties} \thanks{{\it 2000 Mathematics Subject
Classification.} 14M17, 14R20}
\begin{document}
\begin{abstract}
Let $G$ be a connected reductive group. To any irreducible
$G$-variety one assigns a certain linear group generated by
reflections called the {\it Weyl group}. Weyl groups play an
important role in the study of embeddings of homogeneous spaces. We
establish algorithms for computing Weyl groups for homogeneous
spaces and affine homogeneous vector bundles. For some special
classes of $G$-varieties (affine homogeneous vector bundles of
maximal rank, affine homogeneous spaces, homogeneous spaces of
maximal rank with discrete group of central automorphisms) we
compute Weyl groups more or less explicitly.
\end{abstract}
\maketitle \tableofcontents
\section{Introduction}\label{SECTION_Intro}
In the whole paper the base field is the field $\C$ of complex
numbers. In this section $G$ denotes a connected reductive group. We
fix a Borel subgroup $B\subset G$ and a maximal torus $T\subset B$.
$X$ denotes an irreducible $G$-variety.

\subsection{Definition of the Weyl group of a $G$-variety}\label{SUBSECTION_Intro1}
The main object considered in the present paper is the Weyl group of
an irreducible $G$-variety. Before giving the definition we would
like to make some historical remarks.

The Weyl group of $G$ was essentially defined in Herman Weyl's paper
\cite{Weyl}. Little bit later, E. Cartan, \cite{Cartan}, generalized
the notion of the Weyl group to symmetric spaces (so called, {\it
little Weyl groups}). From the algebraic viewpoint a symmetric space
is a homogeneous space $G/H$, where $(G^\sigma)^\circ\subset
H\subset G^\sigma$ for some involutory automorphism $\sigma$ of $G$.

To move forward we need the notion of complexity.

\begin{defi}\label{Def:0.1.1}
The {\it complexity} of $X$ is the codimension of a general
$B$-orbit in $X$, or, equivalently, $\td \C(X)^{B}$. We denote the
complexity of $X$ by $c_{G}(X)\index{cgx@$c_{G}(X)$}$. A normal
irreducible $G$-variety of complexity $0$ is said to be {\it
spherical}.
\end{defi}

In particular, every symmetric space is a spherical $G$-variety,
\cite{Vust1}. In \cite{Brion} Brion constructed the Weyl group for a
spherical homogeneous space $G/H$ with $\#N_G(H)/H<\infty$. Brion's
Weyl group generalizes that of symmetric spaces. In view of results
of \cite{BP}, the restriction $\#N_G(H)/H<\infty$ is not essential.
Knop, \cite{Knop1}, found another way to define the Weyl group for
an arbitrary irreducible $G$-variety also generalizing the Weyl
group of a symmetric space. In \cite{Knop4} he extended Brion's
definition to arbitrary $G$-varieties. Finally, in \cite{Knop3} Knop
gave a third definition of the Weyl group and proved the equivalence
of all three definitions.

Now we are going to introduce the definition of the Weyl group
following \cite{Knop4}.

 Consider
the sublattice $\X_{G,X}\subset \X(T)\index{xgx@$\X_{G,X}$}$
consisting of the weights of $B$-semiinvariant functions from
$\C(X)$. It is called the {\it weight lattice} of $X$.

\begin{defi}\label{Def:0.1.2} Put $\a_{G,X}=\X_{G,X}\otimes_\Z\C\index{agx@$\a_{G,X}$}$. We call the subspace
$\a_{G,X}\subset \t^*$  the {\it Cartan space} of $X$. The dimension
of $\a_{G,X}$ is called the {\it rank} of $X$ and is denoted by
$\rank_G(X)\index{rgx@$\rank_{G}(X)$}$.
\end{defi}

 Fix a $W(\g)$-invariant scalar product on $\t(\R)$. This induces
the scalar product on $\a_{G,X}(\R):=\X_{G,X}\otimes_\Z\R$  and on
$\a_{G,X}(\R)^*$. The Weyl group of $X$ will act on $\a_{G,X}$
preserving the weight lattice and the scalar product. To define the
action we will describe its Weyl chamber. To this end we need the
notion of a central $G$-valuation.

\begin{defi}\label{Def:0.1.3}
By a $G$-valuation of $X$ we mean a discrete $\R$-valued
$G$-invariant valuation of $\C(X)$. A $G$-valuation  is called {\it
central} if it vanishes on $\C(X)^B$.
\end{defi}

In particular,  if $X$ is spherical, then any $G$-invariant
valuation is central. A central valuation $v$ determines an element
$\varphi_v\in \a_{G,X}(\R)^*$ by $\langle
\varphi_v,\lambda\rangle=v(f_\lambda)$, where $\lambda\in \X_{G,X}$,
$f_\lambda\in \C(X)^{(B)}_\lambda\setminus \{0\}$.  The element
$\varphi_v$ is well-defined because $v$ is central.

\begin{Thm}\label{Thm:0.1.4}
\begin{enumerate}
\item The map $v\mapsto \varphi_v$ is injective. Its image is a finitely
generated rational convex cone in $\a_{G,X}(\R)^*$ called the {\it
central valuation cone} of $X$ and denoted by
$\V_{G,X}\index{vgx@$\V_{G,X}$}$.
\item The central valuation cone is simplicial (that is, there are linearly
independent vectors $\alpha_1,\ldots,\alpha_s\in \a_{G,X}(\R)$ such
that the cone coincides with $\{x| \langle\alpha_i,x\rangle\geqslant
0, i=\overline{1,s}\}$). Moreover, the reflections corresponding to
its facets generate a finite group. This group is called the {\it
Weyl group} of $X$  and is denoted by
$W_{G,X}\index{wgx@$W_{G,X}$}$.
\item The lattice $\X_{G,X}\subset \a_{G,X}(\R)$ is $W_{G,X}$-stable.
\end{enumerate}
\end{Thm}

The proof of the first part of the theorem is relatively easy. It is
obtained (in a greater generality) in \cite{Knop4}, Korollare 3.6,
4.2, 5.2, 6.5. The second assertion is much more complicated. It was
proved by Brion in \cite{Brion} in the spherical case. Knop,
\cite{Knop4}, used Brion's result to prove the assertion in general
case. Later, he gave an alternative  proof in \cite{Knop3}. The
third assertion of Theorem \ref{Thm:0.1.4} follows easily from the
construction of the Weyl group in \cite{Knop3}.

Note that the Weyl group does not depend on the scalar product used
in its definition. Indeed, the set of $W$-invariant scalar products
on $\t(\R)^*$ is convex. The Weyl group fixes $\X_{G,X}$ so does not
change under small variations of the scalar product.

\subsection{Main problem}\label{SUBSECTION_Intro2}
Our general problem is to find an algorithm  computing
$\a_{G,X},W_{G,X}$ for an irreducible $G$-variety $X$. However, for
such an algorithm to exist, the variety $X$ should have some good
form. It is reasonable to restrict ourselves to the following two
classes of $G$-varieties:
\begin{enumerate}
\item Homogeneous spaces $G/H$, where $H$ is an algebraic subgroup of $G$.
\item Homogeneous vector bundles over affine homogeneous spaces (=affine homogeneous
vector bundles) $G*_HV$. Here $H$ is a reductive subgroup of $G$ and
$V$ is an $H$-module.
\end{enumerate}

There are several reasons to make these restrictions. First of all,
these $G$-varieties have  "group theoretic" and "representation
theoretic" structure, so one may hope to find algorithms with
"group-" or "representation theoretic" steps. Secondly, the
computation of the Cartan space and the Weyl group of an arbitrary
$G$-variety can be reduced to the computation for homogeneous
spaces. Namely, for an irreducible $G$-variety $X$ and a point $x\in
X$ in general position the equalities $\a_{G,X}=\a_{G,Gx}$,
$W_{G,X}=W_{G,Gx}$ hold (see Proposition \ref{Prop:1.2.1}).
Moreover, if $X$ is affine and smooth and $x\in X$ is a point with
closed $G$-orbit, then $\a_{G,X}=\a_{G,G*_HV}, W_{G,X}=W_{G,G*_HV}$,
where $H=G_x$ and $V$ is the slice module at $x$, that is
$V=T_xX/\g_*x$ (Corollary \ref{Cor:1.2.3}).

So the main results of the paper are algorithms computing the Cartan
spaces and the Weyl groups for $G$-varieties of types (1) and (2).
Moreover, we compute the Weyl groups of affine homogeneous spaces
more or less explicitly.

Our algorithms are quite complicated so we do not give them here.
They will be presented (in a brief form) in
Section~\ref{SECTION_algorithm}. Roughly speaking, all our steps
consist of  computing  some "structure characteristics" for pairs
(an algebraic Lie algebra, a subalgebra), (a reductive Lie algebra,
a module over this algebra). The computation of the normalizer or
the unipotent radical is an example of an operation from the first
group. An operation from  the second group is, for instance, the
decomposition of the restriction of an irreducible representation to
a Levi subalgebra together with the determination of all highest
vectors of the restriction.

\subsection{Motivations and known results}\label{SUBSECTION_Intro3}
Our main motivation comes from the theory of embeddings of
homogeneous spaces.

One may say that the theory of algebraic transformation groups
studies the category of varieties acted on by  some  algebraic
group. Because of technical reasons, one usually considers actions
of connected reductive groups on normal irreducible $G$-varieties.
The first problem in the study of a category is the classification
of its objects up to an isomorphism. In our case, the problem may be
divided into two parts, birational and regular. The birational part
is the classification of $G$-varieties up to a birational
equivalence, or, in the algebraic setting, the classification of all
finitely generated fields equipped with an action of $G$ by
automorphisms. An important special case here is the birational
classification of quasi-homogeneous $G$-varieties, i.e. those
possessing an open $G$-orbit. Of course, an equivalent problem is
the classification of algebraic subgroups of $G$ up to conjugacy.

The regular part is the classification of $G$-varieties in a given
class of birational equivalence. In the quasi-homogeneous case this
is equivalent to the classification of all open embeddings of a
given homogeneous space into normal varieties. The program to
perform this classification was proposed by Luna and Vust
(\cite{LV}). Note that in that paper only the quasi-homogeneous case
was considered. However, the Luna-Vust theory can be generalized to
the general case, see \cite{Timashev2}. Using the Luna-Vust theory
one obtains a combinatorial (in a certain sense) description of
$G$-varieties of complexity not exceeding 1. The spherical case was
considered already in \cite{LV}. More self-contained and plain
exposition is given, for example, in \cite{Knop5}. The case of
complexity 1 is due to Timashev, \cite{Timashev2}. The
classification for complexity greater than 1 seems to be wild.

Now we sketch  the classification theory  of spherical variety.
Clearly, a spherical $G$-variety has an open orbit. So the
birational part of the classification is just the classification of
all spherical homogeneous spaces. To describe all embeddings of a
given spherical homogeneous space one needs to know the following
data:

\begin{enumerate}
\item The rational vector space $\a_{G,X}(\Q):=\X_{G,X}\otimes_\Z\Q$.
\item The (central) valuation cone $\mathcal{V}_{G,X}\subset \a_{G,X}(\R)^*$.
\item Certain {\it colored} vectors in $\a_{G,X}(\Q)^*$ that are in one-to-one correspondence
with $B$-divisors of the spherical homogeneous space. Namely, given
a prime $B$-divisor $D$ we define the colored vector
$\varphi_D\index{zzzv@$\varphi_D$}$ by $\langle
\varphi_D,\lambda\rangle=\operatorname{ord}_D(f_\lambda)$, where
$\lambda\in \X_{G,X}$, $f_\lambda\in \C(X)^{(B)}_\lambda\setminus
\{0\}$.
\end{enumerate}

It turns out that normal embeddings of a spherical homogeneous space
$X$ are in one-to-one correspondence with certain {\it admissible}
sets of cones in $\a_{G,X}(\R)^*$. Every cone from an admissible set
is generated by colored vectors and elements of $\V_{G,X}$ and
satisfies certain combinatorial requirements. An admissible set is
one that satisfies some additional requirements of combinatorial
nature. So the solution of the regular part of the classificational
problem consists in the determination of $\a_{G,X}(\Q)$, $\V_{G,X}$
and colored vectors.

Note that  $\a_{G,X}(\Q)=\a_{G,X}\cap \t(\Q)^*$. The computation of
$\a_{G,X}$ is not very difficult. After $\a_{G,X}$ is computed one
can proceed to the computation of the valuation cone $\V_{G,X}$.
Despite of the fact that a group generated by reflections has
several Weyl chambers, the Weyl group determines the cone $\V_{G,X}$
uniquely. Namely, $\V_{G,X}$ is a unique Weyl chamber of $W_{G,X}$
containing the image of the negative Weyl chamber of $\t$ under the
projection $\t(\R)\rightarrow \a_{G,X}(\R)^*$.

Now we discuss results concerning the computation of Weyl groups and
Cartan spaces of $G$-varieties.

D.I. Panyushev in \cite{Panyushev2}, see also
\cite{Panyushev_thesis}, reduced
 the computation of weight lattices for $G$-varieties of two types
mentioned above to that for affine homogeneous spaces (in fact,
together with some auxiliary datum). In \cite{ranks}  Cartan spaces
for affine homogeneous spaces were computed.

Proceed to results on the computation of Weyl groups. They are
formulated in three possible ways:
\begin{enumerate}
\item In terms of the Weyl group itself.
\item In terms of the central valuation cone.
\item In terms of primitive linearly independent elements
$\beta_1,\ldots,\beta_r\in \X_{G,X}$ such that the central valuation
cone is given by the inequalities $\beta_i\leqslant 0$. We denote
the set $\{\beta_1,\ldots,\beta_r\}$ by
$\Pi_{G,X}\index{zzzq@$\Pi_{G,X}$}$. For spherical $X$ such elements
$\beta_i$ are called  spherical roots of $X$.
\end{enumerate}
The Weyl group is unit iff $X$ is {\it horospherical}, i.e., the
stabilizer of any point contains a maximal unipotent subgroup. In
full generality it was proved by Knop, \cite{Knop1}.

The other results (at least known to the author) on the computation
of Weyl group relate to the spherical case.

In \cite{Brion}, Brion proposed a technique allowing to extract the
Weyl group $W_{G,G/H}$, where $G/H$ is a spherical homogeneous
space, from the algebra $g\h g^{-1}$ for some special $g\in G$. It
is an open problem to describe the set of all suitable $g$.

As we have already mentioned above, the Weyl group of a symmetric
space coincides with its little Weyl group, see \cite{Knop1},
\cite{Brion}, \cite{Vust2}. The Weyl groups of spherical $G$-modules
were computed by Knop in \cite{Knop11}. Note that in that paper the
notation $W_V$ is used for the Weyl group of $V^*$.

There are also computations of Weyl groups for some other special
classes of spherical homogeneous spaces. Spherical roots for {\it
wonderful varieties}  of rank $2$ were computed in \cite{Wasserman}.
The computation is based on some structure theorems on wonderful
varieties. The computation for other interesting class of
homogeneous spaces can be found in \cite{Smirnov}. It uses the
method of formal curves, established in \cite{LV}.

Of all results mentioned above we use only Wasserman's (see,
however, Remark \ref{Rem:8.0.1}).

Finally, let us make a remark on the classification of spherical
varieties. The first step of the classification is describing all
spherical homogeneous spaces. Up to now there is only one approach
to this problem due to Luna, who applied it to classify spherical
subgroups in groups of type $A$ (a connected reductive group is said
to be of type $A$ if all simple ideals of its Lie algebra are of
type $A$). Using Luna's approach, the full classification for groups
of type $A-D$ (\cite{Bravi}) and a partial one for type $A-C$
(\cite{Pezzini}) were obtained. The basic idea of the Luna
classification is to establish a one-to-one correspondence between
spherical homogeneous spaces and certain combinatorial data that are
almost equivalent to those listed above. We note, however, that
Luna's approach does not allow to obtain the full classification
even for groups of type $A-C$. Besides, the computation of
combinatorial data (in particular, the Weyl group) for certain
homogeneous spaces plays an important role in this approach.

\subsection{The structure of the paper}\label{SUBSECTION_Intro4}






Every section is divided into subsections. Theorems, lemmas,
definitions etc. are numbered within each subsection, while formulae
and tables within each section. The first subsection of Sections
\ref{SECTION_prelim}-\ref{SECTION_Weyl_homogen} describes their
content in detail.

{\bf Acknowledgements.}   I wish to thank Prof. E.B. Vinberg, who
awoke my interest to the subject. Also I would like to thank D.A.
Timashev for stimulating discussions. Some parts of this paper were
written during my visit to Fourier University of Grenoble in June
2006. I express my gratitude to this institution and especially to
Professor M. Brion for hospitality.

\section{Notation and conventions}\label{SECTION_notation}
For an algebraic group denoted by a capital Latin letter we denote
its Lie algebra by the  corresponding small German letter. For
example, the Lie algebra of  $\widetilde{L}_{0}$ is denoted by
$\widetilde{\l}_0$.

By a unipotent Lie algebra we mean the Lie algebra of a unipotent
algebraic group.

{\bf $H$-morphisms, $H$-subvarieties, etc.} Let $H$ be an algebraic
group. We say that a variety $X$ is an $H$-variety if an action of
$H$ on $X$ is given. By an $H$-subset (resp., subvariety) in a given
$H$-variety we mean an $H$-stable subset (resp., subvariety). A
morphism  of $H$-varieties is said to be an $H$-morphism if it is
$H$-equivariant. The term "$H$-bundle" means a principal bundle with
the structure group $H$.

{\bf Borel subgroups and maximal tori.} While considering a
reductive group $G$, we always fix its Borel subgroup $B$ and a
maximal torus $T\subset B$. In accordance with this choice, we fix
the root system $\Delta(\g)$ and the system of simple roots
$\Pi(\g)$ of $\g$. The Borel subgroup of $G$ containing $T$ and
opposite to $B$ is denoted by $B^-$.

If $G_1,G_2$ are reductive groups with the fixed Borel subgroups
$B_i\subset G_i$ and maximal tori $T_i\subset B_i$, then we take
$B_1\times B_2,T_1\times T_2$ for the fixed Borel subgroup and
maximal torus in $G_1\times G_2$.

Suppose $G_1$ is a reductive algebraic group. Fix an embedding
$\g_1\hookrightarrow \g$ such that $\t\subset \n_\g(\g_1)$. Then
$\t\cap\g_1$ is a Cartan subalgebra and $\b\cap\g_1$ is a Borel
subalgebra of $\g_1$. For fixed Borel subgroup and maximal torus in
$G_1$ we take those with the Lie algebras $\b_1,\t_1$.

{\bf Homomorphisms and representations.} All homomorphisms of
reductive algebraic Lie algebras (for instance, representations) are
assumed to be differentials of homomorphisms of the corresponding
algebraic groups.

{\bf Identification $\g\cong \g^*$}. Let $G$ be a reductive
algebraic group. There is a $G$-invariant symmetric bilinear form
$(\cdot,\cdot)$ on $G$ such that its restriction to $\t(\R)$ is
positively definite. For instance, if $V$ is a locally effective
$G$-module, then $(\xi,\eta)=\tr_V(\xi\eta)$ has the required
properties. Note that if $H$ is a reductive subgroup of $G$, then
the restriction of $(\cdot,\cdot)$ to $\h$ is nondegenerate, so one
may identify $\h$ with $\h^*$.


{\bf Parabolic subgroups and Levi subgroups}. A parabolic subgroup
of $G$ is called {\it standard} (resp., {\it antistandard}) if it
contains $B$ (resp., $B^-$). It is known that any parabolic subgroup
is $G$-conjugate to a unique standard (and antistandard) parabolic.
Standard  (as well as antistandard) parabolics are in one-to-one
correspondence with subsets of $\Pi(\g)$. Namely, one assigns to
$\Sigma\subset \Pi(\g)$ the standard  (resp., antistandard)
parabolic subgroup, whose Lie algebra is generated by $\b$ and
$\g^{-\alpha}$ for $\alpha\in \Sigma$ (resp., by $\b^-$ and
$\g^\alpha,\alpha\in \Sigma$).

By a standard Levi subgroup in $G$ we mean the Levi subgroup
containing $T$ of a standard (or an antistandard) parabolic
subgroup.

{\bf Simple Lie algebras, their roots and weights.}

Simple roots of a simple Lie algebra $\g$ are denoted by $\alpha_i$.
The numeration is described below. By $\pi_i$ we denote the
fundamental weight corresponding to $\alpha_i$.

{\it Classical algebras}. In all cases for $\b$ (resp. $\t$) we take
the algebra of all upper triangular (resp., diagonal) matrices in
$\g$.

$\g=\sl_n$. Let $e_1,\ldots,e_n$ denote the standard basis in
$\C^{n}$ and $e^1,\ldots,e^n$ the dual basis in $\C^{n*}$.  Choose
the generators $\varepsilon_i, i=\overline{1,n},$ given by
$\langle\varepsilon_i,\operatorname{diag}(x_1,\ldots,x_n)\rangle=x_i$.
Put $\alpha_i=\varepsilon_i-\varepsilon_{i+1}, i=\overline{1,n-1}$.

$\g=\so_{2n+1}$. Let $e_1,\ldots,e_{2n+1}$ be the standard basis in
$\C^{2n+1}$. We suppose  $\g$ annihilates the form
$(x,y)=\sum_{i=1}^{2n+1}x_iy_{2n+2-i}$. Define
$\varepsilon_i\in\t^*, i=\overline{1,n},$ by $\langle \varepsilon_i,
\operatorname{diag}(x_1,\ldots,x_n,0,-x_n,\ldots,-x_1)\rangle$
$=x_i$. Put $\alpha_i=\varepsilon_i-\varepsilon_{i+1},
i=\overline{1,n-1}, \alpha_n=\varepsilon_n$.

$\g=\sp_{2n}$. Let $e_1,\ldots,e_{2n}$ be the standard basis in
$\C^{2n}$. We suppose that $\g$ annihilates the form
$(x,y)=\sum_{i=1}^{n}(x_iy_{2n+1-i}-y_ix_{2n+1-i})$. Let us define
$\varepsilon_i\in\t^*, i=\overline{1,n},$ by $\langle \varepsilon_i,
\operatorname{diag}(x_1,\ldots,x_n,-x_n,\ldots,-x_1)\rangle=x_i$.
Put $\alpha_i=\varepsilon_i-\varepsilon_{i+1}, i=\overline{1,n-1},
\alpha_n=2\varepsilon_n$.

$\g=\so_{2n}$. Let $e_1,\ldots,e_{2n}$ be the standard basis in
$\C^{2n}$. We suppose that $\g$ annihilates the form
$(x,y)=\sum_{i=1}^{2n}x_iy_{2n+1-i}$. Define $\varepsilon_i\in\t^*,
i=\overline{1,n},$ in the same way as for $\g=\sp_{2n}$. Put
$\alpha_i=\varepsilon_i-\varepsilon_{i+1}, i=\overline{1,n-1},
\alpha_n=\varepsilon_{n-1}+\varepsilon_n$.

{\it Exceptional algebras}. For roots and weights of exceptional Lie
algebras we use the notation from \cite{VO}. The numeration of
simple roots is also taken from \cite{VO}.

{\bf Subalgebras in semisimple Lie algebra.} For semisimple
subalgebras of exceptional Lie algebras we use the notation from
\cite{Dynkin}. Below we explain the notation for classical algebras.

Suppose $\g=\sl_n$. By $\sl_k,\so_k,\sp_{k}$ we denote the
subalgebras of $\sl_n$ annihilating a subspace  $U\subset \C^n$
 of dimension $n-k$, leaving its complement $V$ invariant, and (for
$\so_k,\sp_k$) annihilating a nondegenerate orthogonal or symplectic
form on $V$.

The subalgebras $\so_k\subset\so_n, \sp_k\subset \sp_n$ are defined
analogously. The subalgebra $\gl_k^{diag}$ is embedded into
$\so_n,\sp_n$ via the direct sum of $\tau,\tau^*$ and a trivial
representation (here $\tau$ denotes the tautological representation
of  $\gl_k$). The subalgebras
$\sl_k^{diag},\so_k^{diag},\sp_k^{diag}\subset \so_n,\sp_n$ are
defined analogously. The subalgebra $G_2$ (resp., $\spin_7$) in
$\so_n$ is the image of $G_2$ (resp., $\so_7$) under the direct sum
of the 7-dimensional irreducible (resp., spinor) and the trivial
representations.

Finally, let $\h_1,\h_2$ be subalgebras of $\g=\sl_n,\so_n,\sp_n$
described above. While writing $\h_1\oplus\h_2$, we always mean that
$(\C^n)^{\h_1}+(\C^n)^{\h_2}=\C^n$.

The description above determines a subalgebra uniquely up to
conjugacy in $\Aut(\g)$.

Now we list some notation used in the text.

\begin{longtable}{p{4.5cm} p{10cm}}
\\$\sim_G$& the equivalence relation induced by an action of group
$G$.
\\$A^{(B)}$& the subset of all $B$-semiinvariant functions in a $G$-algebra $A$.\\
$A^\times$& the group of all invertible elements of an algebra $A$.
\\ $\Aut(\g)$& the group of
automorphisms of a Lie algebra $\g$.
\\$e_\alpha$& a nonzero element of the root subspace $\g^\alpha$.
\\ $(G,G)$& the commutant of a group $G$.
\\ $[\g,\g]$& the commutant of a Lie algebra $\g$.
\\ $G^{\circ}$& the connected component of unit of an algebraic group $G$.
\\ $G*_HV$& the homogeneous bundle over $G/H$ with the fiber $V$.
\\ $[g,v]$& the equivalence class of $(g,v)$ in $G*_HV$.
\\ $G_x$& the stabilizer of $x\in X$ under an action
$G:X$.
\\ $\g^{\alpha}$& the root subspace of $\g$ corresponding to a root  $\alpha$.
\\ $\g^{(A)}$& the subalgebra $\g$ generated by $\g^{\alpha}$ with $\alpha\in A\cup
-A$.
\\ $G^{(A)}$& the connected subgroup of $G$ with Lie algebra
$\g^{(A)}$.
\\$\Gr(V,d)$& the Grassmanian of $d$-dimensional subspaces of a
vector space $V$.
\\ $\Int(\h)$& the group of inner automorphisms of a Lie algebra
$\g$.
\\ $N_G(H)$ ($N_G(\h)$)& the normalizer of a  subgroup $H$ (subalgebra $\h\subset\g$) in
a group $G$.\\ $\n_\g(\h)$& the normalizer  of a subalgebra $\h$ in
a Lie algebra $\g$.
\\ $\Quot(A)$& the fraction field of $A$.
\\ $\rank(G)$& the rank of an algebraic group $G$.
\\ $R(\lambda)$& the irreducible representation of a reductive algebraic group (or a
reductive Lie algebra) corresponding to a highest weight $\lambda$.
\\ $\Rad_u(H)$ ($\Rad_u(\h)$)& the unipotent radical of an
algebraic group $H$  (of an algebraic Lie algebra $\h$).
\\
 $s_\alpha$& the reflection in a Euclidian space corresponding to a vector $\alpha$.
\\ $\td A$& the transcendence degree of an algebra $A$.
\\ $U^{\skewperp}$& the skew-orthogonal complement to a subspace $U\subset V$ of
a symplectic vector space $V$.
\\ $V^\g$& $=\{v\in V| \g v=0\}$, where $\g$ is a Lie algebra and
$V$ is a $\g$-module.
\\ $V(\lambda)$& the irreducible module of the highest weight $\lambda$
over a reductive algebraic group or a reductive Lie algebra.\\
$W(\g)$& the Weyl group of a reductive Lie algebra $\g$.\\
$\X(G)$& the character lattice of an algebraic group $G$.\\
$\X_{G}$& the weight lattice of a reductive algebraic group $G$.\\
$X^G$& the fixed point set for an action $G:X$.\\
$X\quo G$& the categorical quotient for an action $G:X$, where $G$
is a reductive group and $X$ is an affine  $G$-variety.\\
$\#X$& the number of elements in a set $X$.
\\   $Z_G(H)$, ($Z_G(\h)$)& the centralizer of a subgroup $H$ (of a subalgebra $\h\subset\g$)
in an algebraic group $G$.
\\ $Z(G)$&$:=Z_G(G)$.
\\ $\z_\g(\h)$& the centralizer of a subalgebra $\h$ in $\g$.
\\ $\z(\g)$& $:=\z_\g(\g)$.
\\ $\alpha^\vee$& the dual root to $\alpha$.
\\  $\Delta(\g)$& the root system of a reductive Lie algebra $\g$.
\\ $\lambda^*$& the dual highest weight to $\lambda$.
\\ $\Lambda(\g)$& the root lattice of a reductive Lie algebra $\g$.
\\ $\Pi(\g)$& the system of simple roots for a reductive Lie algebra $\g$.
\\$\pi_{G,X}$& the (categorical) quotient morphism $X\rightarrow X\quo G$.
\end{longtable}

\section{Known results and constructions}\label{SECTION_prelim}
In this section $G$ is a connected reductive group and $X$ is an
irreducible $G$-variety. We fix a Borel subgroup $B\subset G$ and a
maximal torus $T\subset B$. Put $U=\Rad_u(B)$.

\subsection{Introduction}\label{SUBSECTION_prelim_intro}
In this section we quote known results and constructions related to
Cartan spaces, Weyl groups and weight  lattices. We also prove some
more or less easy results for which it is difficult to give a
reference.


In Subsection \ref{SUBSECTION_prelim1} we present results on the
computation of Cartan spaces. In the beginning of the subsection we
establish an important notion  of a tame inclusion of a subgroup of
$G$ into a parabolic. Then we present a reduction of the computation
from the general case to the case of affine homogeneous spaces. This
reduction  belongs to Panyushev. Then we partially present results
of \cite{ranks} on the computation of the Cartan spaces for affine
homogeneous spaces. Finally, at the end of the section we study the
behavior of Cartan spaces, Weyl groups, etc. under the twisting of
the action $G:X$ by an automorphism.

The most important part of this section is Subsection
\ref{SUBSECTION_prelim4}, where we review some definitions and
results related to Weyl groups of $G$-varieties. We start with
results of F. Knop, \cite{Knop1},\cite{Knop4},\cite{Knop3},
\cite{Knop12}. Then we quote results of \cite{comb_ham} that provide
certain reductions for computing Weyl groups. These results allow to
reduce the computation of the groups $W_{G,X}$ to the case when $G$
is simple and $\rank_G(X)=\rank G$.

Till the end of the subsection  we deal with that special case. Here
we have two types of restrictions on the Weyl group. Restrictions of
the first type are valid for smooth affine varieties. They are
derived from results of \cite{fibers}. The computation in Section
\ref{SECTION_Weyl_aff} is based on these restrictions. Their main
feature is that they describe the class of conjugacy of $W_{G,X}$
and do not answer the question whether a given reflection lies in
$W_{G,X}$.

On the other hand, we have some restrictions on the form of a
reflection lying in $W_{G,X}$. They are used in Section
\ref{SECTION_Weyl_homogen}.  These restrictions are derived from the
observation that the Weyl group of an arbitrary $X$ coincides with a
Weyl group of a certain {\it wonderful variety}. This observation
follows mainly from results of Knop, \cite{Knop4}.

\subsection{Computation of Cartan spaces}\label{SUBSECTION_prelim1}
In this subsection $G$ is a connected reductive group. The
definitions of the Cartan spaces and the Weyl group of $X$ given in
Subsection \ref{SUBSECTION_Intro1} are compatible with those given
in \cite{Knop1}, \cite{Knop3}, see \cite{Knop3}, Theorem 7.4 and
Corollary 7.5. In particular,  $W_{G,X}$ is a subquotient of $W(\g)$
(i.e., there exist subgroups $\Gamma_1,\Gamma_2\subset W(\g)$ such
that $\a_{G,X}$ is $\Gamma_1$-stable, $\Gamma_2$ is the inefficiency
kernel of the action $\Gamma_1:\a_{G,X}$, and
$W_{G,X}=\Gamma_1/\Gamma_2$). The following proposition describes
functorial properties of Cartan spaces and Weyl groups.

\begin{Prop}[\cite{Knop1}, Satz 6.5]\label{Prop:1.2.1}
Let $X_1,X_2$ be irreducible $G$-varieties and
$\varphi:X_1\rightarrow X_2$ a $G$-morphism.
\begin{enumerate}
\item Suppose $\varphi$ is dominant. Then $\a_{G,X_2}\subset \a_{G,X_1}$ and
$W_{G,X_2}$ is a subquotient of $W_{G,X_1}$.
\item Suppose $\varphi$ is generically finite. Then $\a_{G,X_1}\subset \a_{G,X_2}$
and $W_{G,X_1}$ is a subquotient of $W_{G,X_2}$.
\item If $\varphi$ is dominant and generically finite (e.g. etale), then
$\a_{G,X_1}=\a_{G,X_2}$ and $W_{G,X_1}=W_{G,X_2}$.
\item Let $X$ be an irreducible $G$-variety. There is an open $G$-subvariety $X^0\subset X$
such that $\a_{G,Gx}=\a_{G,X}, W_{G,Gx}=W_{G,X}$ for any $x\in X^0$.
\end{enumerate}
\end{Prop}

\begin{Cor}\label{Cor:1.2.2}
Let $H_1\subset H_2$ be algebraic subgroups in $G$. Then
$\a_{G,G/H_2} \subset \a_{G,G/H_1}$ and $W_{G,G/H_2}$ is a
subquotient of $W_{G,G/H_1}$. If $H_1^\circ=H_2^\circ$, then
$\a_{G,G/H_2}=\a_{G,G/H_1}$ and $W_{G,G/H_2}= W_{G,G/H_1}$.
\end{Cor}

In the sequel we write $\a(\g,\h)\index{agh@$\a(\g,\h)$}$ instead of
$\a_{G,G/H}$.

\begin{Cor}\label{Cor:1.2.3}
Let $X$ be  smooth and affine, $x$ a point of $X$ with  closed
$G$-orbit. Put $H=G_x, V=T_xX/\g_*x, X'=G*_HV$. Then
$\a_{G,X}=\a_{G,X'}$, $W_{G,X}=W_{G,X'}$.
\end{Cor}
\begin{proof}
This is a direct consequence of the Luna slice theorem for smooth
points (\cite{VP}, Subsection 6.5) and assertion 3 of Proposition
\ref{Prop:1.2.1}.
\end{proof}

\begin{Cor}\label{Cor:1.2.4}
Let $H$ be a reductive subgroup of $G$ and $V$  an $H$-module. Then
\begin{enumerate}
\item $\a_{G,G*_HV}\subset \a_{G,G*_{H_0}V}$ for any normal subgroup
$H_0\subset H$. If $\a_{G,G*_HV}= \a_{G,G*_{H_0}V}$, then
$W_{G,G*_HV}\subset W_{G,G*_{H_0}V}$. If $H^\circ\subset H_0$, then
$W_{G,G*_{H}V}= W_{G,G*_{H_0}V}$.
\item $\a_{G,G*_H(V/V^H)}=\a_{G,G*_HV}$, $W_{G,G*_H(V/V^H)}=W_{G,G*_HV}$.
\end{enumerate}
\end{Cor}
\begin{proof}
Assertion 1 follows from assertions 1,3 of Proposition
\ref{Prop:1.2.1}. To prove assertion 2 note that there is an
isomorphism of $G$-varieties $G*_HV\cong G*_H(V/V^H)\times V^H$ ($G$
is assumed to act trivially on $V^H$). It remains to apply assertion
4 of Proposition \ref{Prop:1.2.1}.
\end{proof}

In the sequel we write
$\a(\g,\h,V)\index{aghv@$\a(\g,\h,V)$},W(\g,\h,V)\index{Wgh@$W(\g,\h,V)$}$
instead of $\a_{G,G*_HV},W_{G,G*_HV}$.

Now we reduce the computation of Cartan spaces for homogeneous
spaces to that for affine homogeneous vector bundles. To this end we
need one fact about subgroups in $G$ due to Weisfeller,
\cite{Weisfeller}.

\begin{Prop}\label{Prop:1.5.1}
Let $H$ be an algebraic subgroup of $G$. Then there exists a
parabolic $Q\subset G$ such that $\Rad_u(H)\subset \Rad_u(G)$.
\end{Prop}

\begin{defi}\label{Def:1.5.2}
Under the assumptions of the previous proposition, we say that the
inclusion $H\subset Q$ is {\it tame}.
\end{defi}

Algorithm \ref{Alg:9.1.1} allows one to construct a tame inclusion.

Fix Levi decompositions $H=\Rad_u(H)\leftthreetimes S,
Q=\Rad_u(Q)\leftthreetimes M$. Conjugating $H$ by an element of $Q$,
one may assume that $S\subset M$. Besides, conjugating $Q,M,H$ by an
element of $G$, one may assume that  $Q$ is an antistandard
parabolic and $M$ is its standard Levi subgroup.

The following lemma and remark seem to be standard.

\begin{Lem}\label{Lem:1.1.1}
Let $Q,M,H,S$ be such as in the previous discussion and $V$ an
$H$-module, $Q^-:=BM$. Put $X=G*_HV, \underline{X}:=Q*_HV$. Then the
fields $\C(\underline{X}),\C(X)^{\Rad_u(Q^-)}$ are $M$-equivariantly
isomorphic.
\end{Lem}
\begin{proof}
Consider the map $\iota:X^0:=\Rad_u(Q^-)\times
\underline{X}\rightarrow X, (q,[m,x])\mapsto [qm,x]$. Define the
action of $Q^-=M\rightthreetimes \Rad_u(Q^-)$ on $X^0$ as follows:
$q_1.(q,[g,v])=(qq_1,[g,v]), $ $
m_1.(q,[g,v])=(m_1qm_1^{-1},[m_1g,v]), q_1,q\in \Rad_u(Q^-), m_1\in
M, g\in Q,v\in V$. The morphism $\iota$ becomes $Q$-equivariant. One
easily checks that $\iota$ is injective. Since $\dim X^0=\dim X$,
the morphism $\iota$ is dominant. Taking into account that $X$ is
smooth, we see that $\iota$ is an open embedding. So $\C(X^0),\C(X)$
are $Q^-$-equivariantly isomorphic whence the claim of the lemma.
\end{proof}

\begin{Rem}\label{Rem:1.5.3}
Let $Q,M,H,S,V$ be such as in Lemma \ref{Lem:1.1.1}.
 Let us construct an $M$-isomorphism of  $Q*_HV$ and
$M*_S((\Rad_u(\q)/\Rad_u(\h))\oplus V)$. Consider the decreasing
chain of ideals
$\Rad_u(\q)=\q^{(1)}\supset\q^{(2)}\supset\ldots\supset\q^{(m)}=\{0\}$,
where $\q^{(i+1)}=[\Rad_u(\q),\q^{(i)}]$. Choose an
 $S$-submodule $V_i\subset \q^{(i)},i=\overline{1,m-1},$
 complementary to
$\q^{(i+1)}+(\q^{(i)}\cap\h)$. The map $M*_S(V_1\oplus\ldots
V_{m-1}\oplus V)\rightarrow Q*_HV, [m,(v_1,\ldots,v_{m-1},v)]\mapsto
[m\exp(v_1)\ldots\exp(v_{m-1}),v]$ is a well-defined
$M$-isomorphism.
\end{Rem}

The following proposition  stems from Lemma \ref{Lem:1.1.1} and
Remark \ref{Rem:1.5.3}. It is also a direct generalization of a part
of Theorem 1.2 from \cite{Panyushev2} (see also
\cite{Panyushev_thesis}, Theorem 2.5.20).

\begin{Prop}\label{Prop:1.5.4}
Let $Q$ be an antistandard parabolic subgroup of $G$, $M$ its
standard Levi subgroup and $H$ an algebraic subgroup of $Q$ such
that the inclusion $H\subset Q$ is tame and $S:=M\cap H$ is a Levi
subgroup in $H$. Then
$\X_{G,G*_HV}=\X_{M,M*_{S}((\Rad_u(\q)/\Rad_u(\h))\oplus V)}$.
\end{Prop}

Next, we reduce the case of affine homogeneous vector bundles to
that of affine homogeneous spaces. To state the main result
(Proposition \ref{Prop:1.5.8}) we need the notion of the
distinguished component.

%
%

First of all, set
\begin{equation}\label{eq:0.4:1}L_{G,X}:=Z_G(\a_{G,X}).\end{equation}
\begin{equation}\label{eq:0.4:2}L_{0\,G,X}:=\{g\in L|\chi(g)=1, \forall \chi\in \X_{G,X}
\}.\end{equation}

\begin{Prop}[\cite{comb_ham}, Proposition 8.4]\label{Prop:1.5.6}
Suppose $X$ is  smooth and quasiaffine.
\begin{enumerate}
\item
Let $L_1$ be a normal subgroup of $L_{0\,G,X}$.  Then there exists a
unique irreducible component $\underline{X}\subset X^{L_1}$ such
that $\overline{U\underline{X}}=X$.
\item
Set $P:=L_{G,X}B$. Let $S$ be a locally closed $L_{G,X}$-stable
subvariety of $X$ such that $(L_{G,X},L_{G,X})$ acts trivially on
$S$ and the map $P*_{L_{G,X}}S\rightarrow X, [p,s]\mapsto ps,$ is an
embedding (such $S$ always exists, see \cite{Knop3}, Section 2,
Lemma 3.1).  Then $\underline{X}:=\overline{FS}$, where
$F:=\Rad_u(P)^{L_1}$.
\end{enumerate}
\end{Prop}

\begin{defi}\label{Def:1.5.7}
The component $\underline{X}\subset X^{L_1}$ satisfying the
assumptions of Proposition \ref{Prop:1.5.6} is said to be {\it
distinguished}.
\end{defi}

The distinguished component for $L_1=L_{0\,G,X}$ was considered by
Panyushev in \cite{Panyushev3}.

\begin{Prop}\label{Prop:1.5.8} Let $H$ be a reductive subgroup in $G$, $V$ an $H$-module
and $\pi$ the natural projection $G*_HV\rightarrow G/H$. Put
$L_1:=L_{0\,G,G/H}^\circ$. Let $x$ be a point from the distinguished
component of $(G/H)^{L_1}$. Then
$\l_{0\,G,G*_HV}=\l_{0\,L_1,\pi^{-1}(x)}$.
\end{Prop}
\begin{proof}
Let $x_1$ be a canonical point in general position in the sense of
\cite{Panyushev_thesis}, Definition 5 of Subsection 2.1. This means
that $B_{x_1}=B\cap L_1$. It follows from Theorem 2.5.20 from
\cite{Panyushev_thesis} that $L_{0\,G,G*_HV}^\circ\cap B=
L_{0\,L_1,\pi^{-1}(x)}^\circ\cap B$. Thus
$L_{0\,G,G*_HV}^\circ=L_{0\,L_1,\pi^{-1}(x)}^\circ$. Note that the
$L_1$-module $\pi^{-1}(x)$ does not depend on  the choice of a point
$x$ from the distinguished component. Now it remains to note that
the distinguished component of  $(G/H)^{L_1}$ contains a canonical
point in general position.  Let $S$ be an $L_{G,G/H}$-subvariety in
$G/H$ mentioned in assertion 2 of Proposition \ref{Prop:1.5.6} (for
$X=G/H$). Such $S$ consists of canonical points in general position.
\end{proof}

Algorithm \ref{Alg:9.1.1} computes the subgroups $L_{0\,G,V}$ for a
$G$-module $V$, see \cite{Panyushev_thesis}.

Thus the computation of  $\a_{G,X}$ is reduced to the  computation
of the following data:
\begin{enumerate}
\item The  spaces $\a(\g,\h)$, where $H$ is a reductive subgroup of $G$.
\item A point from the distinguished component of $(G/H)^{L_1}$,
where $L_1=L_{0\,G,G/H}^\circ$, for a reductive subgroup $H\subset
G$.
\end{enumerate}

Now we are going  to present results of the paper \cite{ranks}
concerning the computation of $\a(\g,\h)$.  Until a further notice
$H$ denotes a reductive subgroup in $G$.

To state our main results we need some definitions. We begin with a
standard one.

\begin{defi}\label{Def:1.7.1}
A subalgebra $\h\subset\g$ is said to be {\it indecomposable}, if
for $(\h\cap\g_1)\oplus (\h\cap\g_2)\subsetneq\h$ any pair of ideals
$\g_1,\g_2\subset \g$ with $\g=\g_1\oplus\g_2$.
\end{defi}

Since
$\a(\g_1\oplus\g_2,\h_1\oplus\h_2)=\a(\g_1,\h_1)\oplus\a(\g_2,\h_2)$,
the computation of $\a(\g,\h)$ can be easily reduced to the case
when $\h\subset\g$ is indecomposable.

In virtue of Corollary~\ref{Cor:1.2.2}, $\a(\g,\h)\subset
\a(\g,\h_1)$ for any ideal $\h_1\subset \h$. This observation
motivates the following definition.

\begin{defi}\label{Def:1.7.2}
A reductive subalgebra   $\h\subset \g$ is called {\it
$\a$-essential} if  $\a(\g,\h)\subsetneq \a(\g,\h_1)$  for any ideal
$\h_1\subsetneq \h$.
\end{defi}

To make the presentation of our results more convenient, we
introduce one more class of subalgebras.

\begin{defi}\label{Def:1.7.3} An
$\a$-essential subalgebra $\h\subset\g$ is said to be {\it
saturated} if $\a(\g,\h)=\a(\g,\widetilde{\h})$, where
$\widetilde{\h}:=\h+\z(\n_\g(\h))$.
\end{defi}

Finally, we  note that Lemmas \ref{Lem:1.5.10},\ref{Lem:1.5.11}
below allow to perform the computation of $\a(\g,\h)$ just for one
subalgebra $\h$ in  a given class of $\Aut(\g)$-conjugacy.

The next proposition is a part of Theorem 1.3 from \cite{ranks}.

\begin{Prop}\label{Prop:1.7.4}
\begin{enumerate}
\item There is a unique ideal $\h^{ess}\subset \h\index{hess@$\h^{ess}$}$ such that
$\h^{ess}$ is an $\a$-essential subalgebra of $\g$ and
$\a(\g,\h)=\a(\g,\h^{ess})$. The ideal $\h^{ess}$ is maximal (with
respect to inclusion) among all ideals of $\h$ that are
$\a$-essential subalgebras of $\g$.
\item All semisimple indecomposable  $\a$-essential subalgebras in $\g$
up to $\Aut(\g)$-conjugacy are listed in Table~\ref{Tbl:1.7.4}.
\item All nonsemisimple saturated indecomposable $\a$-essential subalgebras in $\g$
up to  $\Aut(\g)$-conjugacy are listed in Table~\ref{Tbl:1.7.6}. In
all cases  $\h=[\h,\h]\oplus \z(\z_\g([\h,\h]))$ and $\dim
\h/[\h,\h]=1$.
\end{enumerate}
\end{Prop}

In Theorem 1.3 from \cite{ranks} all essential subalgebras are
classified and a way to compute the Cartan spaces for them is given.
Note that, as soon as this is done, assertion 1 of
Proposition~\ref{Prop:1.7.4} provides an effective method of the
determination of  $\h^{ess}$.

\setlongtables
\begin{longtable}{|c|c|c|c|}
\caption{Semisimple indecomposable  $\a$-essential subalgebras
$\h\subset\g$}\label{Tbl:1.7.4}\\\hline
N&$\g$&$\h$&$\a(\g,\h)$\\\endfirsthead\hline
N&$\g$&$\h$&$\a(\g,\h)$\\\endhead\hline 1&$\sl_n,n\geqslant
2$&$\sl_k, \frac{n+2}{2}\leqslant k\leqslant
n$&$\langle\pi_i,\pi_{n-i}; i\leqslant n-k\rangle$\\\hline
2&$\sl_n,n\geqslant 4$&$\sl_k\times\sl_{n-k}, \frac{n}{2}\leqslant
k\leqslant n-2$&$\langle\pi_i+\pi_{n-i},\pi_k,\pi_{n-k}; i<
n-k\rangle$
\\\hline 3&$\sl_{2n},n\geqslant
2$&$\sp_{2n}$&$\langle\pi_{2i}; i\leqslant n-1\rangle$\\\hline
4&$\sp_{2n},n\geqslant 2$&$\sp_{2k},\frac{n+1}{2}\leqslant
k\leqslant n$&$\langle\pi_i; i\leqslant 2(n-k)\rangle$\\\hline
5&$\sp_{2n},n\geqslant 2$&$\sp_{2k}\times\sp_{2(n-k)},
\frac{n}{2}\leqslant k< n $&$\langle\pi_{2i}; i\leqslant
n-k\rangle$\\\hline 6&$\sp_{2n},n\geqslant 4
$&$\sp_{2n-4}\times\sl_2\times\sl_2$&$\langle
\pi_2,\pi_4,\pi_1+\pi_3\rangle$\\\hline
7&$\sp_6$&$\sl_2\times\sl_2\times\sl_2$&$\langle\pi_2,\pi_1+\pi_3\rangle$\\\hline
 8&$\so_n,n\geqslant
7$&$\so_k, \frac{n+2}{2}\leqslant k\leqslant
n$&$\langle\pi_i,i\leqslant n-k\rangle$\\\hline
9&$\so_{4n},n\geqslant 2$&$\sl_{2n}$&$\langle\pi_{2i}; i\leqslant
n\rangle$\\\hline 10&$\so_{4n+2},n\geqslant 2
$&$\sl_{2n+1}$&$\langle \pi_{2i},\pi_{2n+1}; i\leqslant
n\rangle$\\\hline
11&$\so_9$&$\spin_7$&$\langle\pi_1,\pi_4\rangle$\\\hline
12&$\so_{10}$&$\spin_7$&$\langle
\pi_1,\pi_2,\pi_4,\pi_5\rangle$\\\hline
13&$\so_7$&$G_2$&$\langle\pi_3\rangle$\\\hline
14&$\so_8$&$G_2$&$\langle \pi_1,\pi_3,\pi_4\rangle$\\\hline
15&$G_2$&$A_2$&$\langle\pi_1\rangle$\\\hline
16&$F_4$&$B_4$&$\langle\pi_1\rangle$\\\hline
17&$F_4$&$D_4$&$\langle\pi_1,\pi_2\rangle$\\\hline
18&$E_6$&$F_4$&$\langle\pi_1,\pi_5\rangle$\\\hline
19&$E_6$&$D_5$&$\langle\pi_1,\pi_5,\pi_6\rangle$\\\hline
20&$E_6$&$B_4$&$\langle\pi_1,\pi_2,\pi_4,\pi_5,\pi_6\rangle$\\\hline
21&$E_6$&$A_5$&$\langle \pi_1+\pi_5,\pi_2+\pi_4,
\pi_3,\pi_6\rangle$\\\hline
22&$E_7$&$E_6$&$\langle\pi_1,\pi_2,\pi_6\rangle$\\\hline
23&$E_7$&$D_6$&$\langle \pi_2,\pi_4,\pi_5,\pi_6\rangle$\\\hline
24&$E_8$&$E_7$&$\langle \pi_1,\pi_2,\pi_3,\pi_7\rangle$\\\hline
25&$\h\times\h$&$\h$&$\langle\pi^*_i+\pi'_i;i\leqslant
\rank\h\rangle$\\\hline 26&$\sp_{2n}\times\sp_{2m},
m>n>1$&$\sp_{2n-2}\times\sl_2\times\sp_{2m-2}$&$\langle\pi_2,\pi_2',\pi_1+\pi_1'\rangle$\\\hline
27& $\sp_{2n}\times \sl_2,n>1$&$\sp_{2n-2}\times \sl_2$&$\langle
\pi_2, \pi_1+\pi_1'\rangle$\\\hline
\end{longtable}

If  $\g$ has two simple ideals (rows 25-27), then by $\pi_i$ (resp.,
$\pi_i'$) we denote fundamental weights of the first (resp., the
second) one.

\begin{longtable}{|c|c|}
\caption{Nonsemisimple saturated $\a$-essential indecomposable
subalgebras $\h\subset\g$}\label{Tbl:1.7.6}\\\hline
$(\g,[\h,\h])$&$\a(\g,\h)$\\\hline $(\sl_{n},\sl_k),
k>\frac{n}{2}$&$\{\sum_{i=1}^{n-k}(x_i\pi_i+x_{n-i}\pi_{n-i});
\sum_{i=1}^{n-k}i(x_i-x_{n-i})=0\}$\\\hline
$(\sl_n,\sl_{k}\times\sl_{n-k}),k> \frac{n}{2}$&$\langle
\pi_i+\pi_{n-i}; i\leqslant n-k\rangle$\\\hline
$(\sl_{2n+1},\sp_{2n})$&$\{\sum_{i=1}^{2n}x_i\pi_i;
\sum_{i=0}^{n-1}(n-i)x_{2i+1}-\sum_{i=1}^nix_{2i}=0\}$\\\hline
$(\so_{4n+2},\sl_{2n+1})$&$\langle \pi_{2i},\pi_{2n}+\pi_{2n+1};
i\leqslant n-1\rangle$\\\hline $(E_6,D_5)$&$\langle
\pi_1+\pi_5,\pi_6\rangle$\\\hline
\end{longtable}

\begin{Rem}\label{Rem:1.7.7}
Note that for all subalgebras $\h$ from Tables
\ref{Tbl:1.7.4},\ref{Tbl:1.7.6} except NN8 ($n=8$), 9, 25 the class
of $\Aut(\g)$-conjugacy of $\h$ coincides with the class of
$\Int(\g)$-conjugacy. In case
 $8, n=8,$ (resp., 9,25) the class of $\Aut(\g)$-conjugacy is the union of   3 (resp., 2,$\#
\Aut(\h)/\Int(\h)$) classes of $\Int(\g)$-conjugacy.
\end{Rem}

At the end of this subsection we consider the behavior of Cartan
spaces, Weyl groups etc. under the twisting of the action $G:X$ by
an automorphism.

Let $\tau\in \Aut(G)$. By $^\tau\! X$ we denote the $G$-variety
coinciding with $X$ as a variety, the action of $G$ being defined by
$(g,x)\mapsto \tau^{-1}(g)x$. The identity map is an isomorphism of
$^\sigma\!(^\tau\! X)$ and $^{\sigma\tau}\!X$ for $\sigma,\tau\in
\Aut(G)$. If $\tau$ is an inner automorphism, $\tau(g)=hgh^{-1}$,
then $x\mapsto hx:^\tau\!\! X\mapsto X$ is a $G$-isomorphism. Hence
the $G$-variety $^\tau\! X$ is determined up to isomorphism by the
image of $\tau$ in $\Aut(G)/\Int(G)$. In particular,  considering
$G$-varieties of the form $^\tau\! X$, one may assume that
$\tau(B)=B,\tau(T)=T$.

\begin{Lem}\label{Lem:1.5.10}
Let $\tau\in \Aut(G), \tau(B)=B,\tau(T)=T$. Then $\X_{G,^\tau\!
X}=\tau(\X_{G,X})$, $\a_{G,^\tau\! X}=\tau(\a_{G,X})$, $L_{G,^\tau
\!X}=\tau(L_{G,X})$, $L_{0\,G,^\tau\! X}=\tau(L_{0\,G,X})$,
$W_{G,^\tau\! X}=\tau W_{G,X}\tau^{-1}$. Further, if $L_1$ is a
normal subgroup in $L_{0\,G,X}$, then the distinguished components
in $X^{L_1}$ and $^\tau\! X^{\tau(L_1)}$ coincide.
\end{Lem}
\begin{proof}
Let $f\in \C(X)$ be a $B$-semiinvariant function of  weight $\chi$.
Then $f$ considered as an element of $\C(^\tau X)$ is
$B$-semiinvariant of weight $\tau(\chi)$. Assertions on
$\a_{G,\bullet}, \X_{G,\bullet}, L_{G,\bullet}, L_{0\,G,\bullet}$
follow immediately from this observation. $U$-orbits of the actions
$G:X, G: ^\tau\!\!\!X$ coincide whence the equalities for the
distinguished components. Now let $v$ be a central valuation of $X$.
Since $\C(X)^{(B)}=\C(^\tau\! X)^{(B)}$, we see that $v$ is a
central valuation of $^\tau\! X$. Let $\varphi_v,^\tau\!\varphi_v$
be the corresponding elements in $\a_{G,X}(\R)^*,\a_{G,^\tau\!
X}(\R)^*$. Then $\langle \varphi_v,\lambda\rangle=v(f)=\langle
^\tau\!\varphi_v,\tau(\lambda)\rangle$, where $f\in
\C(X)^{(B)}_\lambda$, whence the equality for the Weyl groups.
\end{proof}

\begin{Lem}\label{Lem:1.5.11}
Let $\tau\in \Aut(G)$.
\begin{enumerate}
\item If $H$ is an algebraic subgroup of $G$, then $^\tau\!(G/H)\cong G/\tau(H)$.
\item Let $H$ be a reductive subgroup of $G$ and $V$ be an $H$-module. Then
$^\tau\!(G*_HV)=G*_{\tau(H)}V$, where $\tau(H)$ acts on $V$ via the
isomorphism $\tau^{-1}:\tau(H)\rightarrow H$.
\end{enumerate}
\end{Lem}
\begin{proof}
Required isomorphisms are given by $gH\mapsto \tau(g)\tau(H)$ and
$[g,v]\mapsto [\tau(g),v]$.
\end{proof}

\subsection{Results about Weyl groups}\label{SUBSECTION_prelim4}


\begin{Prop}\label{Prop:1.2.9}
Let $\h$ be an algebraic subalgebra of $\g$ and $\h_0$ a subalgebra
of $\g$ lying in the closure of $G\h$ in $\Gr(\g,\dim\h)$. Then
$\h_0$ is an algebraic subalgebra of $\g$, $\a(\g,\h)=\a(\g,\h_0)$
and $W(\g,\h_0)\subset W(\g,\h)$.
\end{Prop}
\begin{proof}
This stems from \cite{Knop12}, Lemmas 3.1, 4.2.
\end{proof}

Let $T_0$ be a torus and $\pi:\widetilde{X}\rightarrow X$ be a
principal locally trivial $T_0$-bundle, where $T_0$ is a torus. In
particular, $T_0$ acts freely on $\widetilde{X}$ and $X$ is a
quotient for this action. Now let $H$ be an algebraic group acting
on $X$. The bundle $\pi:\widetilde{X}\rightarrow X$ is said to be
{\it $H$-equivariant} if $\widetilde{X}$ is equipped with an action
$H:\widetilde{X}$ commuting with that of $T_0$ and such that $\pi$
is an $H$-morphism. We can consider $\widetilde{X},X$ as $H\times
T_0$-varieties ($T_0$ acts trivially on $X$) and $\pi$ as an
$H\times T_0$-morphism.

\begin{Prop}[\cite{Knop3}, Theorem 5.1]\label{Prop:1.2.5}
Let $X$ be an irreducible $G$-variety, $T_0$ a torus and
$\pi:\widetilde{X}\rightarrow X$ a  $G$-equivariant principal
locally  trivial $T_0$-bundle. Set $\widetilde{G}:=G\times T_0$.
Then $\a_{\widetilde{G},\widetilde{X}}=\a_{G,X}\oplus \t_0$,
$\t_0,\a_{G,X}$ are $W_{\widetilde{G},\widetilde{X}}$-subspaces in
$\a_{\widetilde{G},\widetilde{X}}$ and
$W_{\widetilde{G},\widetilde{X}}$ acts trivially on $\t_0$ and as
$W_{G,X}$ on $\a_{G,X}$.
\end{Prop}

Now we want to establish a relation between the linear part of the
cone $\V_{G,X}$, which coincides with $\a_{G,X}^*(\R)^{W_{G,X}}$,
and a certain subgroup of $\Aut^G(X)$.

\begin{defi}\label{Def:1.2.13} A $G$-automorphism $\varphi$ of $X$ is said to be {\it central}
if $\varphi$ acts on $\C(X)^{(B)}_\lambda$ by the multiplication by
a constant for any $\lambda\in \X_{G,X}$. Central automorphisms of
$X$ form the subgroup of $\Aut^G(X)$ denoted by $\A_G(X)$.
\end{defi}

\begin{Lem}[\cite{Knop8}, Corollary 5.6]\label{Lem:1.2.11}
A central automorphism commutes with any $G$-automorphism of $X$.
\end{Lem}

It turns out that $\A_G(X)$ is not a birational invariant of $X$.
However, by Corollary 5.4 from \cite{Knop8}, there is an open
$G$-subvariety $X_0$ such that $\A_{G}(X_0)=\A_{G}(X_1)$ for any
open $G$-subvariety $X_1\subset X_0$. We denote $\A_{G}(X_0)$ by
$\A_{G,X}\index{agx@$\A_{G,X}$}$.

Put
$A_{G,X}:=\operatorname{Hom}(\X_{G,X},\C^\times)\index{agx@$A_{G,X}$}$.
The group $\A_{G,X}$ is embedded into $A_{G,X}$ as follows. We
assign $a_{\varphi,\lambda}$ to $\varphi\in \A_{G,X},\lambda\in
\X_{G,X}$ by the formula $\varphi
f_\lambda=a_{\varphi,\lambda}f_\lambda$, $f\in \C(X)^{(B)}_\lambda$.
The map $\iota_{G,X}:\A_{G,X}\rightarrow A_{G,X}$ is given by
$\lambda(\iota_{G,X}(\varphi))=a_{\varphi,\lambda}$. Clearly,
$\iota_{G,X}$ is a well-defined homomorphism.

\begin{Lem}[\cite{Knop8}, Theorem 5.5]\label{Lem:1.2.12} $\iota_{G,X}$ is
injective and its image is closed. In particular, $\A_{G,X}^\circ$
is a torus.
\end{Lem}

In the sequel we identify $\A_{G,X}$ with $\im\iota_{G,X}$.

The following proposition is a corollary of Satz 8.1 and (in fact,
the proof of) Satz 8.2 from \cite{Knop4}.

\begin{Prop}\label{Prop:1.2.10}
The Lie algebra of  $\A_{G,X}$ coincides with $\a_{G,X}^{W_{G,X}}$.
\end{Prop}

Until a further notice we assume that $X$ is normal. Let us study a
relation between the central valuation cones of $X$ and of certain
$G$-divisors of $X$. Our goal is to prove Corollary \ref{Cor:1.2.8},
which will be used in the end of the subsection to obtain
restrictions on possible Weyl groups of homogeneous spaces
(Proposition \ref{Prop:8.2.1}).

Let $v_0$ be a nonzero $\R$-valued discrete geometric valuation of
$\C(X)^B$  (by definition, $v_0$ is  {\it geometric} if it is a
nonnegative multiple of the valuation induced by a divisor on some
model of $\C(X)^B$). We denote by ${\mathcal V}_{v_0}$ the set of
all geometric $G$-valuations $v$ of $\C(X)$ such that
$v|_{\C(X)^B}=kv_0$ for some $k\geqslant 0$.

Now we construct a map from ${\mathcal V}_{v_0}$ to a finite
dimensional vector space similar to that from Subsection
\ref{SUBSECTION_Intro1}. Namely, we have an exact sequence of
abelian groups
$$ \{1\}\rightarrow \C(X)^{B\times}\rightarrow \C(X)^{(B)}\setminus
\{0\}\rightarrow \X_{G,X}\rightarrow \{0\}.
$$
This sequent splits because $\X_{G,X}$ is  free. Fix a splitting
$\lambda\mapsto f_{\lambda}$. We assign an element
$(\varphi_v,k_v)\in \a_{G,X}\times\R_{\geqslant 0}$  to $v\in
{\mathcal V}_{v_0}$ by the following formula:
$$\langle \varphi_v,\lambda\rangle=v(f_\lambda), v|_{\C(X)^{B}}=k_vv_0.$$

There is a statement similar to Theorem \ref{Thm:0.1.4}.
\begin{Prop}\label{Prop:1.2.6}
\begin{enumerate}
\item The map $v\mapsto (\varphi_v,k_v)$ is injective, so we may identify
${\mathcal V}_{v_0}$ with its image.
\item ${\mathcal V}_{v_0}\subset \a_{G,X}\times\R_{\geqslant 0}$
is a simplicial cone. One of its faces is the central valuation cone
considered as a subspace in $\a_{G,X}\times \{0\}\subset
\a_{G,X}\times \R_{\geqslant 0}$.
\end{enumerate}
\end{Prop}

In particular, there is $v\in \V_{v_0}$ such that $\V_{v_0}$
coincides with the cone spanned by $v$ and $\V_{G,X}$. Such $v$ is
defined uniquely up to rescaling and the shift by an element of
$\V_{G,X}\cap -\V_{G,X}$.

The first part of this proposition stems from \cite{Knop4},
Koroll\"{a}re 3.6,4.2. The second one is a reformulation of the
second part of Satz 9.2 from \cite{Knop4}.

The following result is a special case of Satz 7.5 from
\cite{Knop4}.

\begin{Prop}\label{Prop:1.2.7} Let $v$ be as above and $D$ be a prime $G$-divisor on $X$ such that
its valuation is a positive multiple of $v$. Then
$\a_{G,D}=\a_{G,X}, W_{G,D}=W_{G,X}$.
\end{Prop}

\begin{Cor}\label{Cor:1.2.8}
Let $X$ be an irreducible $G$-variety. Then there exists a spherical
$G$-variety $X'$ such that $\a_{G,X}=\a_{G,X'}, W_{G,X}=W_{G,X'}$.
\end{Cor}
\begin{proof}
Replacing $X$ with some birationally equivalent $G$-variety, we may
assume that there is a divisor $D\subset X$ satisfying the
assumptions of Proposition \ref{Prop:1.2.7}. Since the valuation
corresponding to $D$ is noncentral, we have $c_G(D)=c_G(X)-1$
(\cite{Knop4},Satz 7.3). Now it remains to use the induction on
complexity.
\end{proof}

Now we quote some results from \cite{comb_ham} providing some
reduction procedures for computing Weyl groups.

Until a further notice, $X$ is a smooth quasiaffine $G$-variety. Set
$L_0:=L_{0\,G,X}^\circ$.  Let $\underline{X}$ denote the
distinguished component of $X^{L_0^\circ}$, see Definition
\ref{Def:1.5.7}. Note that $\underline{X}_0$ is a smooth quasiaffine
variety. Its smoothness is a standard fact, see, for example,
\cite{comb_ham}, Lemma 8.6. Put
$\underline{G}:=N_G(L_0^\circ,\underline{X})/L_0^\circ$. It is a
reductive group acting on $\underline{X}$. Note that its tangent
algebra $\underline{\g}$ can be naturally identified with
$\l_0^\perp\cap \n_\g(\l_0)=\g^{L_0^\circ}\cap \z(\l_0)^\perp\subset
\g$. Further, note that $\t\subset \n_\g(\g_1)$. Thus there are the
distinguished maximal torus $\underline{T}$ and the Borel subgroup
$\underline{B}$ in $\underline{G}$, see
Section~\ref{SECTION_notation}. Let $\Gamma$ denote the image of
$N_{\underline{G}}(\underline{T})\cap
N_{\underline{G}}(\underline{B})$ in $\GL(\underline{\t})$.

\begin{Thm}[\cite{comb_ham}, Theorem 8.7, Proposition 8.1]\label{Thm:3.0.3}
 $\a_{G,X}=\a_{\underline{G}^\circ,\underline{X}}=\underline{\t}$,
$W_{G,X}=W_{\underline{G}^\circ,\underline{X}}\leftthreetimes
\Gamma$.
\end{Thm}

It is easy to prove, see Section \ref{SECTION_distinguished}, that
if $X$ is a homogeneous space (resp., affine homogeneous vector
bundle), then $\underline{X}$ is a homogeneous space (resp., affine
homogeneous vector bundle) with respect to the action of
$\underline{G}^\circ$.

\begin{Prop}\label{Cor:3.4.3}
Let $H$ be an algebraic subgroup of $G$,  $H=S\rightthreetimes
\Rad_u(H)$ its Levi decomposition, $Q$ an antistandard parabolic in
$G$  and $M$ its standard Levi subgroup. Suppose that
$\Rad_u(H)\subset \Rad_u(Q), S\subset M$. Then
\begin{enumerate} \item If $G/H$ is quasiaffine and $\a(\g,\h)=\t$, then
$\a(\m,\s,\Rad_u(\q)/\Rad_u(\h))=\t$ and $W(\g,\h)\cap
M/T=W(\m,\s,\Rad_u(\q)/\Rad_u(\h))$.
\item Suppose  $H=S$. Let $V$ be
an $H$-module such that $\a(\g,\h,V)=\t$. Then
$\a(\m,\h,\Rad_u(\q)\oplus V)=\t$ and $W(\g,\h,V)\cap M/T=
W(\m,\h,\Rad_u(\q)\oplus V)$.
\end{enumerate}
\end{Prop}
\begin{proof}
This follows from \cite{comb_ham}, Proposition 8.13, and Lemma
\ref{Lem:1.1.1}.
\end{proof}

The following proposition may be considered as a weakened version of
Proposition \ref{Prop:1.2.10}

\begin{Prop}[\cite{comb_ham}, Proposition 8.3]\label{Prop:3.4.3}
Let $T_0$ be a torus acting on  $X$ by $G$-equivariant
automorphisms. Put $\widetilde{G}=G\times T_0$. Then
$\a_{\widetilde{G},X}\cap\g\subset \a_{G,X}$ and
$\a_{G,X}\cap(\a_{\widetilde{G},X}\cap\g)^\perp\subset
\a_{G,X}^{W_{G,X}}$.
\end{Prop}

\begin{Cor}\label{Cor:3.4.6}
Let $H_1,H_2$ be subgroups of $G$ such that $H_2\subset N_G(H_1)$
and $H_2/H_1$ is a torus. Suppose that $G/H_1$ is quasiaffine. Then
$W(\g,\h_2)\cong W(\g,\h_1)$, the inclusion
$\a(\g,\h_2)\hookrightarrow \a(\g,\h_1)$ (induced by the natural
epimorphism $G/H_1\rightarrow G/H_2$) is $W(\g,\h_2)$-equivariant
and the action of $W(\g, \h_2)$ on $\a(\g,\h_2)^\perp\cap
\a(\g,\h_1)$ is trivial.
\end{Cor}
\begin{proof}
The morphism $G/H_1\rightarrow G/H_2$ is a  $G$-equivariant
principal $H_2/H_1$-bundle. Put $\widetilde{G}=G\times H_2/H_1$.
According to Proposition~\ref{Prop:1.2.5},  $W(\g,\h_2)$ is
identified with $W_{\widetilde{G},G/H_1}$ and the embedding
$\a(\g,\h_2)\hookrightarrow
\a(\g,\h_2)\oplus\t_0=\a_{\widetilde{G},G/H_1}$ is
$W(\g,\h_2)$-equivariant, where the action $W(\g,\h_2):\t_0$ is
trivial. The embedding $\a(\g,\h_2)\hookrightarrow \a(\g,\h_1)$ is
the composition of the embedding $\a(\g,\h_2)\hookrightarrow
\a_{\widetilde{G},G/H_1}$ and the orthogonal projection
$\a_{\widetilde{G},G/H_1}\rightarrow \a(\g,\h_1)$. The required
claims follow now from Proposition~\ref{Prop:3.4.3} applied to the
action $\widetilde{G}:G/H_1$.
\end{proof}

The next well-known statement describes the behavior of the Weyl
group under a so called  parabolic induction. The only case we need
is that of homogeneous spaces of rank equal to $\rank G$. We give
the proof only to illustrate our techniques.

\begin{Cor}\label{Cor:3.4.7}
Let $Q$ be an antistandard parabolic subgroup of $G$ and $M$ the
standard Levi subgroup of $Q$. Further, let $\h$ be a subalgebra of
$\g$ such that $\Rad_u(\q)\subset\h,\a(\g,\h)=\t$. Finally, assume
that $G/H$ is quasiaffine. Then $\a(\m,\m/\h)=\t$,
$W(\g,\h)=W(\m,\h/\Rad_u(\q))$.
\end{Cor}
\begin{proof}
 Using Proposition~\ref{Cor:3.4.3}, we obtain $\a(\m,\m\cap
\h)=\t, W(\m,\h/\Rad_u(\q))\subset W(\g,\h)$. On the other hand,
$Z(M)^\circ$ acts on $G/H$ by $G$-automorphisms. By
Proposition~\ref{Prop:1.5.4},  $\a_{G\times
Z(M)^\circ,G/\Rad_u(Q)}\cap\g=\t\cap[\m,\m]$. Therefore $\a_{G\times
Z(M)^\circ,G/H}\cap\g\subset \t\cap[\m,\m]$. The inclusion
$W(\m,\h/\Rad_u(\q))\supset W(\g,\h)$ follows from
Proposition~\ref{Prop:3.4.3}.
\end{proof}

Now let $G_1,\ldots, G_k$ be all simple normal subgroups in $G$ so
that $G=Z(G)^\circ G_1\ldots G_k$. Put $T_i=T\cap G_i$. This is a
maximal torus in $G_i$.

\begin{Prop}\label{Prop:3.4.5}
\begin{enumerate}
\item Let $H$ be an algebraic subgroup of $G$ such that $\a(\g,\h)=\t$
and $G/H$ is quasiaffine. Then $G_i/(G_i\cap H)$ is quasiaffine too
and $\a(\g_i,\g_i\cap \h)=\t_i, W(\g,\h)=\prod_{i}W(\g_i,\g_i\cap
\h)$.
\item Let $H$ be a reductive subgroup of $G$ and $V$ an
$H$-module. Suppose that $\a(\g,\h,V)=\t$. Then
$\a(\g_i,\g_i\cap\h,V)=\t_i$ and $W(\g,\h,V)=\prod_i W(\g_i,\g_i\cap
\h,V)$.
\end{enumerate}
\end{Prop}
\begin{proof}
To prove the first assertion we note that the stabilizer of any
point of $G/H$ in $G_i$ is conjugate to $G_i\cap H$ and use
assertion 4 of Proposition~\ref{Prop:1.2.1} and \cite{comb_ham},
Proposition 8.8.

Proceed to assertion 2. All  $G_i$-orbits in  $G/H$ are of the same
dimension whence closed. Choose $x\in G/H$ with $(G_i)_x=G_i\cap H$.
The $G_x$-module $T_xX/\g^{(i)}_*x$ is isomorphic to $V\oplus V_0$,
where $V_0$ denotes a trivial $G_x$-module. By
Lemma~\ref{Cor:1.2.3}, $W_{G_i,G*_HV}=W(\g_i,\g_i\cap \h,V)$. It
remains to use Corollary~\ref{Cor:1.2.4}.
\end{proof}

Till the end of the section $G$ is simple, $X$ is a quasiaffine
irreducible smooth $G$-variety such that $\rank_G(X)=\rank G$. In
this case $W_{G,X}\subset W(\g)$.

At first, we consider the case when $X$ is affine. Here we present
some results on $W_{G,X}$  obtained in \cite{fibers}. Those results
can be applied here because  $W_{G,X}$ coincides up to conjugacy
with the Weyl group of the Hamiltonian $G$-variety $T^*X$ (see
\cite{Knop3}).

\begin{defi}\label{Def:1.2.14}
A subgroup $\Gamma\subset W(\g)$ is said to be {\it large} if for
any roots $\alpha,\beta\in \Delta(\g) $ such that $\beta\neq
\pm\alpha, (\alpha,\beta)\neq 0$ there exists  $\gamma\in
\R\alpha+\R\beta$ with $s_\gamma\in \Gamma$.
\end{defi}

For a subgroup $\Gamma\subset W(\g)$ we denote by $\Delta_\Gamma$
the subset of $\Delta(\g)$ consisting of all $\alpha$ with
$s_\alpha\in \Gamma$.

\begin{Prop}[\cite{fibers}, Corollaries
4.16,4.19]\label{Prop:1.2.11}
\begin{enumerate}
\item The subgroup $W_{G,X}\subset W(\g)$ is large.
\item Suppose $\g$ is a simple classical Lie algebra. Then $\Gamma\subset
W(\g)$ is large iff $\Delta_\Gamma$ is listed in
Table~\ref{Tbl:5.2.8}.
\end{enumerate}
\end{Prop}

 \begin{longtable}{|c|l|}\caption{Subsets $\Delta_\Gamma$
 for large subgroups
 $\Gamma\subset W(\g)$  when $\g$ is classical}\label{Tbl:5.2.8}\\\hline $\g$&$\Delta_\Gamma$\\\hline
$A_l,l\geqslant 2$& $\{\varepsilon_i-\varepsilon_j|i\neq j, i,j\in
I\text{ or }i,j\not\in I\}, I\subsetneq \{1,\ldots,n+1\},
I\neq\varnothing$\\\hline $B_l,l\geqslant 3$&(a) $\{\pm
\varepsilon_i\pm\varepsilon_j|i\neq j, i,j\in I \text{ or
}i,j\not\in I\}\cup \{\pm\varepsilon_i| i\in I\},
I\subsetneq\{1,\ldots,n\}$\\&(b) $\{\pm
\varepsilon_i\pm\varepsilon_j|i\neq j, i,j\in I \text{ or
}i,j\not\in I\}\cup \{\pm\varepsilon_i| i\in \{1,2,\ldots,n\}\},
I\subsetneq\{1,\ldots,n\},I\neq\varnothing$\\&(c)
$\{\varepsilon_i-\varepsilon_j|i\neq j, i,j\in I \text{ or }
i,j\not\in I\}\cup \{\pm(\varepsilon_i+\varepsilon_j), i\in
I,j\not\in I)\}, I\subset \{1,\ldots,n\}$
\\\hline $C_l,l\geqslant 2$&(a) $\{\pm \varepsilon_i\pm\varepsilon_j|i\neq j, i,j\in
I \text{ or }i,j\not\in I\}\cup \{\pm 2\varepsilon_i| i\in I\},
I\subsetneq\{1,\ldots,n\}$\\&(b) $\{\pm
\varepsilon_i\pm\varepsilon_j|i\neq j, i,j\in I \text{ or
}i,j\not\in I\}\cup \{\pm 2\varepsilon_i| i\in \{1,2,\ldots,n\}\},
I\subsetneq\{1,\ldots,n\}, I\neq\varnothing$\\&(c)
$\{\varepsilon_i-\varepsilon_j|i\neq j, i,j\in I \text{ or }
i,j\not\in I\}\cup \{\pm(\varepsilon_i+\varepsilon_j), i\in
I,j\not\in I)\}, I\subset \{1,\ldots,n\}$ \\\hline $D_l,l\geqslant
3$& (a) $\{\pm \varepsilon_i\pm\varepsilon_j|i\neq j, i,j\in I
\text{ or }i,j\not\in I\}, I\neq\{1,\ldots,n\},\varnothing$\\&(b)
$\{\varepsilon_i-\varepsilon_j|i\neq j, i,j\in I \text{ or }
i,j\not\in I\}\cup \{\pm(\varepsilon_i+\varepsilon_j), i\in
I,j\not\in I)\}, I\subset \{1,\ldots,n\}$
\\\hline
\end{longtable}

Note that some subsets $\Delta_\Gamma$ appear in
Table~\ref{Tbl:5.2.8} more than once.

Now we obtain a certain restriction on $W_{G,X}$ in terms of the
action $G:T^*X$.  To do this we introduce the notions of a {\it
$\g$-stratum} and a {\it completely perpendicular} subset of
$\Delta(\g)$.

\begin{defi}\label{Def:1.4.1}
A pair $(\h,V)$, where $\h$ is a reductive subalgebra of $\g$ and
$V$ is an $\h$-module, is said to be a $\g$-stratum. Two $\g$-strata
$(\h_1,V_1)$, $(\h_2,V_2)$ are called  {\it equivalent} if there
exist $g\in G$ and a linear isomorphism
$\varphi:V_1/V_1^{\h_1}\rightarrow V_2/V_2^{\h_2}$ such that
$\Ad(g)\h_1=\h_2$ and $ (Ad(g)\xi)\varphi(v_1)=\varphi(\xi v_1)$ for
all $\xi\in\h_1, v_1\in V_1/V_1^{\h_1}$.
\end{defi}

\begin{defi}\label{Def:1.4.2} Let $Y$ be a smooth affine variety and
$y\in Y$ a point with closed $G$-orbit. The pair $(\g_y,
T_yY/\g_*y)$ is called the $\g$-stratum of $y$. We say that $(\h,V)$
is a $\g$-stratum of $Y$ if $(\h,V)$ is equivalent to a $\g$-stratum
of a point of $Y$. In this case we write $(\h,V)\rightsquigarrow_\g
Y$.
\end{defi}

\begin{Rem}\label{Rem:1.4.3} Let us justify the terminology. Pairs $(\h,V)$
do define some stratification of $Y\quo G$ by varieties with
quotient singularities. Besides, analogous objects were called
"strata" in \cite{Schwarz3}, where the term is borrowed from.
\end{Rem}

\begin{defi}\label{defi:5.2.2}
A subset $A\subset \Delta(\g) $ is called  {\it completely
perpendicular} if the following two conditions take place:
\begin{enumerate}
\item $(\alpha,\beta)=0$ for any different $\alpha,\beta\in A$.
\item $\Span_{\R}(A)\cap \Delta(\g) =A\cup -A$.
\end{enumerate}
\end{defi}

For example, any one-element subset of $\Delta(\g) $ is completely
perpendicular.

Let $A$ be a nonempty completely perpendicular subset of $\Delta(\g)
$. By $S^{(A)}$ we denote the $\g$-stratum $(\g^{(A)},
\sum_{\alpha\in A}V^{\alpha})$, where $V^\alpha$ is, by definition,
the direct sum of two copies of the two-dimensional irreducible
$\g^{(A)}/\g^{(A\setminus\{\alpha\})}$-module.

\begin{Prop}[\cite{fibers}, Corollary 4.14]\label{Cor:5.2.3}
 If $W_{G,X}\cap W(\g^{(A)})=\{1\}$, then
$S^{(A)}\rightsquigarrow_\g X$.
\end{Prop}

In particular, if $W(\g)$ contains all reflections $W(\g)$-conjugate
to $s_\alpha,\alpha\in \Delta(\g),$ then
$S^{(\alpha)}\not\rightsquigarrow T^*X$.

Till the end of the subsection $G$ is simple, $X$ is a homogeneous
$G$-space of rank $\rank(G)$ (not necessarily quasiaffine). In this
case any element of $\Pi_{G,X}$ is a positive multiple of a unique
positive root from $\Delta(\g)$. The set of all positive roots
arising in this way is denoted by
$\widehat{\Pi}_{G,X}\index{zzzq@$\widehat{\Pi}_{G,X}$}$. Note that
$\widehat{\Pi}_{G,G/H}$ depends only on $(\g,\h)$. Therefore in the
sequel we write
$\widehat{\Pi}(\g,\h)\index{zzzq@$\widehat{\Pi}(\g,\h)$}$ instead of
$\widehat{\Pi}_{G,G/H}$.

\begin{Prop}\label{Prop:8.2.1}
Let $G,X$ be such as above and  $G\neq G_2$. Then the pair
$(\Supp(\alpha),\alpha)$, where $\alpha\in\widehat{\Pi}(\g,\h)$ is
considered as an element of the root system associated with
$\Supp(\alpha)$, are given in the following list: \begin{enumerate}
\item $(A_1,\alpha_1)$.
\item $(A_2, \alpha_1+\alpha_2)$.
\item $(B_2, \alpha_1+\alpha_2)$.
\end{enumerate}
\end{Prop}
\begin{proof}[Proof of Proposition~\ref{Prop:8.2.1}]
For an arbitrary $G$-variety $X$ let $P_{G,X}$ denote the
intersection of the stabilizers of all $B$-stable prime divisors of
$X$.

Thanks to Corollary~\ref{Cor:1.2.8}, we may assume that $X$ is
spherical. There is a subgroup $\widetilde{H}\subset G$ containing
$H$ such that $\Pi_{G,G/H}=\Pi_{G,G/\widetilde{H}}$ and
$\X_{G,G/\widetilde{H}}=\Span_\Z(\Pi_{G,G/H})$. Furthermore, the
homogeneous space $G/\widetilde{H}$ possesses a so called wonderful
embedding $G/\widetilde{H}\hookrightarrow\overline{X}$ (these two
facts follow from results of \cite{Knop8}, Sections 6,7; for
definitions and results concerning wonderful varieties see
\cite{Luna5} or \cite{Timashev_rev}, Section 30).  Since
$\a_{G,G/H}=\t$, we see that $P_{G,X}=B$. So a unique closed
$G$-orbit on $X$ is isomorphic to $G/B$. For $\alpha\in \Pi_{G,X}$
let $X_\alpha$ denote the wonderful subvariety of $X$ of rank 1
corresponding to $\alpha$. Since the closed $G$-orbit in $X_\alpha$
is isomorphic to $G/B$, we see that  $P_{G,X_\alpha}=B$. Now the
proposition stems from the classification of wonderful varieties of
rank 1 and the computation of their spherical roots. These results
are gathered in \cite{Wasserman}, Table 1.
\end{proof}

\section{Determination of distinguished components}\label{SECTION_distinguished}
\subsection{Introduction}\label{SUBSECTION_distinguished_intro}
In this section we find an algorithm to determine the distinguished
component of $X^{L_{0\,G,X}^\circ}$ in the case when $X$ is a
homogeneous space or an affine homogeneous vector bundle. This will
complete the algorithm computing $\a_{G,X}$ for the indicated
classes of varieties, see Subsection~\ref{SUBSECTION_prelim1}.
Besides, this makes possible  to use Theorem~\ref{Thm:3.0.3}.

Put $L_0=L_{0\,G,X}^\circ$. By $\underline{X}$ we denote the
distinguished component of $X^{L_0}$.

Our first task is to describe the structure of distinguished
components.

\begin{Prop}\label{Prop:6.0.2}
Here $L_0$ denotes one of the groups $L_{0,G,X},L_{0,G,X}^\circ$ and
$\underline{X}$ is the distinguished component  $X^{L_0}$.
\begin{enumerate}
\item Let  $X=G/H$ be a quasiaffine homogeneous space.  Then the
action $N_G(L_0)^\circ:\underline{X}$ is transitive. If $eH\in
\underline{X}$, then $N_G(L_0,\underline{X})=N_G(L_0)^\circ
N_H(L_0)$.
\item Let $H$ be a reductive subgroup of $G$, $V$ an $H$-module,  $X=G*_HV$.
The distinguished component  $Y$ of $(G/H)^{L_{0\,G,G/H}^\circ}$ is
contained in $\underline{X}$. If $eH\in Y$, then
$\underline{X}=N_G(L_0)^\circ*_{H\cap N_G(L)^\circ}V^{L_0}$,
$N_G(L_0,\underline{X})=N_G(L_0)^\circ N_H(L_0)$.
\end{enumerate}
\end{Prop}

So in the both cases of interest the distinguished component is
recovered from an arbitrary point of the appropriate distinguished
component in a homogeneous space. The next problem is to reduce the
determination of distinguished components for homogeneous spaces to
the case of affine homogeneous vector bundles. Let $H\subset G$ be
an algebraic subgroup. Recall that we may assume that $H$ is
contained in an antistandard parabolic subgroup  $Q\subset G$ and
for the standard Levi subgroup $M\subset Q$ the intersection $M\cap
H$ is a Levi subgroup in $H$. By Remark \ref{Rem:1.5.3}, $Q/H$ is an
affine homogeneous vector bundle over $M/(M\cap H)$.

\begin{Prop}\label{Prop:6.0.3}
Let $H,Q,M$ be such as above, $X=G/H$ be quasiaffine. The
distinguished component of the $M$-variety
$(Q/H)^{L_{0\,M,Q/H}^\circ}$ is contained in $\underline{X}$.
\end{Prop}

Proposition \ref{Prop:6.0.3} together with assertion 2 of
Proposition~\ref{Prop:6.0.2} allows to reduce  the determination of
distinguished components to the case of affine homogeneous spaces.

To state the result concerning affine homogeneous space we need some
notation. Let $\h$ be an $\a$-essential subalgebra of $\g$ (see
Definition~\ref{Def:1.7.2}), $\widetilde{\h}:=\h+\z_\g(\n_\g(\h))$.
Put $\h^{sat}=\widetilde{\h}^{ess}$. By properties of the mapping
$\h\mapsto \h^{ess}$, see \cite{ranks}, Corollary 2.8,
$\h\subset\h^{sat}$. It follows directly from the construction that
$\h^{sat}$ is a saturated subalgebra of $\g$
(Definition~\ref{Def:1.7.3}).

The next proposition allows one to find a point from a distinguished
component of an affine homogeneous space.

\begin{Prop}\label{Prop:6.0.1}
Let $H$ be a reductive subgroup of $G$ and  $X=G/H$.
\begin{enumerate}
\item Let $H^{ess}$ denote the connected subgroup of   $G$ corresponding to
$\h^{ess}$,  $\pi$  the projection $G/H^{ess}\rightarrow G/H,
gH^{ess}\mapsto gH,$ and $\underline{X}'$  the distinguished
component of $(G/H^{ess})^{L_{0\,G,G/H^{ess}}^\circ}$. Then
$\underline{X}\supset\pi(\underline{X}')$.
\item Suppose  $\h$ is $\a$-essential and $H$ is connected. Let $H^{sat}$ denote the connected
subgroup of $G$ with   Lie algebra  $\h^{sat}\subset \g$,  $\pi$ be
the projection $G/H\rightarrow G/H^{sat}, gH\mapsto gH^{sat},$ and
$\underline{X}'$ the distinguished component of
$(G/H^{sat})^{L_{0\,G,G/H^{sat}}^\circ}$. Then
$\pi^{-1}(\underline{X}')\subset \underline{X}$.
\item Let $(\g,\h)=(\g_1,\h_1)\oplus (\g_2,\h_2)$,
 $G_1,G_2,H_1,H_2$  denote the connected subgroups of  $G$
 corresponding to
 $\g_1,\g_2,\h_1,\h_2$, respectively,
$\pi$ denote the projection $G_1/H_1\times G_2/H_2\rightarrow G/H,
(g_1H_1,g_2H_2)\mapsto g_1g_2H,$ and $\underline{X}_i,i=1,2,$ denote
the distinguished component of
$(G_i/H_i)^{L_{0\,G_i,G_i/H_i}^\circ}$. Then
$\underline{X}=\pi(\underline{X}_{1}\times \underline{X}_{2})$.
\item  Suppose $\h$ is one of the subalgebras listed in Tables~\ref{Tbl:1.7.4},\ref{Tbl:1.7.6},
and $H$ is connected. Then $\underline{X}\subset X$ is stable under
the action of $N_G(L_0)$.
\item Let $\h,H$ be such as indicated in Subsection \ref{SUBSECTION_disting_compI2}.
 Then $eH\in \underline{X}$.
\end{enumerate}
\end{Prop}

The first three assertions reduce the problem of finding a point
from $(G/H)^{L_0}$ to the case when $H$ is connected and $\h$ is
$\a$-essential, saturated and indecomposable. All such subalgebras
$\h$ are listed in Tables~\ref{Tbl:1.7.4},\ref{Tbl:1.7.6}, see
Proposition~\ref{Prop:1.7.4}. Assertion 5 solves the problem in this
case. Assertion 4 is auxiliary.

Let us give a brief description of the section. In
Subsection~\ref{SUBSECTION_disting_compI1} we prove
Propositions~\ref{Prop:6.0.2}, \ref{Prop:6.0.3}. In
Subsection~\ref{SUBSECTION_disting_compI4} we find some necessary
condition (Proposition \ref{Prop:6.1.2}) for a component of
$X^{L_0}$, where $X$ is smooth and affine, to be distinguished.
Subsection~\ref{SUBSECTION_disting_compI2} describes some embeddings
$\h\hookrightarrow\g$, where $\h$ is presented in
Tables~\ref{Tbl:1.7.4},\ref{Tbl:1.7.6}. For these embeddings the
point $eH$ lies in the component of $(G/H)^{L_{0\,G,G/H}^\circ}$,
satisfying the necessary condition of Proposition~\ref{Prop:6.1.2}.
In Subsection~\ref{SUBSECTION_disting_compI3} we complete the proof
of Proposition~\ref{Prop:6.0.1}.

\subsection{Proofs of Propositions~\ref{Prop:6.0.2},\ref{Prop:6.0.3}}\label{SUBSECTION_disting_compI1}
\begin{proof}[Proof of Proposition~\ref{Prop:6.0.2}]
Suppose  $F$ is a reductive subgroup of $G$.  Then for any subgroup
$H\subset G$ the variety $(G/H)^F$ is a finite union of
$N_G(F)^\circ$-orbits. Indeed, let $x\in (G/H)^F$. Then for any
$x\in (G/H)^F$
$$T_x\left((G/H)^F\right)\subset
(T_x(G/H))^{F}=(\g/\g_x)^{F}=\g^{F}/(\g_x)^{F}=T_x(N_G(F)^\circ x
).$$ Therefore $T_x\left((G/H)^F\right)=T_x(N_G(F)^\circ x)$. To
prove the claim we remark that $N_G(F)^\circ x$ is a smooth
irreducible subvariety of $(G/H)^F$. Thence in assertion 1
$\underline{X}$  is a homogeneous $N_G(L_0)^\circ$-space. The claim
describing $N_G(L_0,\underline{X})$ is obvious.

Proceed to the proof of assertion 2. Clearly, the projection
$X^{L_0}\rightarrow (G/H)^{L_0}$ is surjective, its fiber over $x\in
(G/H)^{L_0}$ coincides with  $\pi^{-1}(x)^{L_0}$. Here $\pi$ denotes
the natural projection $X\rightarrow G/H$. It follows that, as an
$N_G(L_0,\underline{X})$-variety, $\underline{X}$ is a homogeneous
vector bundle over some component of $(G/H)^{L_0}$. Note also that
$N_G(L_0,\underline{X})=N_G(L_0,\pi(\underline{X}))$. The last
equality yields the claim on the structure of
$N_G(L_0,\underline{X})$.

It remains to verify that the distinguished component $Y$ of
$(G/H)^{L_1}$, where $L_1=L_{0\,G,G/H}^\circ$, is contained in
$\underline{X}$. This claim will follow if we check that
$U(\pi^{-1}(Y)^{L_0})$ is dense in $X$. Thanks to
Proposition~\ref{Prop:1.5.8},  $(U\cap L_1)\pi^{-1}(x)^{L_0}$ is
dense in $\pi^{-1}(x)$ for all $x\in Y$. It remains to recall that
$\overline{UY}=G/H$.
\end{proof}

\begin{proof}[Proof of Proposition~\ref{Prop:6.0.3}]
By Proposition~\ref{Prop:1.5.4}, $L_0=L_{0\,M,Q/H}^\circ$. Let $Y$
denote the distinguished component in $(Q/H)^{L_0}$. By definition,
$(U\cap M)Y$ is dense $Q/H$. To complete the proof note that
$U(Q/H)$ is dense in $X$.
\end{proof}

\subsection{An auxiliary result}\label{SUBSECTION_disting_compI4}

Here we study components of the variety $X^{L_0}$ in the case when
$X$ is smooth and affine.

\begin{defi}\label{Def:6.1.3}
Let  $H$ be an algebraic group acting on an irreducible variety $Y$.
A subgroup $H_0\subset H$ is said to be the {\it stabilizer in
general position} (shortly, s.g.p.) for the action $H:Y$ if there is
an open subset $Y^0\subset Y$ such that $H_y\sim_H H_0$ for all
$y\in Y^0$. By the {\it stable subalgebra in general position}
(shortly, s.s.g.p.) we mean  the Lie algebra of the s.g.p.
\end{defi}

\begin{Prop}\label{Prop:6.1.2}
Let $X$ be a smooth affine  $G$-variety,  $L_0$ one of the groups
$L_{0\,G,X},L_{0\,G,X}^\circ$, and $\underline{X}$ the distinguished
component of $X^{L_0}$. Then $\dim X^{L_0}= \dim X-\frac{\dim G-\dim
N_G(L_0)}{2}$. A component of $X^{L_0}$ is $N_G(L_0)$-conjugate to
$\underline{X}$ iff its dimension is equal to $\dim X^{L_0}$.
\end{Prop}
\begin{proof}
Proposition~\ref{Prop:1.5.6} implies (in the notation of that
proposition) $\dim \underline{X}=\dim X-\dim \Rad_u(P)+\dim F$.
Since $F$ is a maximal unipotent subgroup in $N_G(L_0)/L_0$, the
equality $\dim \underline{X}=\dim X-(\dim G-\dim N_G(L_0))/2$ holds.

 By \cite{Knop1}, Korollar 8.2,   $L_{0\,G,X}$ is the s.g.p. for
the action $G:T^*X$. Let us show that the action $G:T^*X$ is stable,
that is, an orbit in general position is closed. Choose a point
$x\in X$ with closed $G$-orbit. Put $H=G_x$, $V=T_xX/\g_*x$. Then
$T_x(T^*X)/\g_*x\cong V\oplus V^*\oplus \h^\perp$. Using the Luna
slice theorem, we see that the action $G:T^*X$ is stable iff so is
the representation $H:\h^\perp\oplus V\oplus V^*$. Since
$\h^\perp\oplus V\oplus V^*$ is an orthogonal $H$-module, we are
done by results of \cite{Lun_ort}.

Thanks to results of \cite{LR},   $N_G(L_0)$ permutes transitively
irreducible (=connected) components of $(T^*X)^{L_0}$ whose image in
$(T^*X)\quo G$ is dense. A component of $(T^*X)^{L_0}$ has dense
image in $(T^*X)\quo G$ if it contains a point $x$ with closed
$G$-orbit and  $G_x=L_0$ (when $L_0=L_{0\,G,X}$) or
$(G_x)^\circ=L_0$ (when $L_0=L_{0\,G,X}^\circ$).  Let $Z$ be such a
component and $G_0=N_{G}(L_0,Z)$. Then the morphism
$G*_{G_0}Z\rightarrow T^*X$ is birational. Indeed, this morphism is
dominant because a point in general position in $T^*X$ is
$G$-conjugate to a point of $Z$. On the other hand, for $z\in Z$ in
general position the equality $G_z=L_0$ (for $L_0=L_{0\,G,X}$) or
$(G_z)^\circ=L_0$ (for $L_0=L_{0\,G,X}^\circ$) holds. Therefore the
inclusion $gz\in Z$ implies $g\in G_0$.  We conclude that $\dim
Z=2\dim \underline{X}$. Note that the unit component of the s.g.p.
for the action $N_G(L_0)^\circ:Z$ coincides with $L_0^\circ$.

Denote by $\pi$ the projection $T^*X\rightarrow X$. Clearly,
$\pi((T^*X)^{L_0})= X^{L_0}$. By Corollary in \cite{VP}, Subsection
6.5, $X^{L_0}$ is smooth. Let $Z_0$ be a component of $X^{L_0}$.
Then $Z:=\pi^{-1}(Z_0)^{L_0}$ is an irreducible component of
$(T^*X)^{L_0}$.

It remains to prove that $\dim Z_0\leqslant \dim X- (\dim G-\dim
N_G(L_0))/2$ and that if the equality holds, then $Z$ contains a
point $z$ with closed $G$-orbit and $G_z=L_0$ (or $G_z^\circ=L_0$).

Note that the restriction of $\pi$ to $Z$ is the cotangent bundle
$T^*Z_0\rightarrow Z_0$. Therefore the action $N_G(L_0)^\circ:Z$ is
stable. Let $L_Z$ denote the s.g.p. for this action. Obviously,
$L_0\subset L_Z$. One also has the equality $\dim Z\quo
N_G(L_0)^\circ=\dim Z-\dim N_G(L_0)+\dim L_Z$.

The morphism $Z\quo N_G(L_0)^\circ\rightarrow T^*X\quo G$ induced by
the embedding $Z\hookrightarrow T^*X$ is finite (see~\cite{VP},
Theorem 6.16). It follows that $\dim Z-\dim N_G(L_0)+\dim
L_Z\leqslant \dim T^*X-\dim G+\dim L_0,$ or, equivalently, $$\dim
Z_0\leqslant \dim X-\frac{\dim G-\dim N_G(L_0)}{2}-\frac{\dim
L_Z-\dim L_0}{2}.$$ This completes the proof.
\end{proof}

\subsection{Embeddings of subalgebras}\label{SUBSECTION_disting_compI2}
In this subsection we construct embeddings  $\h\hookrightarrow \g$
for the pairs $(\g,\h)$ from Tables~\ref{Tbl:1.7.4},\ref{Tbl:1.7.6}.
In the next subsection we will see that the corresponding points
$eH\in G/H$ lie in the distinguished components of $(G/H)^{L_0},
L_0=L_{0\,G,G/H}^\circ$. If $eH$ lies in this distinguished
component, then, obviously, $\l_0\subset \h$ and, by
Proposition~\ref{Prop:6.1.2},

\begin{equation}\label{eq:6.2:1}
2(\dim\h-\dim\n_\h(\l_0))=\dim\g-\dim\n_\g(\l_0).
\end{equation}

\begin{Prop}\label{Prop:6.2.1}
Suppose that one of the following two assumptions  hold
\begin{enumerate}
\item  $(\g,\h)$ is one of the pairs from Table \ref{Tbl:1.7.4}.
\item  $\h=[\h,\h]\oplus \z(\z_\g([\h,\h]))$, and
$(\g,[\h,\h])$ is listed in Table \ref{Tbl:1.7.6}.
\end{enumerate}
If $[\h,\h]\hookrightarrow \g$ is embedded into $\g$ as described
below in this subsection, then $\l_0\subset \h$ and (\ref{eq:6.2:1})
holds.
\end{Prop}

We check this assertion case by case.

{\it The case  $\g=\sl_n$.}

1) The subalgebra $[\h,\h]=\sl_{k}, \frac{n}{2}\leqslant k<n,$ (the
first rows of Tables~\ref{Tbl:1.7.4},\ref{Tbl:1.7.6}) is embedded
into $\g$ as the annihilator of
$e_1,\ldots,e_{n-k},e^{1},\ldots,e^{n-k}$. In the both cases in
interest ($\dim\z(\h)=0$ or $1$) we get
$\n_\g(\l_0)=4(n-k)^2+(2k-n)^2-1,\dim\n_\h(\l_0)=(n-k)^2+(2k-n)^2-1+\dim\z(\h).$

2) The subalgebra $[\h,\h]=\sl_k\times \sl_{n-k},
\frac{n}{2}\leqslant k<n-1$ (the second rows of the both tables) is
embedded into $\g$ as $[\n_\g(\f),\n_\g(\f)]$, where $\f\cong \sl_k$
is embedded as described in 1). We have
$\dim\n_\g(\l_0)=4(n-k)+(2k-n)^2-1,
\dim\n_\h(\l_0)=2(n-k)+(2k-n)^2-1+\dim\z(\h).$

3) Let $n=2m$. The subalgebra  $\h=\sp_{2m}$ (row N3 of
Table~\ref{Tbl:1.7.4}) is embedded into   $\sl_{n}$ as the
annihilator of $e^1\wedge e^2+e^3\wedge e^4+\ldots+e^{2m-1}\wedge
e^{2m}$. We have $\dim\n_\g(\l_0)=4m-1, \dim\n_\h(\l_0)=3m$.

4) Let $n=2m+1$. The subalgebra $[\h,\h]=\sp_{2m}$ (the third row of
Table~\ref{Tbl:1.7.6}) is embedded into $\g=\sl_{2m+1}$ as the
annihilator of $e_{1},e^1$ and the  2-form $e^2\wedge
e^3+\ldots+e^{2m}\wedge e^{2m+1}$. In this case $\dim
\n_\g(\l_0)=2m^2+2m, \dim\n_\h(\l_0)=m^2$.

{\it The case $\g=\sp_{2n}$.}

1) The subalgebra $\sp_{2k}, k\geqslant \frac{n}{2}$ (row N$4$ of
Table~\ref{Tbl:1.7.4}) is embedded into
 $\sp_{2n}$ as the annihilator of $2(n-k)$ vectors
$e_1+e_{2n-1},e_2-e_{2n},e_3+e_{2n-3},e_4-e_{2n-2},\ldots,e_{2(n-k)}-
e_{2(n+k+1)}$. Here $\dim \n_\g(\l_0)=8(n-k)^2+2(2k-n)^2+n, \dim
\n_\h(\l_0)=2(n-k)^2+2(2k-n)^2+k.$

2) The subalgebra $\sp_{2n-2k}\times \sp_{2k}, k\geqslant
\frac{n}{2},$ is embedded into $\sp_{2n}$ as $\n_\g(\f)$, where
$\f\cong \sp_{2k}$ is embedded as described in 1). We have
$\dim\n_\g(\l_0)=4(n-k)+2(2k-n)^2+(2k-n),
\n_\h(\l_0)=3(n-k)+2(2k-n)^2+(2k-n).$

\begin{Rem}\label{Rem:6.2.2}
Let us explain why we choose this strange (at the first glance)
embedding. The pair $(\g,\h)=(\sp_{2n},\sp_{2k}\times \sp_{2(n-k)})$
is symmetric. The corresponding involution $\sigma$ acts on the
annihilator of $\sp_{2(n-k)}\subset\h$ (resp. $\sp_{2k}\subset \h$)
in $\C^{2n}$ identically (resp., by $-1$). Under the chosen
embedding, $\sigma$ acts on $\a(\g,\h)$ by $-1$, in other words,
$\a(\g,\h)$ is the Cartan space of $(\g,\h)$ in the sense of the
theory of symmetric spaces.

Further, we note that there is an embedding $\h\hookrightarrow \g$
such that $\l_0\subset\h$ but
$2(\dim\h-\n_\h(\l_0))<\dim\g-\dim\n_\g(\l_0)$. Indeed, take the
(most obvious) embedding $\sp_{2k}\times\sp_{2(n-k)}\hookrightarrow
\sp_{2n}$ such that $\sp_{2k}$  is the annihilator of $\langle
e_1,\ldots, e_{n-k},e_{2n},\ldots, e_{n-k+1}\rangle$.
\end{Rem}

3) The subalgebra  $\sl_2\times \sl_2\times \sp_{2n-4}, n\geqslant
3$, (rows NN 6,7 of Table~\ref{Tbl:1.7.4}) is embedded into
$\sp_{2n}$ as follows:  the ideal $\sp_{2n-4}\subset\h$ (resp., one
of the ideals $\sl_2\subset\h$) is the annihilator of $
e_1,e_2,e_{2n},e_{2n-1}$ (resp., of $ e_2,e_3,\ldots,e_{2n-1}$). If
$n\geqslant 8$, then
$\dim\n_\g(\l_0)=16+2(n-4)^2+(n-4),\dim\n_\h(\l_0)=6+2(n-4)^2+(n-4)$.
For  $n=3$ we have $\dim\n_\g(\l_0)=9,\dim\n_\h(\l_0)=3$.

{\it The case $\g=\so_n, n\geqslant 7$.}

1) The subalgebra $\so_k, k\geqslant \frac{n+2}{2},$ (row N8 of
Table~\ref{Tbl:1.7.4}) is embedded into $\so_n$ as the annihilator
of the vectors $e_1,e_2,\ldots, e_{l}, e_{n},\ldots, e_{n-l+1}$ (for
$n-k=2l$) or of the vectors  $e_1,e_2,\ldots, e_l,e_{n},\ldots,
e_{n-l+1}, e_{l+1}+e_{n-l}$ (for  $n-k=2l+1$). We have $\dim
\n_\g(\l_0)=\frac{(2k-n)(2k-n-1)}{2}+(n-k)(2n-2k-1),
\dim\n_\h(\l_0)=\frac{(2k-n)(2k-n-1)}{2}+\frac{(n-k)(n-k-1)}{2}.$

2) Let $n=2m$. The subalgebra  $[\h,\h]=\sl_{m}$ (rows NN 9,10 of
Table~\ref{Tbl:1.7.4} and the forth row of Table~\ref{Tbl:1.7.6}) is
embedded into $\g=\so_{2m}$ as
$\g^{(\alpha_1,\ldots,\alpha_{m-1})}$. If $m$ is even, then
$\dim\n_\g(\l_0)=3m,\dim\n_\h(\l_0)=2m-1$. For odd $m$ we have
$\dim\n_\g(\l_0)=3m-2, \dim\n_\h(\l_0)=2m-2+\dim\z(\h)$.

3) The subalgebra $\spin_7$ is embedded into  $\g=\so_9$ (row 11 of
Table~\ref{Tbl:1.7.4}) as the annihilator of the sum of a highest
vector and a lowest vector in the $\so_9$-module $V(\pi_4)$. Under
this choice of the embedding,
$\dim\n_\g(\l_0)=12,\dim\n_\h(\l_0)=9$.

4) For the embedding  $\spin_7\hookrightarrow\so_{10}$ (row 12 of
Table \ref{Tbl:1.7.4}) we take the composition of the embeddings
$\spin_7\hookrightarrow \so_9,\so_9\hookrightarrow \so_{10}$ defined
above. In this case $\n_\g(\l_0)=21,\n_\h(\l_0)=9$.

5)  $G_2$ is embedded into $\so_7$ (row 13 of Table~\ref{Tbl:1.7.4})
as the annihilator of the sum of a highest vector and a lowest
vector of the
 $\so_7$-module $V(\pi_3)$. The equalities
 $\dim\n_\g(\l_0)=9,\dim \n_\h(\l_0)=8$ hold.

6)  For the embedding $G_2\hookrightarrow\so_8$ (row 14 of
Table~\ref{Tbl:1.7.4}) we take the composition of the embeddings
$G_2\hookrightarrow \so_7,\so_7\hookrightarrow \so_8$ described
above. Here $\dim\n_\g(\l_0)=12, \dim\n_\h(\l_0)=6$.

{\it The case $\g=G_2$.}

The subalgebra $A_2$ (row 15 of Table~\ref{Tbl:1.7.4}) is embedded
into $G_2$ as $\g^{(\Delta_{max})}$, where $\Delta_{max}$ denotes
the subsystem of all long roots. In this case
$\dim\n_\g(\l_0)=6,\dim\n_\h(\l_0)=4$.

{\it The case $\g=F_4$.}

1) The subalgebra $B_4$ is embedded into $F_4$ (row 16 of
Table~\ref{Tbl:1.7.4}) as follows. The equality $\dim\n_\g(\l_0)=22$
holds.  (\ref{eq:6.2:1}) takes place iff $\l_0\cong \so_7$ is
embedded into $\h\cong \so_9$ as $\spin_7$. Put
$\alpha=\frac{\varepsilon_1+\varepsilon_2+\varepsilon_3+\varepsilon_4}{2},$
$ x=e_{\alpha}+e_{-\alpha},
\h_0:=\g^{(\{\varepsilon_i,i=\overline{1,4}\})}\cong\sp_9$. Let
$x\in \h_0^\perp$ be the sum of a highest vector and lowest vector
of the $\h_0$-module $\h_0^\perp$, so that $\z_{\g_0}(x)\cap\h_0$ is
the subalgebra $\spin_7$ in $\h_0\cong\so_9$. Now put
$\h=\Ad(g)\h_0$, where $g\in G$ is such that
$\Ad(g)x=\varepsilon_1^\vee$. One may take for such an element $g$
the product $g=g_2g_1$, where $g_1\in N_G(\t)$ is such that
$g_1\frac{\varepsilon_1+\varepsilon_2+\varepsilon_3+\varepsilon_4}{2}=\varepsilon_1$
and $g_2$ is an appropriate element of the connected subgroup of $G$
corresponding to the $\sl_2$-triple
$(e_{\varepsilon_1},\varepsilon_1^\vee,e_{-\varepsilon_1})$.

2) The subalgebra $D_4$ (row N17 of Table~\ref{Tbl:1.7.4}) is
embedded into $\g$ as $\g^{(\Delta_{max})}$, where $\Delta_{max}$
denotes the subsystem of all long roots. The equalities
$\dim\n_\g(\l_0)=16,\dim\n_\h(\l_0)=10$ hold.

{\it The case $\g=E_6$.}

1) The subalgebra $F_4$ (row N18 of Table~\ref{Tbl:1.7.4}) is
embedded into $\g$ as the annihilator of the vector
$v_{\pi_1}+v_{-\pi_5}+v_{\pi_5-\pi_1}\in V(\pi_1)$ (where
$v_\lambda$ denotes a nonzero vector of weight $\lambda$). We have
$\dim\n_\g(\l_0)=30, \dim\n_\h(\l_0)=28$.

2) The subalgebra $D_5$ (row N19 of Table~\ref{Tbl:1.7.4} and row N5
of Table~\ref{Tbl:1.7.6}) is embedded into $E_6$ as
$\g^{(\{\alpha_i\}_{i=\overline{2,6}})}$. Here $\dim
\n_\g(\l_0)=22,\dim\n_\h(\l_0)=17+\dim\z(\h)$.

3) For the embedding $B_4\hookrightarrow E_6$ (row N20 of
Table~\ref{Tbl:1.7.4}) we take the composition of the embeddings
$B_4\hookrightarrow D_5$, $D_5\hookrightarrow E_6$ described above.
The equalities $\dim\n_\g(\l_0)=38, \dim\n_\h(\l_0)=16$ take place.

4) The subalgebra $A_5$ (row 21 of Table~\ref{Tbl:1.7.4}) is
embedded into $E_6$ as $\g^{(\{\alpha_i\}_{i=\overline{1,5}})}$. We
have $\dim\n_\g(\l_0)=30,\dim\n_\h(\l_0)=11$.

{\it The case $\g=E_7$.}

1) The subalgebra $\h=E_6$ (row N22 of Table~\ref{Tbl:1.7.4}) is
embedded into $\g$ as $\g^{(\{\alpha_i\}_{i=\overline{2,7}})}$.
Under this choice of an embedding, we get
$\dim\n_\g(\l_0)=37,\dim\n_\h(\l_0)=30$.

2) The subalgebra $\h=D_6$ (row N23 of Table~\ref{Tbl:1.7.4}) is
embedded into $\g$ as $\g^{(\{\alpha_i\}_{i=\overline{1,6}})}$. We
have $\dim\n_\g(\l_0)=37,\dim\n_\h(\l_0)=18$.

{\it The case $\g=E_8$.}

 $\h=E_7$ (row $N24$ in Table~\ref{Tbl:1.7.4}) is embedded into $\g$
 as
$\g^{(\{\alpha_i\}_{i=\overline{2,8}})}$. The equalities
$\dim\n_\g(\l_0)=56,\dim\n_\h(\l_0)=37$ hold.

{\it The case $\g=\h\times\h$.}

Fix a Borel subgroup $B_0\subset H$ and a maximal torus $T_0\subset
B_0$ and construct from  them the Borel subgroup and the maximal
torus of $G$ (Section~\ref{SECTION_notation}). The embedding
 $\h\hookrightarrow\h\times\h$ is given by
$\xi\mapsto (\xi,w_0\xi)$. Here $w_0$ is an element $N_H(T_0)$
mapping into the element of maximal length of the Weyl group. We
have $\dim\n_\g(\l_0)=2\rank\h,\dim\n_\h(\l_0)=\rank\h$.

{\it The case $\g=\sp_{2n}\times \sp_{2m}, m\geqslant 1,n>1$.}

Let $e_1,\ldots,e_{2n},e_1',\ldots,e_{2m}'$ be the standard frames
in  $\C^{2n},\C^{2m}$, resp., and the symplectic forms on
$\C^{2n},\C^{2m}$ are such as indicated in
Section~\ref{SECTION_notation}. We embed $\h=\sp_{2n-2}\times
\sl_2\times \sp_{2m-2}$ into $\g$ as the stabilizer of the subspace
$\Span_\C(e_1+e_1',e_{2n}+e_{2m}')$. We get
$\dim\n_\g(\l_0)=2(n-2)^2+n-2+2(m-2)^2+(m-2)+8,
\dim\n_\h(\l_0)=2(n-2)^2+n-2+2(m-2)^2+(m-2)+3$ for $m>1$ and
$\dim\n_\g(\l_0)=2(n-2)^2+n-2+5,
\dim\n_\h(\l_0)=2(n-2)^2+n-2+2(m-2)^2+(m-2)+2$ for $m=1$.

In the following table we list the pairs $\n_\g(\l_0)/\l_0,
\n_\h(\l_0)/\l_0$ for the embeddings $\h\hookrightarrow \g$ of
subalgebras from Table \ref{Tbl:1.7.4} constructed above. In the
first column the number of a pair in Table~\ref{Tbl:1.7.4} is given.
The fourth column contains elements of $\a(\g,\h)$ constituting a
system of simple  roots in $\n_\g(\l_0)/\l_0$. The table will be
used in Subsection \ref{SUBSECTION_Weyl_aff5} to compute the Weyl
groups of affine homogeneous spaces.

\begin{longtable}{|c|c|c|c|}
\caption{The pairs
$(\n_\g(\l_0)/\l_0,\n_\h(\l_0)/\l_0)$}\label{Tbl:6.2.2}\\\hline
N&$\n_\g(\l_0)/\l_0$&$\n_\h(\l_0)/\l_0$&simple roots of
$\n_\g(\l_0)/\l_0$\\\endfirsthead\hline
N&$\n_\g(\l_0)/\l_0$&$\n_\h(\l_0)/\l_0$&simple roots of
$\n_\g(\l_0)/\l_0$
\\\endhead\hline 1&$\sl_{2(n-k)}$&$\sl_{n-k}$&$\varepsilon_{i}-\varepsilon_{i+1},i<n-k$ or $i\geqslant k$,
$\varepsilon_{n-k}-\varepsilon_k$\\\hline 2&$\sl_2^{n-k}\oplus
\C$&$\C^{n-k}:(\C,\ldots,\C,0)$&$\varepsilon_i-\varepsilon_{n+1-i},
i\leqslant n-k.$
\\\hline 3&$\C^n$&$0$&\\\hline
4&$\sp_{4(n-k)}$&$\sp_{2(n-k)}$&$\varepsilon_i-\varepsilon_{i+1},
i<2(n-k),\varepsilon_{2(n-k)}$\\\hline 5&$\C^{n-k}$&$0$&\\\hline
6&$\sl_{4}$&$\sl_2\times\C^2$&$\varepsilon_1-\varepsilon_4,\varepsilon_3+\varepsilon_4,
\varepsilon_2-\varepsilon_3$\\\hline
7&$\sl_3$&$\C^2$&$\varepsilon_1+\varepsilon_3,\varepsilon_2-\varepsilon_3$\\\hline
 8&$\so_{2(n-k)}$&$\so_{n-k}$&$\varepsilon_i-\varepsilon_{i+1},i<n-k,\varepsilon_{n-k+1}+\varepsilon_{n-k}$\\\hline
9&$\sl_2^n$&$\C^n$&$\varepsilon_{2i-1}+\varepsilon_{2i},
i=\overline{1,n}$
\\\hline 10&$\sl_2^n\oplus\C
$&$\C^{n+1}$&$\varepsilon_{2i-1}+\varepsilon_{2i},i=\overline{1,n}$\\\hline
11&$\sl_2\times \C$&$\C$:$(\C,\C)$&$\varepsilon_1$\\\hline
12&$\so_6\times\sl_2$&$\sl_2\times\sl_2$:$(\so_4,\sl_2)$&$\varepsilon_1-\varepsilon_2,
\varepsilon_2-\varepsilon_5, -\varepsilon_1-\varepsilon_2,
\varepsilon_3+\varepsilon_4$\\\hline 13&$\C$&$0$&\\\hline
14&$\sl_2^3$&$\sl_2$:
$(\sl_2,\sl_2,\sl_2)$&$\varepsilon_1-\varepsilon_4,\varepsilon_1+\varepsilon_4,
\varepsilon_2+\varepsilon_3$\\\hline
15&$\sl_2$&$\C$&$\varepsilon_1$\\\hline
16&$\sl_2$&$\C$&$\varepsilon_1$\\\hline
17&$\sl_3$&$\C^2$&$(\varepsilon_1\pm(\varepsilon_2+\varepsilon_3\pm\varepsilon_4))/2$\\\hline
18&$\C^2$&$0$&\\\hline 19&$\sl_2^2\times \C$&$\C^2:
(\C,\C,\C)$&$\varepsilon_1-\varepsilon_6,2\varepsilon$\\\hline
20&$\sl_6$&$\sp_4\times\sl_2$&$\varepsilon-\varepsilon_1-\varepsilon_3-\varepsilon_4, \varepsilon_1-\varepsilon_2,$\\
&&&$\varepsilon_2-\varepsilon_5,\varepsilon_5-\varepsilon_6,
\varepsilon+\varepsilon_3+\varepsilon_4+\varepsilon_6$\\\hline
21&$\so_8$&$\sl_2^3$&$\varepsilon_1-\varepsilon_6,\varepsilon+\varepsilon_4+\varepsilon_5+\varepsilon_6,
\varepsilon_2-\varepsilon_5,\varepsilon_3-\varepsilon_4$\\\hline
22&$\sl_2^3$&$\C^2: (\C,\C,\C)$&$\varepsilon_1-\varepsilon_2,
\varepsilon_3+\varepsilon_4+\varepsilon_5+\varepsilon_6,\varepsilon_8-\varepsilon_7$\\\hline
23&$\so_8$&$\sl_2^3$&$\varepsilon_3+\varepsilon_4+\varepsilon_7+\varepsilon_8,\varepsilon_6-\varepsilon_7,\varepsilon_5
-\varepsilon_6,\varepsilon_1+\varepsilon_2+\varepsilon_7+\varepsilon_8$\\\hline
24&$\so_8$&$\sl_2^3$&$\varepsilon_2-\varepsilon_3,\varepsilon_1-\varepsilon_2,-\varepsilon_1-\varepsilon_2-\varepsilon_9,
\varepsilon_2+\varepsilon_3+\varepsilon_8$\\\hline 25&$\C^{\rank
\h}$&$0$&\\\hline 26&$\sl_2^2\times \C$&$\C^2:
(\C,\C,\C)$&$\varepsilon_1+\varepsilon_2,
\varepsilon_1'+\varepsilon_2'$\\\hline 27& $\sl_2\times\C$&$\C:
(\C,\C)$&$\varepsilon_1+\varepsilon_2$\\\hline
\end{longtable}

In the brackets after $\n_\h(\l_0)/\l_0$ in the third columns we
indicate the projections of $\n_\h(\l_0)/\l_0$ to simple ideals of
$\n_\g(\l_0)/\l_0$ provided the last algebra is not simple and
$\n_\h(\l_0)/\l_0\neq 0$.

\subsection{Proof of Proposition~\ref{Prop:6.0.1}}\label{SUBSECTION_disting_compI3}
\begin{proof}
{\it Assertion  1.} Note that
$L_{0\,G,G/H^{ess}}^\circ=L_{0\,G,G/H}^\circ$. It is clear that
$\pi(\underline{X}')\subset X^{L_0}$. Moreover,
$\overline{U\pi(\underline{X}')}=\pi(\overline{U\underline{X}'})=\pi(G/H^{ess})=X$.
Therefore $\pi(\underline{X}')\subset \underline{X}$.

{\it Assertion 2.}  Put $L=L_{G,G/H},P=LB,
\widetilde{L}_0=L_{0\,G,G/H^{sat}}^\circ$. Since
$L_{G,G/H^{sat}}=L$, there is an $L$-stable subvariety $S\subset
G/H^{sat}$ such that $(L,L)$ acts trivially on $S$ and the natural
morphism $P*_L S\rightarrow G/H^{sat}$ is an embedding, see
Proposition \ref{Prop:1.5.6}. Thence there is a $P$-embedding $P*_L
\pi^{-1}(S)\hookrightarrow X$. It follows that
$\X_{G,X}=\X_{L,\pi^{-1}(S)}$. Hence the action $(L,L):\pi^{-1}(S)$
is trivial. By Proposition~\ref{Prop:1.5.6},
$\underline{X}'=\overline{\Rad_u(P)^{\widetilde{L}_{0}}S},
\underline{X}=\overline{\Rad_u(P)^{L_0}\pi^{-1}(S)}$. Therefore
$\Rad_u(P)^{\widetilde{L}_0}\pi^{-1}(S)$ is a dense subset of
$\pi^{-1}(\underline{X}')$.  This proves assertion 2.

{\it Assertion 3.} This is obvious.

{\it Assertion 4.} Put $L=L_{G,G/H}$. By
Proposition~\ref{Prop:6.0.2}, the assertion will follow if we prove
\begin{equation}\label{eq:6.3:1}N_G(L_0)=N_G(L_0)^\circ N_H(L_0).\end{equation}
Note that  $Z(G)\subset L\subset N_G(L_0)^\circ N_H(L_0)$. Therefore
in the proof we may always consider any  group $G$ with Lie algebra
$\g$.

The proof of (\ref{eq:6.3:1}) will be carried out in two steps. Put
$N_G(L_0)_0=Z_G(L_0)L_0, N_H(L_0)_0=Z_H(L_0)L_0$. These are the
subgroups in the corresponding normalizers consisting of all
elements acting on $\l_0$ by inner automorphisms. On the first step
we show that
 $N_G(L_0)_0=N_H(L_0)_0 N_G(L_0)^\circ$ and on the second one that
$N_G(L_0)=N_G(L_0)_0N_H(L_0)$. The claim of the first step will
follow if we check that $Z_G(\l_0)/Z_H(\l_0)$ is connected.

{\it Step  1.}
\begin{Lem}\label{Lem:6.3.1}
Let   $\m_0$ be a subalgebra of   $\g$ contained in some Levi
subalgebra $\m\subset \g$ and containing $[\m,\m]$. Suppose that any
element  of $\Delta(\m)$ is a long root in $\Delta(\g)$ (more
precisely, a long root in $\Delta(\g_i)$, where $\g_i$ is some
simple ideal of $\g$). Then $Z_G(\m_0)M_0$ is connected, where $M_0$
denotes the connected subgroup of $G$ corresponding to $\m_0$.
\end{Lem}
\begin{proof}[Proof of Lemma~\ref{Lem:6.3.1}]
The claim is well-known in the case when $\m_0$ is commutative
(see~\cite{VO}, ch.3, $\S3$). Using induction on $\dim\m_0$, we may
assume that $\m_0$ is simple. Choose a Cartan subalgebra $\t\subset
\m$. Any connected component of $Z_G(M_0)M_0$ contains an element
from $N_G(\t)\cap Z_G(\t\cap \m_0)$. An element of the last group
acts on $\t$ as a composition of  reflections corresponding to
elements of $\Delta(\m)^\perp\cap\Delta(\g) $. Since $\beta$ is a
long root,  we see that $\alpha+\beta\not\in\Delta(\g) $ for any
$\alpha\in \Delta(\m)^\perp, \beta\in \Delta(\m)$. Therefore
$\Delta(\m)^{\perp}\cap \Delta(\g) \subset \Delta(\z_\g(\m_0))$. In
particular, $Z_G(\m_0\cap\t)\cap N_G(\t) \subset Z_G(\m_0)^\circ
M_0$.
\end{proof}

It remains to consider only pairs $(\g,\h)$ such that $\g$ contains
a simple ideal isomorphic to $B_l,C_l,F_4,G_2$ and $\Delta(\l)$
contains a short root. Such a pair $(\g,\h)$ is on of the following
pairs from Table~\ref{Tbl:1.7.4}: NN4 ($2k-n>1$),5,6 ($n>5$),8 ($n$
is odd),16,26 ($n\geqslant 3$),27 ($n\geqslant 3$). In cases 4,6,27
the subalgebra  $[\l_0,\l_0]$ coincides with $\sp_{2l}$, while in
case 26 it coincides with the direct sum of the subalgebras
$\sp_{2l_1},\sp_{2l_2}$ embedded into different simple ideals of
$\g$. The centralizer of $\sp_{2l}\subset \sp_{2n}$ in $\Sp(2n)$ is
connected for all
  $l,n$. It remains to consider cases 5,8,16.

Put $\g=\sp_{2m}, \h=\sp_{2k}\times \sp_{2(m-k)}, k\geqslant n/2$.
As in the previous paragraph, it is enough to consider the case
$n=2k$. Let us describe the embedding $L_0\hookrightarrow \Sp(2n)$.
The space $\C^{4k}$ is decomposed into the direct sum of
two-dimensional spaces $V_1,\ldots,V_{2k}$. We may assume that
$H=\Sp(V_1\oplus V_3\oplus\ldots\oplus V_{2k-1})\times
\Sp(V_{2}\oplus V_4\oplus \ldots\oplus V_{2k})$. The subgroup
$L_0\subset \Sp_{2n}$ is isomorphic to the direct product of  $k$
copies of $\SL_2$. The $i$-th copy of $\SL_2$ (we denote it by
$L_0^i$) acts diagonally on $V_{2i-1}\oplus V_{2i}$ and trivially on
$V_j, j\neq 2i-1,2i$. It can be seen directly that
$Z_G(L_0)=Z_{\Sp(V_1\oplus V_2)}(L_0^1)\times \ldots\times
Z_{\Sp(V_{2k-1}\oplus V_{2k})}(L_0^k)$. Therefore
$Z_G(L_0)/Z_H(L_0)$ is the $k$-dimensional torus.

In case 8 both groups $Z_G(L_0),Z_H(L_0)$ contain two connected
components. The components of unit in the both cases consist of all
elements acting trivially on $\C^{2n+1}/(\C^{2n+1})^{\l_0}$. This
proves the required claim.

In case 16 we have $Z_H(\l_0)/Z_H(\l_0)^\circ\cong \Z_2$. An element
$\sigma\in(Z_H(\l_0)\setminus Z_H(\l_0)^\circ)\cap N_G(\t)$ acts on
$\t\cap \l_0^\perp$ by -1. Any element of $Z_G(\l_0)\cap N_G(\t)$
acts on $\t\cap\l_0^\perp$ by  $\pm 1$. Since
$Z_G(\l_0)=Z_G(\l_0)^\circ(Z_G(\l_0)\cap N_G(\t))$, we see that the
group $Z_G(L_0)/Z_H(L_0)$ is connected.

{\it Step  2.} It remains to check that $N_G(L_0)=N_G(L_0)_0
N_H(L_0)$, equivalently, the images of $N_G(\l_0),N_H(\l_0)$ in
$\GL(\l_0)$ coincide. For the pair N1 in Table \ref{Tbl:1.7.4} the
group
 $N_G(\l_0)$ is connected. If $(\g,\h)$ is one of the pairs NN 4,8 ($n$ is odd), 12,14,15,16,20,
then the group of outer automorphisms of $\l_0$ is trivial. For the
pairs NN11,13,17,18,19,22,24 from Table~\ref{Tbl:1.7.4} the algebra
$\l_0$ is simple and $N_H(\l_0)$ contains all automorphisms of
$\l_0$ (in cases  22,24 one proves the last claim using the
embeddings $D_4\subset F_4\subset E_6\subset E_7\subset E_8$). In
cases 3,9,22, 5$(n=2k)$ the algebra $\l_0$ is isomorphic to the
direct product of some copies of $\sl_2$. The group of outer
automorphisms of $\l_0$ is the symmetric group on the set of simple
ideals of $\l_0$. The image of $N_H(\l_0)$ in
$\Aut(\l_0)/\Int(\l_0)$ coincides with the whole symmetric group.
One considers case 5, $n>2k$, analogously. Here
$\Aut(\l_0)/\Int(\l_0)$ is isomorphic to  the symmetric group on the
set of simple ideals of $\l_0$ that are isomorphic to $\sl_2$ and
have Dynkin index 2. In case 8 (with even $n$) the images of both
$N_H(\l_0),N_G(\l_0)$ in $\Aut(\l_0)$ contain an outer automorphism
of $\l_0$ induced by an element from $\O_{2k-n}\setminus \SO_{2k-n}$
and do not contain elements from the other nontrivial connected
components of $\Aut(\l_0)$ (which exist only when $2k-n=8$).

It remains to consider the pairs    NN2,6,7,21,25,26,27 from Table
\ref{Tbl:1.7.4} and the pairs $(\g,\h)$, where
$\z(\h)=\z(\z_\g([\h,\h]))$ and $(\g,[\h,\h])$ is listed in
Table~\ref{Tbl:1.7.6}. In all these cases $\z(\l_0)\neq \{0\}$. If
$(\g,[\h,\h])\neq (\so_{4n+2},\sl_{2n+1}),(E_6,D_5)$, then
$N_G([\l_0,\l_0])$ is connected. As we have seen above, in these two
cases the equality $N_G([\l_0,\l_0])=N_G([\l_0,\l_0])_0
N_H([\l_0,\l_0])$ holds. Therefore it is enough to show that
\begin{equation}\label{eq:6.3:2}N_{Z_G([\l_0,\l_0])}(\z(\l_0))=
Z_G(\l_0)N_{Z_H([\l_0,\l_0])}(\z(\l_0)).\end{equation}

If a pair $(\g,[\h,\h])$ is contained in Table~\ref{Tbl:1.7.6},
then, using Lemma \ref{Lem:6.3.1}, one can see that
$Z_G([\l_0,\l_0])$ is connected. If, additionally,
 $(\g,[\h,\h])\neq (\sl_n,\sl_{n-k}\times \sl_k)$,  then the group $N_{Z_G([\l_0,\l_0])}(\z(\l_0))$
is connected too. In cases 6,7,26,27 the algebra $\z(\l_0)$ is
one-dimensional and the groups in the both sides of (\ref{eq:6.3:2})
act on  $\z(\l_0)$ as $\Z_2$. This observation yields
(\ref{eq:6.3:2}) in these cases. It remains to check
(\ref{eq:6.3:2}) for
$(\g,[\h,\h])=(\h\times\h,\h),(E_6,A_5),(\sl_n,\sl_{n-k}\times\sl_k)$.
The first case is obvious.

{\it The case $(\g,\h)=(E_6,A_5)$.} We suppose that the invariant
symmetric form on $\g$ is chosen in such a way that the length of a
root equals 2.

We have $\l_0=\langle \alpha_1-\alpha_5,\alpha_2-\alpha_4\rangle$.
The action of  $N_H(\l_0)$  on $\l_0$ coincides with the action of
the symmetric group $S_3$ on its unique 2-dimensional module. To see
this note that $\l_0$ is embedded into $\sl_6$ as
$\{diag(x,y,-x-y,-x-y,y,x)\}$.

Any element of length 4 lying in the intersection of $\l_0$ and the
root lattice of $\g$ is one of the following elements:
$\pm(\alpha_1-\alpha_5),\pm(\alpha_2-\alpha_4),\pm(\alpha_1+\alpha_2-\alpha_4-\alpha_5)$.
Assume that the images of $N_G(\l_0)$ and $N_H(\l_0)$ in $\GL(\l_0)$
differ. The image $N$ of $N_G(\l_0)\cap N_G(\t)$ in $\GL(\l_0)$
permutes the six elements of length 4 listed above and preserves the
scalar product on $\l_0\cap \t(\R)$. It follows that $-id\in N$. Let
us check that $g\in N_G(\t)\cap N_G(\l_0)$ cannot act on $\l_0$ by
$-1$.

Assume the converse, let $g$ be such an element. Clearly, $g\in
N_G(\z_\g(\alpha_1-\alpha_5))$. Note that
$\z_\g(\alpha_1-\alpha_5)\cong\C\times \sl_2\times \sl_2\times
\sl_4$. Simple roots of the ideals isomorphic to $\sl_2$ are
$\varepsilon_1-\varepsilon_6,\varepsilon_2-\varepsilon_5$. Thus $g$
preserves the pair of lines
$\{\langle\varepsilon_1-\varepsilon_6\rangle,\langle\varepsilon_2-\varepsilon_5\rangle\}$.
Replacing   $\alpha_1-\alpha_5$ with $\alpha_2-\alpha_4$ and
$\alpha_1+\alpha_2-\alpha_4-\alpha_5$ in the previous argument, we
see that $g$ preserves the pairs $\{\langle
\varepsilon_2-\varepsilon_5\rangle,\langle \alpha_3 \rangle\}$ and
$\{\langle\varepsilon_1-\varepsilon_6\rangle,\langle\alpha_3\rangle\}$.
Thence  $\langle \varepsilon_1-\varepsilon_6\rangle, \langle
\varepsilon_2-\varepsilon_5\rangle, \langle \alpha_3\rangle$ are
$g$-stable. Since   $g$ is an orthogonal transformation of $\t$
leaving $\langle \alpha_1-\alpha_5,\alpha_2-\alpha_4\rangle$
invariant,  we see that the line $\langle \varepsilon\rangle$ is
$g$-stable. Replacing $g$ with  $gs_{2\varepsilon}$, we may assume
that $g\in Z_G(\langle \varepsilon\rangle)$. The last group is the
product of $A_5\subset E_6$ and a one-dimensional torus and does not
contain an element acting on $\l_0$ by $-1$.

{\it The case $(\g,[\h,\h])=(\sl_n,\sl_k\times \sl_{n-k})$.} The
subalgebra $\l_0\subset\h$ consists of all matrices of the form
$diag(x_1,\ldots,x_{n-k},A,x_{n-k},\ldots x_1)$, where
$x_1,\ldots,x_{n-k}\in \C, A\in \gl_{2k-n}$,
$\tr(A)=-2\sum_{i=1}^{n-k}x_i$ for $\z(\h)\neq 0$ and
$\tr(A)=\sum_{i=1}^{n-k}x_i=0$ for $\z(\h)=0$. The groups
$N_G(\l_0),N_H(\l_0)$ act on $\z(\l_0)$ as $S_{n-k}$ (permuting
$x_i$).

{\it Assertion 5.} This follows directly from assertion 4 and
Propositions~\ref{Prop:6.0.2}, \ref{Prop:6.1.2}.
\end{proof}

In Subsection~\ref{SUBSECTION_Weyl_aff5} we will need to know the
image of $N_G(\l_0)\cap N_G(\t)$ in $\GL(\a(\g,\h))$ for pairs
$(\g,\h)$ from Table~\ref{Tbl:1.7.4}. This information is extracted
mostly from the previous proof and Table~\ref{Tbl:6.2.2}. It is
presented in Table~\ref{Tbl:6.3.2}.

\begin{longtable}{|c|c|}
\caption{The image of $N_G(\l_0)\cap N_G(\t)$ in
$\GL(\a_{\g,\h})$}\label{Tbl:6.3.2}\\\hline N&The
group\\\endfirsthead\hline N&The group
\\\endhead\hline
1&$W$\\\hline
2&$s_{\varepsilon_i-\varepsilon_j-\varepsilon_{n+1-i}+\varepsilon_{n+1-j}},
1\leqslant i<j\leqslant n-k$\\\hline
3&$s_{\varepsilon_{2i-1}+\varepsilon_{2i}\pm
(\varepsilon_{2j-1}+\varepsilon_{2j})}, 1\leqslant i,j\leqslant
n$\\\hline 4&$W$\\\hline
5&$s_{\varepsilon_{2i-1}+\varepsilon_{2i}\pm
(\varepsilon_{2j-1}+\varepsilon_{2j})}, 1\leqslant i,j\leqslant
n-k$\\\hline 6&$A$\\\hline 7&$A$\\\hline 8&$A$\\\hline 9&$A$\\\hline
10&$A$\\\hline 11&$A\oplus \Z_2$\\\hline 12&$W$\\\hline
13&$\Z_2$\\\hline 14&$W$\\\hline 15&$\Z_2$\\\hline 16&$\Z_2$\\\hline
17&$A$\\\hline 18&$s_{\pi_1},s_{\pi_5}$\\\hline 19&$A$\\\hline
20&$W$\\\hline 21&$A$\\\hline 22&$A$\\\hline 23&$A$\\\hline
24&$A$\\\hline 25&$W(\h)$\\\hline 26&$W\oplus\Z_2$\\\hline
27&$W\oplus \Z_2$\\\hline
\end{longtable}

Let us explain the notation used  in the table. In rows 2,3,5,18 a
set of generators of the group is given. The symbols $A,W$ mean the
automorphism group (resp., the Weyl group) of the root system
$\n_\g(\l_0)/\l_0$ (we assume that this group acts trivially on
$\z(\n_\g(\l_0)/\l_0)$). The symbol $\Z_2$ denotes the group acting
by $\pm 1$ on the center of $\n_\g(\l_0)/\l_0$ and trivially on the
semisimple part.

\section{Computation of  Weyl groups for affine homogeneous vector bundles}\label{SECTION_Weyl_aff}
\subsection{Introduction}\label{SUBSECTION_Weyl_aff_intro}
The main goal of this section is to compute the group $W_{G,X}$,
where $X=G*_HV$ is a homogeneous vector bundle over an affine
homogeneous space $G/H$. We recall that $W_{G,X}$ depends only on
the triple $(\g,\h,V)$ (see Corollary~\ref{Cor:1.2.4}), so we write
$W(\g,\h,V)$ instead of $W_{G,G*_HV}$.

Let $\underline{X}$ be the distinguished component of $X^{L_0}$,
where $L_0=L_{0\,G,X}^\circ$, and
$\underline{G}=N_G(L_0,\underline{X})/L_0$. Theorem~\ref{Thm:3.0.3}
allows to recover $W_{G,X}$ from
$W_{\underline{G}^\circ,\underline{X}}$. Besides, the results of
Section~\ref{SECTION_distinguished} show that the
$\underline{G}^\circ$-variety $\underline{X}$ is an affine
homogeneous vector bundle and allow to determine it. So to compute
Weyl groups one may restrict to the case  $\rank_G(X)=\rank G$.
Further, Proposition \ref{Prop:3.4.5} reduces the computation of
$W_{G,X}$ to the case of simple  $G$. Finally, we may assume that
$G$ is simply connected. Below in this subsection we always assume
that these conditions are satisfied.

By an {\it admissible triple} we mean a triple $(\g,\h,V)$, where
$\g$ is a simple Lie algebra, $\h$ its reductive subalgebra, and $V$
is a module over the connected subgroup $H\subset G$ with Lie
algebra $\h$ such that $\a(\g,\h,V)=\t$.

Let us state our main result.  There is a minimal ideal $\h_0\subset
\h$ such that $W(\g,\h,V)=W(\g,\h_0,V/V^{\h_0})$. We say that the
corresponding triple $(\g,\h_0,V/V^{\h_0})$ is a {\it $W$-essential
part} of the triple $(\g,\h,V)$. The problem of computing of
$W(\g,\h,V)$ may be divided into three parts:

a) To find all admissible triples  $(\g,\h,V)$ such that
$V^{\h}=\{0\}$ and  $W(\g,\h,V)\neq W(\g,\h_0,V)$ for any  ideal
$\h_0\subsetneq \h$. Such a triple $(\g,\h,V)$ is called {\it
$W$-essential}. Clearly, a triple is  $W$-essential iff it coincides
with its  $W$-essential part.

b) To compute the groups $W(\g,\h,V)$ for all triples  $(\g,\h,V)$
found on the previous step.

c) To show how one can determine a $W$-essential part of a given
admissible triple.

\begin{defi}\label{Def:5.0.1}
Two triples $(\g,\h_1,V_1),(\g,\h_2,V_2)$ are said to be {\it
isomorphic} (resp., {\it equivalent}), if there exist $\sigma\in
\Aut(\g)$ (resp., $\sigma\in \Int(\g)$) and a linear isomorphism
$\varphi:V_1/V_1^{\h_1}\rightarrow V_2/V_2^{\h_2}$ such that
$\sigma(\h_1)=\h_2$ and $\varphi(\xi v)=\sigma(\xi)\varphi(v)$ for
all $\xi\in\h,v\in V_1/V_1^{\h_1}$.
\end{defi}

Let $H_1,H_2$ be the connected subgroups of $G$ corresponding to
$\h_1,\h_2\subset \g$. Put $X_{01}=G*_{H_1}V_1, X_{02}=G*_{H_2}V_2$.
The triples $(\g,\h_1,V_1),(\g,\h_2,V_2)$ are isomorphic (resp.,
equivalent) iff $X_{01}\cong ^\tau\!\!X_{02}$ for some  $\tau\in
\Aut(G)$ (resp., $X_{01}\cong X_{02}$). Lemma~\ref{Lem:1.5.10}
allows to compute the Weyl group  only for one triple in a given
class of isomorphism.

There is a trivial $W$-essential triple $(\g,0,0)$ and
$W(\g,0,0)=W(\g)$.

\begin{Thm}\label{Thm:5.0.2}
Let $(\g,\h,V)$ be an admissible triple.
 If there exists an ideal
$\h_1\subset \h$ such that the triple $(\g,\h_1,V/V^{\h_1})$ is
isomorphic to a triple from Table~\ref{Tbl:5.0.3}, then
$(\g,\h_1,V/V^{\h_1})$ is a $W$-essential part of  $(\g,\h,V)$ and
the Weyl group is presented in the fifth column of Table
\ref{Tbl:5.0.3}. Otherwise, $(\g,0,0)$ is a $W$-essential part of
$(\g,\h,V)$.
\end{Thm}

\begin{longtable}{|c|c|c|c|c|}\caption{$W$-essential triples $(\g,\h,V)$ and the corresponding Weyl groups}\label{Tbl:5.0.3}
\\\hline N&$\g$&$\h$&$V$&$W(\g,\h,V)$\\\endfirsthead\hline
N&$\g$&$\h$&$V$&$W(\g,\h,V)$\\\endhead\hline
1&$\sl_n$&$\sl_{n-k}$&$l\tau+m\tau^*$&$\varepsilon_{i}-\varepsilon_{i-1},\varepsilon_{k+m}-\varepsilon_{k+m+2}
$\\&$n>1$&$k<\frac{n}{2}$&$l+m=n-2k-1,l\geqslant m$&$i\neq
k+m,k+m+1$
\\\hline 2&
$\sl_n$&$\sl_n$&$\bigwedge^2\tau+\tau$&$\varepsilon_i-\varepsilon_{i+2}$\\&$n>3$&&&$i=\overline{1,n-2}$\\\hline
3& $\sl_{n}$&$\sl_n$&$\bigwedge^2\tau+\tau^*$&N2\\&even $n\geqslant
4 $&&&\\\hline
4&$\sl_{n}$&$\sl_{n}$&$\bigwedge^2\tau+\tau^*$&$\varepsilon_i-\varepsilon_{i+2},\varepsilon_1-\varepsilon_2$\\&
odd $n\geqslant 5 $&&&$i=\overline{2,n-1}$\\\hline
5&$\sl_{n}$&$\sl_{n-1}$&$\bigwedge^2\tau$&$\varepsilon_i-\varepsilon_{i+2},\varepsilon_1-\varepsilon_2$\\&even
$n\geqslant 4$&&&$i=\overline{2,n-1}$\\\hline
6&$\sl_{n}$&$\sl_{n-1}$&$\bigwedge^2\tau$&N2\\& odd $n\geqslant
3$&&&\\\hline 7&$\sl_{n}$&$\sp_{n}$&$\tau$&N2\\ &even $n\geqslant
4$&&&\\\hline 8&$\sl_{n}$&$\sp_{n-1}$&0&N2\\ &odd $n\geqslant 5$&&&
\\\hline 9&$\so_{2n+1}$&$\sl_n^{diag}$&0&$\varepsilon_i-\varepsilon_{i+2}, \varepsilon_{n-1},\varepsilon_n$\\&$n\geqslant 3$&
&&$i=\overline{1,n-2}$\\\hline
10&$\so_7$&$\so_7$&$l\tau+(2-l)R(\pi_3), l=0,1$&N9\\\hline
11&$\so_7$&$\so_6$&$kR(\pi_3)+(2-k)R(\pi_1),k>0$&$k=1:
$N9\\&&&&$k=2:
\varepsilon_1-\varepsilon_2,\varepsilon_2,\varepsilon_3$\\\hline
12&$\so_7$&$G_2$&$R(\pi_1)$&N9\\\hline
13&$\so_9$&$\so_9$&$R(\pi_4)+\tau$&$\varepsilon_1-\varepsilon_2,\varepsilon_2-\varepsilon_4,\varepsilon_3,\varepsilon_4$\\\hline
14&$\so_9$&$\so_8$&$(2-k-l)\tau+l R(\pi_3)+k R(\pi_4)$&N13,
$(k,l)\neq (2,0)$\\&&& $(k,l)=(1,0),(2,0),(1,1)$&N9,
$(k,l)=(2,0)$\\\hline 15&$\so_9$&$\so_7$&$R(\pi_3)$&N13\\\hline
16&$\so_9$&$\mathfrak{spin}_7$&$R(\pi_1)$&N9\\\hline
17&$\so_9$&$\mathfrak{spin}_7$&$R(\pi_3)$&N13\\\hline
18&$\so_9$&$G_2$&0&N13\\\hline
19&$\so_{11}$&$\so_{11}$&$R(\pi_5)$&N9
\\\hline
20&$\so_{11}$&$\so_{10}$&$R(\pi_1)+R(\pi_4)$&$\varepsilon_i-\varepsilon_{i+1},i=1,2,3$\\
&&&&$\varepsilon_4,\varepsilon_5$\\\hline
21&$\so_{11}$&$\so_9$&$R(\pi_4)$&N20\\\hline
22&$\so_{11}$&$\so_8$&$R(\pi_3)$&N20\\\hline
23&$\so_{11}$&$\mathfrak {spin}_7$&0&N20\\\hline
24&$\so_{13}$&$\so_{10}$&$R(\pi_4)$&$\varepsilon_i-\varepsilon_{i+1},i=1,2,3$\\
&&&&$\varepsilon_4,\varepsilon_5$\\\hline
25&$\sp_{2n}$&$\sp_{2k}$&$(2k-n)\tau$&$\varepsilon_i-\varepsilon_{i+1},\varepsilon_{n-1}+\varepsilon_n$\\&$n\geqslant
2$&$k\geqslant \frac{n}{2}$&&$i=\overline{1,n-1}$\\\hline
26&$\so_{4n+2}$&$\sl_{2n+1}^{diag}$&$\tau^*$&$\varepsilon_i+\varepsilon_{i+1}$\\&$n\geqslant
2$&&&$i=\overline{1,2n}$\\\hline
27&$\so_{4n}$&$\sl_{2n}^{diag}$&$\tau^*$&$\varepsilon_i+\varepsilon_{i+1}$\\&$n\geqslant
2$&&&$i=\overline{1,2n-1}$\\\hline
28&$G_2$&$A_2$&$\tau$&$\varepsilon_1,\varepsilon_2$\\\hline
29&$F_4$&$B_4$&$2\tau$&$\varepsilon_3,\frac{\varepsilon_1+\varepsilon_2-\varepsilon_3-\varepsilon_4}{2},$
\\&&&&$\varepsilon_4,\varepsilon_2-\varepsilon_4$\\\hline
30&$F_4$&$D_4$&$\tau$&N29\\\hline 31&$F_4$&$B_3$&0&N29\\\hline
32&$\sp_{4n}$&$\sp_{2n}\times
\sp_{2n}$&$R(\pi_1)$&$\varepsilon_i+\varepsilon_{i+1}$\\&$n\geqslant
1 $&&&$i=\overline{1,2n-1}$\\\hline
33&$\sp_{4n+2}$&$\sp_{2n+2}\times\sp_{2n}$&$R(\pi_1)$&$\varepsilon_i+\varepsilon_{i+1}$\\&$n\geqslant
1 $&&&$i=\overline{1,2n}$\\\hline
\end{longtable}

Let us explain the notation used in the table.  In the fourth column
the representation of $\h$ in $V$ is given.  $\tau$ denotes the
tautological representation of $\h$. In column 5 we list roots such
that the corresponding reflections generate $W(\g,\h,V)$. "N$k$" in
column 5 means that the corresponding Weyl group coincide with that
from row N$k$. If in row 1 the lower index $j$ of $\varepsilon_j$ is
less than  1 or bigger than $n$, then the corresponding root should
be omitted. In row 27 we suppose that
$\h=\g^{(\alpha_1,\ldots,\alpha_{2n-1})}$.

\begin{Rem}
Inspecting Table~\ref{Tbl:5.0.3}, we deduce from Theorem
\ref{Thm:5.0.2} that a
 $W$-essential part of  $(\g,\h,V)$ is uniquely determined.
\end{Rem}

\begin{Rem}\label{Rem:5.0.4}
Let $(\g,\h,V)$ be a triple listed in Table~\ref{Tbl:5.0.3}. The
isomorphism class of $(\g,\h,V)$ consists of more than one
equivalence class precisely for the following triples:
 N1 ($l\neq m$), NN2-6, 26,27.
In these cases an isomorphism class consists of  two different
equivalence classes.
\end{Rem}

Now we describe the content of this section. In
Subsection~\ref{SUBSECTION_Weyl_aff2} we classify
 {\it $W$-quasiessential triples}.

\begin{defi}\label{Def:5.0.5}
An admissible triple $(\g,\h,V)$ is called  {\it $W$-quasiessential}
if for any proper ideal $\h_1\subset \h$ there exists a root
$\alpha\in \Delta(\g) $ such that $S^{(\alpha)}\rightsquigarrow_\g
T^*(G*_HV)$ but $S^{(\alpha)}\not\rightsquigarrow_\g
T^*(G*_{H_1}V)$. Here $H,H_1$ are connected subgroups of $G$
corresponding to $\h,\h_1$.
\end{defi}

Below we will see  that $W$-quasiessential triples are precisely
those listed in Table~\ref{Tbl:5.0.3}.

In Subsection~\ref{SUBSECTION_Weyl_aff3} we will compute the Weyl
groups for all triples listed in Table~\ref{Tbl:5.0.3}.
Subsection~\ref{SUBSECTION_Weyl_aff4} completes the proof of
Theorem~\ref{Thm:5.0.2}. Finally, in
Subsection~\ref{SUBSECTION_Weyl_aff5} we compute the Weyl groups of
affine homogeneous spaces (without restrictions on the rank) more or
less explicitly.

\subsection{Classification of $W$-quasiessential triples}\label{SUBSECTION_Weyl_aff2}
In this subsection $\g$ is a simple Lie algebra.

\begin{Prop}\label{Prop:5.3.1}
\begin{enumerate}
\item
An admissible triple $(\g,\h,V)$ is $W$-quasiessential iff it is
listed in Table~\ref{Tbl:5.0.3}.
\item Let $(\g,\h,V)$ be a $W$-quasiessential triple. If $\h$ is
simple and  $S^{(\alpha)}\rightsquigarrow_\g T^*(G*_HV)$, then
$\alpha$ is a long root. If $\h$ is not simple, then
$S^{(\alpha)}\rightsquigarrow_\g T^*(G*_HV)$ for all roots $\alpha$.
\end{enumerate}
\end{Prop}

To prove Proposition~\ref{Prop:5.3.1} we need some technical
results.

Let us introduce some notation. Let $H$ be a reductive algebraic
group, $\s$ a subalgebra of $\h$ isomorphic to $\sl_2$. We denote by
$S^\s\index{ss@$S^\s$}$ the $\h$-stratum consisting of  $\s$ and the
direct sum of two copies of the two-dimensional irreducible
$\s$-module.

\begin{Rem}\label{Rem:5.3.2}
Let $H$ be a reductive subgroup of  $G$, $U$ an $H$-module, $(\s,V)$
a $\g$-stratum.  Then $(\s,V)\rightsquigarrow_\g G*_HU$ iff there
exists $g\in G$ such that $\Ad(g)s\subset\h$ and
$(\Ad(g)\s,V)\rightsquigarrow_\h U$ (the algebra $\Ad(g)\s$ is
represented in $V$ via the isomorphism
$\Ad(g^{-1}):\Ad(g)\s\rightarrow \s$). Conversely, if
$(\s,V)\rightsquigarrow_\h U$, then $(\s,V)\rightsquigarrow_\g
G*_HU$.
\end{Rem}

Now we recall the definition of the Dynkin index (\cite{Dynkin}).
Let $\h$ be a simple subalgebra of $\g$. We fix an invariant
non-degenerate symmetric bilinear form $K_\g$ on $\g$ such that
$K_\g(\alpha^{\vee},\alpha^{\vee})=2$ for a root
$\alpha\in\Delta(\g)$ of the maximal length. Analogously define a
form $K_\h$ on $\h$. The  {\it Dynkin index} of the embedding
$\iota:\h\hookrightarrow \g$  is, by definition,
$K_\g(\iota(x),\iota(x))/K_\h(x,x)$ (the last fraction does not
depend on the choice of $x\in\h$ such that $K_\h(x,x)\neq 0$). For
brevity, we denote the Dynkin index of $\iota$ by $i(\h,\g)$. It
turns out that $i(\h,\g)$ is a positive integer (see~\cite{Dynkin}).

The following lemma seems to be standard.

\begin{Lem}\label{Lem:5.3.3}
Let $\h$ be a simple Lie algebra and $\s$ a subalgebra of $\h$
isomorphic to $\sl_2$. Then  the following conditions are
equivalent:
\begin{enumerate}
\item
$\iota(\s,\h)=1$. \item  $\s\sim_{\Int(\h)}\h^{(\alpha)}$ for a long
root $\alpha\in \Delta(\h)$.
\end{enumerate}
\end{Lem}
\begin{proof}
Clearly, $(2)\Rightarrow (1)$. Let us check $(1)\Rightarrow (2)$.
Choose the standard basis $e,h,f$ in $\s$. We may assume that $h$
lies in the fixed Cartan subalgebra of $\h$. It follows from the
representation theory of $\sl_2$ that $\langle \pi,h \rangle$ is an
integer for any weight $\pi$ of $\h$. Thus $h\in Q^\vee$, where
$Q^\vee$ denotes the dual root lattice. Since
$i(\s,\h)=i(\h^{(\alpha)},\h)$ for a long root $\alpha\in
\Delta(\h)$, the lengths of $h,\alpha^\vee\in Q^\vee$ coincide. It
can be seen directly from the constructions of the root systems,
that all elements  $h\in Q^\vee$ with
$(h,h)=(\alpha^\vee,\alpha^\vee)$ are short dual roots,
see~\cite{Bourbaki}. Thus we may assume that $h=\alpha^\vee$. It
follows from the standard theorems on the conjugacy of
$\sl_2$-triples, see, for example,~\cite{McG}, that $\s$ and
$\h^{(\alpha)}$ are conjugate.
\end{proof}

In Table~\ref{Tbl:5.3.4} we list all simple subalgebras of index 1
in simple classical Lie algebras.

\begin{longtable}{|c|c|}
\caption{Simple subalgebras $\h$ in classical Lie algebras $\g$ with
$\iota(\h,\g)=1$}\label{Tbl:5.3.4}\\\hline
$\g$&$\h$\\\endfirsthead\hline $\g$&$\h$\\\endhead\hline $\sl_n,
n\geqslant 2$&$\sl_k, k\leqslant n$\\\hline $\sl_n, n\geqslant
4$&$\sp_{2k}, 2\leqslant k\leqslant n/2$\\\hline $\so_n,n\geqslant
7$&$\so_k, k\leqslant n, k\neq 4$\\\hline $\so_n,n\geqslant
7$&$\sl_k^{diag}, k\leqslant n/2$\\\hline $\so_n,n\geqslant
8$&$\sp_{2k}^{diag}, 2\leqslant k\leqslant n/4$\\\hline
$\so_n,n\geqslant 7$&$G_2$\\\hline $\so_n, n\geqslant
9$&$\spin_7$\\\hline $\sp_{2n},n\geqslant 2$&$\sp_{2k},k\leqslant
n$\\\hline
\end{longtable}

\begin{Lem}\label{Lem:5.3.5}
Let $\alpha\in \Delta(\g)$ and $\h$ be a reductive subalgebra of
$\g$ containing $\g^{(\alpha)}$. Suppose that there is no proper
ideal of $\h$ containing $\g^{(\alpha)}$. Then
\begin{enumerate}
\item If $\alpha$ is a long root, then $\h$ is simple and
$\iota(\g^{(\alpha)},\h)=\iota(\h,\g)=1$.
\item Suppose $\alpha$ is a short root and $\h$ is not simple. Then
\begin{enumerate}
\item $\g\cong \so_{2l+1},\sp_{2l},l\geqslant 2, F_4$,
$\h=\h_1\oplus\h_2$, where $\h_1,\h_2$ are simple ideals of $\g$
with $i(\h_1,\g)=i(\h_2,\g)=1$. If $\s_i$ denotes the projection of
$\g^{(\alpha)}$ to $\h_i, i=\overline{1,2}$, then $i(\s_i,\h_i)=1$.
\item If  $\g\cong \so_{2l+1},l\geqslant 2$, then $\h=
\so_4$.
\item If $\g\cong \sp_{2l}$, then $\h= \sp_{2k}\oplus\sp_{2m},
k+m\leqslant l$.
\end{enumerate}
\end{enumerate}
\end{Lem}
\begin{proof}
Since $\g^{(\alpha)}\subset [\h,\h]$, we see that $\h$ is
semisimple. Let $\h_1,\h_3$ be simple Lie algebras and $\h_2$ a
semisimple Lie algebra, $\h_1\subset \h_2\subset \h_3$. Let
$\h_2=\h^1_2\oplus\ldots\oplus\h^k_2$ be the decomposition of $\h_2$
into the direct sum of simple ideals and $\iota^i,i=1,\ldots,k,$ the
composition of the embedding $\h_1\hookrightarrow \h_2$ and the
projection $\h_2\rightarrow \h^i_2$. It is shown in \cite{Dynkin}
that
\begin{equation}\label{eq:5.3:1}
i(\g_1,\g_3)=\sum_{i=1}^k i(\iota_i(\g_1),\g_2^i) i(\g_2^i,\g_3).
\end{equation}
This implies assertion 1.

 In assertion
2a it remains to check that $\g\neq G_2$. Indeed, there is a unique
up to conjugacy semisimple but not simple subalgebra of $G_2$,
namely $\g^{(\alpha)}\times \g^{(\beta)}$, where $\beta$ is a long
root and $(\alpha,\beta)=0$. Since $i(\g^{(\alpha)},\g)=1,
i(\g^{(\beta)},\g)=3$, the claim follows from (\ref{eq:5.3:1}).

Proceed to assertion 2b. Note that the representation of
$\g^{(\alpha)}$ in the tautological $\g$-module $V$ is the sum of
the trivial $2l-2$-dimensional and the 3-dimensional irreducible
representations. Since $\g^{(\alpha)}$ is not contained in a proper
ideal of $\h$, we see that the representation of $\h$ in $V/V^{\h}$
is irreducible. Thus $V/V^{\h}=V_1\otimes V_2$, where $V_i$ is an
irreducible $\h_i$-module, $i=1,2$. Note that the representation of
$\g^{(\alpha)}$ in $V_i$ is nontrivial because the projection of
$\g^{(\alpha)}$ to $\h_i$ is nontrivial. From the equality
$\dim(V_1\otimes V_2)/(V_1\otimes V_2)^{\g^{(\alpha)}}=3$ and  the
Clebsh-Gordan formula it follows that  $\dim V_1=\dim V_2=2$. Thence
$\h_1\cong \h_2\cong \sl_2$, $\h=\so_4$.

Proceed to assertion 2c. It follows from assertion 2a that
$i(\h_1,\g)=i(\h_2,\g)=1$. Then $\h_1=\sp_{2m},\h_2=\sp_{2k}$, see
Table~\ref{Tbl:5.3.4}.
\end{proof}

Let us recall the notion of the  {\it index} of a module over a
simple Lie algebra, see~\cite{AEV}. Let $\h$ be a simple Lie
algebra, $U$ an $\h$-module. We define a symmetric invariant
bilinear form $(\cdot,\cdot)_U$ on $\h$ by $(x,y)_U=\tr_U(xy)$. The
form $(\cdot,\cdot)_U$ is nondegenerate whenever $U$ is nontrivial.
By the {\it index} of $U$ we mean the fraction
$\frac{(x,y)_U}{(x,y)_\h}$. Since $\h$ is simple, the last fraction
does not depend on the choice of $x,y\in\h$ with $(x,y)_\h\neq 0$.
We denote the index of $U$ by $l_\h(U)$.

\begin{Lem}\label{Lem:5.3.6}
Let $H$ be a semisimple  algebraic group and $\s$ a subalgebra in
$\h$ isomorphic to $\sl_2$. Let $V$ be an $H$-module such that
$S^{\s}\rightsquigarrow_\h V$.
\begin{enumerate}
\item
Let $e,h,f$ be the standard basis of $\s$. Then
$(h,h)_V=(h,h)_\h-4$.
\item If $\h$ is simple, then $l_\h(V)=1-\frac{4}{i(\s,\h)k_\h}$,
where $k_\h=(\alpha^\vee,\alpha^\vee)_\h$ for a long root
$\alpha\in\Delta(\h)$.
\end{enumerate}
\end{Lem}
\begin{proof}
There is an isomorphism  $V\cong V^\s\oplus \h/\s+(\C^2)^{\oplus 2}$
of $\s$-modules, where $\C^2$ denotes the tautological $\s$-module.
Thus $\tr_V h^2=\tr_\h h^2-\tr_\s h^2+2\tr_{\C^2}h^2=\tr_\h
h^2+8-4$. Assertion 2 stems from $$\frac{\tr_\h
h^2}{\tr_\h\alpha^{\vee
2}}=\frac{i(\s,\h)}{i(\h^{(\alpha)},\h)}=i(\s,\h).$$
\end{proof}

The numbers $k_\h$ for all simple Lie algebras are given in
Table~\ref{Tbl:5.3.7}.

\begin{longtable}{|c|c|c|c|c|c|c|c|c|c|}
\caption{$k_\h$. }\label{Tbl:5.3.7}\\\hline
$\h$&$A_l$&$B_l$&$C_l$&$D_l$&$E_6$&$E_7$&$E_8$&$F_4$&$G_2$\\\hline
$k_\h$& $4l+4$& $8l-4$&$4l+4$&$8l-8$&48&72& 120&36&16\\\hline
\end{longtable}

Let $H$ be a reductive group, $V$  an $H$-module. There exists the
s.g.p. for the action $H:V$, see~\cite{VP}, Theorem 7.2. Recall that
the action $H:V$ is called stable if an orbit in general position is
closed. In this case the s.g.p. is reductive.

\begin{Lem}\label{Lem:5.3.8}
Let $\h$ be a semisimple Lie algebra, $\s$  a subalgebra of $\h_1$
isomorphic to $\sl_2$, $V$ an $\h$-module, $V_1\subset V$ an
$H$-submodule such that the action $H:V_1$ is stable. Let $H_1$
denote  the s.g.p. for the action $H:V_1$. If
$S^{\s}\rightsquigarrow_{\h_1}V/V_1$, then
$S^{\s}\rightsquigarrow_\h V$.
\end{Lem}
\begin{proof}
Let $v_1\in V_1$ be such that $Hv_1$ is closed and $H_{v_1}=H_1$.
The slice module at $v_1$ is the direct sum of $V/V_1$ and a trivial
$H_{v_1}$-module.  The claim of the lemma follows from the Luna
slice theorem and Remark \ref{Rem:5.3.2}.
\end{proof}

\begin{proof}[Proof of Proposition~\ref{Prop:5.3.1}]
The proof is in three steps.

{\it Step 1.}  Here we suppose that  $\h$ is simple. An admissible
triple $(\g,\h,V)$ is $W$-quasiessential iff
$S^{\s}\rightsquigarrow_\h V\oplus V^*\oplus \g/\h$, where $\s$ is a
subalgebra in $\h$ such that $\s\sim_G\g^{(\alpha)}$ for some
$\alpha\in \Delta(\g) $.

Suppose  $S^{\s}\rightsquigarrow_\h V\oplus V^*\oplus \g/\h$ for
some subalgebra $\s\subset \h$ isomorphic to $\sl_2$. From assertion
2 of Lemma~\ref{Lem:5.3.6} it follows that
\begin{equation}\label{eq:5.3:2}l_\h(V\oplus V^*\oplus
\g/\h)=1-\frac{4}{i(\s,\h)k_\h}.\end{equation} Thence
$l_\h(\g/\h)<1$. In Section 3 of~\cite{ranks} it was shown that
$i(\h,\g)=1$. Equivalently, $i(\s,\h)=i(\s,\g)$. All simple
subalgebras $\h\subset \g$ with $l_\h(\g/\h)<1$ are listed
in~\cite{ranks}, Table 5. The list (up to  $\Aut(\g)$-conjugacy) is
presented in Table~\ref{Tbl:5.3.9}. In column 4 the nontrivial part
of the representation of $\h$ in $\g/\h$ is given. By  $\tau$ we
denote the tautological representation of a classical Lie algebra.
\begin{longtable}{|c|c|c|c|}
\caption{Simple subalgebras $\h\subsetneqq\g$ with $l_\h(\g/\h)<
1$}\label{Tbl:5.3.9}\\\hline
N&$\g$&$\h$&$\g/\h_+$\\\endfirsthead\hline
N&$\g$&$\h$&$\g/\h_+$\\\endhead\hline 1&$\sl_n, n>1$&$\sl_k, n/2<
k<n$&$(n-k)(\tau+\tau^*)$\\\hline 2&$\sl_{2n},n\geqslant
2$&$\sp_{2n}$&$\bigwedge^2\tau$\\\hline 3&$\sl_{2n+1},n\geqslant
2$&$\sp_{2n}$&$2\tau+\bigwedge^2\tau$\\\hline 4&$\sp_{2n},n\geqslant
2$&$\sp_{2k}, n/2\leqslant k< n$&$2(n-k)\tau$\\\hline 5&$\so_n,
n\geqslant 7$&$\so_k, \frac{n+2}{2}< k< n, k\neq
4$&$(n-k)\tau$\\\hline 6&$\so_{2n},n\geqslant
5$&$\sl_n$&$\bigwedge^2\tau+\bigwedge^2\tau^*$\\\hline
7&$\so_{2n+1},n\geqslant 3
$&$\sl_n$&$\tau+\tau^*+\bigwedge^2\tau+\bigwedge^2\tau^*$\\\hline
8&$\so_n, 9\leqslant n\leqslant
11$&$\mathfrak{spin}_7$&$\tau+(n-8)R(\pi_3)$\\\hline 9&$\so_n,
7\leqslant n\leqslant 9$&$G_2$&$(n-3)R(\pi_1)$\\\hline
10&$G_2$&$A_2$&$\tau+\tau^*$\\\hline
11&$F_4$&$B_4$&$R(\pi_4)$\\\hline
12&$F_4$&$D_4$&$\tau+R(\pi_3)+R(\pi_4)$\\\hline
13&$F_4$&$B_3$&$2\tau+2R(\pi_3)$\\\hline
14&$E_6$&$F_4$&$R(\pi_1)$\\\hline
15&$E_6$&$D_5$&$R(\pi_4)+R(\pi_5)$\\\hline
16&$E_6$&$B_4$&$\tau+2R(\pi_4)$\\\hline
17&$E_7$&$E_6$&$R(\pi_1)+R(\pi_5)$\\\hline
18&$E_7$&$D_6$&$2R(\pi_6)$\\\hline
19&$E_8$&$E_7$&$2R(\pi_1)$\\\hline
\end{longtable}

All orthogonal $H$-modules $U$ such that $m_H(U)=\dim H$ and
 $S^{\s}\rightsquigarrow_\h U$ were found by G. Schwarz
in~\cite{Schwarz3}, Tables I-V ("orthogonal representations with an
$\mathbb{S}^3$-stratum", in his terminology). This modules are given
in Table~\ref{Tbl:5.3.10}. Note that one of them is omitted in
\cite{Schwarz3}, namely N11, $k=1$. To see that
$S^{\s}\rightsquigarrow U$, where $\s=\h^{(\alpha)}$ and $\alpha$ is
a long root in $\Delta(\h)$, one applies Lemma~\ref{Lem:5.3.8} to
$\C^{12}\subset U$. The corresponding pair $(\h_1,V/(V_1\oplus
V^{\h_1}))$ is N6 of Table~\ref{Tbl:5.3.10}.

\begin{longtable}{|c|c|c|}
\caption{Orthogonal $H$-modules $U$ with $m_H(U)=\dim H$  and
$S^\s\rightsquigarrow_\h U$}\label{Tbl:5.3.10}\\\hline
N&$\h$&$U/U^\h$\\\endfirsthead\hline N&$\h$&$U$\\\endhead\hline 1&
$\sl_n, n>1,$&$(n-1)(\tau+\tau^*)$\\\hline 2& $\sl_n,
n>3,$&$\tau+\tau^*+\bigwedge^2\tau+\bigwedge^2\tau^*$\\\hline 3&
$\sl_4$& $2(\tau+\tau^*)+\bigwedge^2\tau$\\\hline 4& $\so_7$&
$(4-k)\tau+kR(\pi_3), k>0$\\\hline 5& $\so_9$&
$kR(\pi_4)+(6-2k)\tau,k>0$\\\hline 6& $\so_{11}$& $2
R(\pi_5)$\\\hline 7& $\sp_{2m},m>1$& $2m\tau$\\\hline 8&
$\sp_{2m},m>1$& $2\tau+R(\pi_2)$\\\hline 9& $\so_8$&$k\tau+ l
R(\pi_3)+m R(\pi_4), k+l+m=5, k,l,m<5$\\\hline 10& $\so_{10}$&
$3\tau+ R(\pi_4)+R(\pi_5)$\\\hline 11&$\so_{12}$& $\tau+k
R(\pi_5)+(2-k)R(\pi_6)$\\\hline 12&$G_2$&$3R(\pi_1)$\\\hline
\end{longtable}

It follows from (\ref{eq:5.3:2}) that $\iota(\s,\h)=1$ for all
modules from Table~\ref{Tbl:5.3.10}. This proves assertion 2 of the
proposition for simple $\h$. Inspecting
Tables~\ref{Tbl:5.3.9},\ref{Tbl:5.3.10}, we find all admissible
triples $(\g,\h,V)$ with $S^{\s}\rightsquigarrow_\h V\oplus
V^*\oplus \g/\h$. This completes the proof of assertion 1 when $\h$
is simple.

{\it Step 2.} It remains to consider the situation when  $\h$ is not
simple. On this step we assume that $\h$ possesses the following
property
\begin{itemize}
\item[(*)] $S^{\s}\rightsquigarrow_\h\g/\h\oplus V\oplus V^*$,
where $\s\subset\h$ is such that
\end{itemize}
\begin{enumerate}
\item
$\s\subset \h_1\oplus\h_2$, where $\h_1,\h_2$ are simple ideals of
$\h$.
\item  $\s\sim_G\g^{(\alpha)}$, where $\alpha$ is a short root of
$\g$. \item $i(\s_i,\h_i)=1$, where $\s_i$ denotes the projection of
$\s$ to $\h_i, i=1,2$.
\end{enumerate}

It follows from Lemma~\ref{Lem:5.3.5} that a subalgebra $\s\subset
\h$ satisfying (1)-(3) is defined uniquely up to
$\Int(\h)$-conjugacy.

Let us check that
\begin{equation}\label{eq:5.3:3}
k_\g= k_{\h_1}+k_{\h_2}-2-k_{\h_1}l_{\h_1}(V)-k_{\h_2}l_{\h_2}(V).
\end{equation}

Let $h\in\s$ be a dual root, $h=h_1+h_2, h_i\in\h_i$. By
Lemma~\ref{Lem:5.3.6},
$(h,h)_{V}+(h,h)_{V^*}+(h,h)_{\g/\h}=(h,h)_{\h}-4$. Equivalently,
\begin{equation}\label{eq:5.3:4}2(h,h)_V+(h,h)_\g=2(h,h)_\h-4.\end{equation}
Since $(\cdot,\cdot)_V,(\cdot,\cdot)_\h$ are $\h$-invariant forms on
$\h$, we see that
 $\h_1,\h_2$ are orthogonal with respect to
$(\cdot,\cdot)_{V},(\cdot,\cdot)_\h$. Therefore  (\ref{eq:5.3:4}) is
equivalent to
\begin{equation}\label{eq:5.3:5}
(h,h)_{\g}+2((h_1,h_1)_V+(h_2,h_2)_V)=2((h_1,h_1)_{\h_1}+(h_2,h_2)_{\h_2})-4.
\end{equation}
By assertion  2a of Lemma~\ref{Lem:5.3.5} and Lemma~\ref{Lem:5.3.3},
 $h_i\sim_{\Int(\h_i)}\alpha^\vee$, where $\alpha$ is a long root in $\Delta(\h_i) $. Therefore
 $(h_i,h_i)_{\h_1}=k_{\h_i},i=1,2$. It follows  from the choice of $\s$ that $i(\s,\h)=2$ whence
$(h,h)_{\g}=2k_\g$. So (\ref{eq:5.3:5}) and (\ref{eq:5.3:3}) are
equivalent.

Let us show that  $\g\not\cong F_4,\so_{2l+1},l\geqslant 3$. Assume
that $\g\cong F_4$. By (\ref{eq:5.3:3}), $k_{\h_1}+k_{\h_2}\geqslant
38$. Since
 $\rank \h_1+\rank \h_2\leqslant 4$, this is impossible (see Table~\ref{Tbl:5.3.7}). Assume that $\g\cong \so_{2l+1},l\geqslant
3$. By assertion 2b of Lemma~\ref{Lem:5.3.5}, $\h_1\cong \h_2\cong
\sl_2$. This again contradicts (\ref{eq:5.3:3}).

It follows from assertions 2a,2c of Lemma~\ref{Lem:5.3.5} that
$\g\cong\sp_{2n}$, $\h_1\oplus \h_2=\sp_{2m_1}\oplus
\sp_{2m_2},m_1+m_2\leqslant n,m_1\leqslant m_2$. One can rewrite
(\ref{eq:5.3:3}) as
\begin{equation}\label{eq:5.3:6}
2(m_1+m_2-n)+1=(2m_1+2)l_{\h_1}(V)+(2m_2+2)l_{\h_2}(V).
\end{equation}
From (\ref{eq:5.3:6}) it follows that $m_1+m_2=n$,
$l_{\h_i}(V)\leqslant \frac{1}{2m_i+2}$. Therefore the $\h$-module
$\g/\h$ is the tensor product of the tautological $\sp_{2m_1}$- and
$\sp_{2m_2}$-modules and, see the table from \cite{AEV},  $V$ is the
tautological $\sp_{2m_j}$-module for some $j\in \{1,2\}$. Recall
that $m_H(V\oplus V^*\oplus \g/\h)=\dim H$. In particular,
$m_{H_i}(V\oplus V^*\oplus \g/\h)=\dim H_i,i=1,2$. By above, the
$H_i$-module $V\oplus V^*\oplus \g/\h$ is the direct sum of
tautological $\Sp(2m_i)$-modules. Thus if $m_1<m_2$, then $H_1$ acts
trivially on $V$ and $m_2=m_1+1$. This shows that a
$W$-quasiessential triple $(\g,\h,V)$ satisfying (*) is one of NN
32,33 of Table \ref{Tbl:5.0.3}.

Let us show that the triples NN32,33 satisfy (*). Put $V_0=\g/\h$.
This is an orthogonal whence stable, see, for instance,
\cite{Lun_ort}, $H$-module. The s.s.g.p. $\h_0$  for the $H$-module
$V_0$ is the direct sum of  $m_1$ copies of $\sl_2$ embedded
diagonally into $\h_1\oplus\h_2$ and $m_2-m_1$ copies of $\sl_2$
embedded into $\h_2$ (see \cite{Elash2}). We are done by Lemma
\ref{Lem:5.3.8} applied to $\g/\h\subset \g/\h\oplus V\oplus V^*$.

Finally, we see that
$S^{(\h_2^{(\alpha)})}\rightsquigarrow_{\h_2}V\oplus V^*\oplus
\g/\h$  (if $\h_1\cong\h_2$, then for $\h_2$ we take the ideal of
$\h$ acting on $V$ trivially), where $\alpha$ is a long root of
$\h_2$.

{\it Step  3.} It remains to show that any $W$-quasiessential triple
$(\g,\h,V)$ satisfies (*) provided $\h$ is not simple. Assume the
converse. Let us show that
\begin{itemize}\item[(**)] $S^{\s}\rightsquigarrow_\h \g/\h\oplus V\oplus
V^*$ for some subalgebra $\s\subset \h,\s\cong \sl_2,$ not lying in
any simple ideal of $\h$.\end{itemize} Indeed, let  $\h_1$ be a
simple ideal  of $\h$ and $\s_1$ a subalgebra of $\h_1$ isomorphic
to $\sl_2$ such that $S^{\s_1}\rightsquigarrow_{\h}\g/\h\oplus
V\oplus V^*$, or equivalently,
$S^{\s_1}\rightsquigarrow_{\h_1}\g/\h\oplus V\oplus V^*$. Then
$(\g,\h_i,V/V^{\h_i})$ is one of the triples NN 1-31 from
Table~\ref{Tbl:5.0.3} and $\s_1\sim_G\g^{(\alpha)}$, where $\alpha$
is a long root in $\Delta(\g)$. If (**) does not hold, then, by step
1, $S^{(\alpha)}\rightsquigarrow_\g T^*(G*_HV)$ implies that
$\alpha$ is a long root.  Since $\h$ is not simple, we see that
$(\g,\h,V)$ is not $W$-quasiessential. So (**) is checked.

It follows from assertion  2a of Lemma~\ref{Lem:5.3.6} that there is
a subalgebra $\s\subset \h$ and simple ideals $\h_1,\h_2$ of $\h$
such that $(\g,\h_1\oplus\h_2,V/V^{\h_1\oplus \h_2})$ satisfies
condition (*) of step 2. But in this case $\n_\g(\h_1\oplus
\h_2)=\h_1\oplus\h_2$ whence $\h=\h_1\oplus \h_2$. Contradiction.
\end{proof}

\subsection{Computation of Weyl groups}\label{SUBSECTION_Weyl_aff3}
In this subsection we compute the groups $W(\g,\h,V)$ for triples
$(\g,\h,V)$ listed in Table~\ref{Tbl:5.0.3}.

At first, we reduce the problem to the case when $\C[V]^H=\C$. Let
$(\g,\h,V)$ be one of the triples from Table~\ref{Tbl:5.0.3}. Let
$H_0$ denote the unit component  of the principal isotropy subgroup
for the action  $H:V$. Put $V_0=V/(V^{H_0}+\h v)$, where $v\in V$ is
such that
 $Hv$ is closed and $H_v^\circ=H_0$. Thanks to Corollary
 \ref{Cor:1.2.3},
$W(\g,\h,V)=W(\g,\h_0,V_0)$. The triple $(\g,\h_0,V_0)$ is said to
be {\it reduced} from $(\g,\h,V)$. Clearly,
$(\g,\h,V)=(\g,\h_0,V_0)$ iff $\C[V]^H=\C$. Since  $H$ is
semisimple, we get $\Quot(\C[V]^H)=\C(V)^H$. So $\C[V]^H=\C$ iff $H$
has a dense orbit in $V$. In Table \ref{Tbl:5.4.1} we present all
triples that are reduced from nonreduced triples from Table
\ref{Tbl:5.0.3}. We use the same notation as in
Table~\ref{Tbl:5.0.3}. In column 2 we give the number of $(\g,\h,V)$
in Table~\ref{Tbl:5.0.3} and, in some cases,   restrictions on  the
algebra $\g$ or the  $\h$-module $V$. It turns out that
$(\g,\h_0,V_0)$ is again contained in Table~\ref{Tbl:5.0.3}.

\begin{longtable}{|c|c|c|c|}\caption{Reduced triples}
\label{Tbl:5.4.1}\\\hline
N&$(\g,\h,V)$&$\h_0$&$V_0$\\\endfirsthead\hline
N&$(\g,\h,V)$&$\h_0$&$V_0$\\\endhead\hline 1&N1,
$m>0$&$\sl_{n-k-m}$&$(l-m)\tau$\\\hline 2&N2, even
$n$&$\sp_{n}$&$\tau$\\\hline 3&N2, odd $n$&$\sp_{n-1}$&$0$\\\hline
4&N3&$\sp_n$&$\tau$\\\hline 5&N6, odd $n$&$\sp_{n-1}$&$0$\\\hline
6&N10&$\sl_3$&$0$\\\hline 7&N11, $k=1$&$\sl_3^{diag}$&$0$\\\hline
8&N12&$\sl_3$&$0$\\\hline 9&N13&$G_2$&0\\\hline 10&N14,
$(k,l)=(2,0)$&$\sl_4$&$0$\\\hline 11&N14, $(k,l)\neq
(2,0)$&$G_2$&$0$\\\hline 12&N15&$G_2$&$0$\\\hline
13&N16&$\sl_4^{diag}$&0\\\hline 14&N17&$G_2$&0\\\hline
15&N19&$\sl_5^{diag}$&$0$\\\hline
16&N20&$\mathfrak{spin}_7$&$0$\\\hline
17&N21&$\mathfrak{spin}_7$&$0$\\\hline
18&N22&$\mathfrak{spin}_7$&$0$\\\hline 19&N25, even
$n$&$\sp_n$&$0$\\\hline 20&N25, odd $n$&$\sp_{n+1}$&$\tau$\\\hline
21&N30&$B_3$&$0$\\\hline 22&N31&$B_3$&$0$\\\hline
\end{longtable}

\begin{Lem}\label{Lem:5.4.2}
All triples $(\g,\h,V)$ from Table~\ref{Tbl:5.0.3} such that
$\C[V]^H\neq \C$ are listed in the second column of Table
\ref{Tbl:5.4.1}. The reduced triple for $(\g,\h,V)$ coincides with
$(\g,\h_0,V_0)$.
\end{Lem}
\begin{proof}
The triples NN 32,33 are reduced. So we may assume that $\h$ is
simple.  The list of all simple linear groups with a dense orbit is
well-known, see, for example, \cite{Vinberg_connections}. It follows
from the classification of the paper~\cite{Elash1} that the
s.s.g.p's for the $H$-module $\g/\h\oplus V\oplus V^*$ are simple
for all triples $(\g,\h,V)$, except NN1,2,4,20. By Popov's
criterion, see~\cite{Popov}, the action $H:V$ are stable whenever
the s.s.g.p is reductive. For the remaining four triples the reduced
triples are easily found case by case.
\end{proof}

Below in this subsection we consider only reduced triples
$(\g,\h,V)$ from Table~\ref{Tbl:5.0.3}.

The group $T_1\times T_2$, where $T_1=Z(Z_G(H))^\circ, T_2=
Z(\GL(V)^H)$, is a torus naturally acting on $X=G*_HV$ by
$G$-automorphisms. Namely, we define the action morphism by
$(t_1,t_2,[g,v])\mapsto[gt_1^{-1},t_2v], t_1\in T_1,t_2\in T_2, g\in
G,v\in V$. Put $\widetilde{G}=G\times T_1\times T_2$,
$\widetilde{H}=H\times T_1\times T_2$.

For some triples $(\g,\h,V)$ the inequality
$\rank_{\widetilde{G}}(X)<\rank \widetilde{G}$ holds. These triples
are presented in Table~\ref{Tbl:5.4.3}. The matrix in column 3 of
the first row is the diagonal matrix, whose  $k+1$-th entry is
$(n-1)x$ and the other entries are $-x$ .

\begin{longtable}{|c|c|c|}\caption{The projections of $\l_{0\,\widetilde{G},X}$ to $\g$}\label{Tbl:5.4.3}
\\\hline N&$(\g,\h,V)$& the projection of $\l_{0\,\widetilde{G},X}$ to $\g$\\\endfirsthead\hline
N&$(\g,\h,V)$&the projection of $\l_{0\,\widetilde{G},X}$ to
$\g$\\\endhead\hline 1&N1,
$m=0$&$diag(-x^k,(n-1)x,-x^{n-k-1})$\\\hline
2&N4&$diag(-\frac{x}{n+1},-\frac{x}{n+1},\frac{x}{n-1},-\frac{x}{n+1}\ldots,\frac{x}{n-1},-\frac{x}{n+1})$\\\hline
3&N5&$diag(\frac{x}{n},\frac{x}{n}
-\frac{x}{(n-2)},\frac{x}{n},\ldots,
-\frac{x}{(n-2)},\frac{x}{n})$\\\hline
4&N7&$diag(-x,x,-x,\ldots,x)$\\\hline
5&N8&$diag(\frac{x}{n+1},-\frac{x}{n-1},\frac{x}{n+1},\ldots,-\frac{x}{n+1},\frac{x}{n+1})$\\\hline
6&N26&$x\sum_{i=1}^{2n+1}(-1)^{i}\varepsilon_i$\\\hline
7&N27&$x\sum_{i=1}^{2n}(-1)^{i}\varepsilon_i$\\\hline
8&N32&$x\sum_{i=1}^{2n}(-1)^{i}\varepsilon_i$\\\hline 9&N33&
$x\sum_{i=1}^{2n+1}(-1)^{i}\varepsilon_i$
\\\hline
\end{longtable}

\begin{Lem}\label{Lem:5.4.4}
If $(\g,\h,V)$ is a triple from the first column of Table
\ref{Tbl:5.4.3}, then $W(\g,\h,V)$ coincides with the group
indicated in Table~\ref{Tbl:5.0.3}.
\end{Lem}
\begin{proof}
At first, we will check that the projections of
$\l_{0\,\widetilde{G},\widetilde{X}_0}$ to $\g$ for triples in
consideration are given in column 3 of Table \ref{Tbl:5.4.3}.

Let us note that  $\dim T_1=1$ for rows NN 1 (with $k\neq 0$),3,5-8
of Table \ref{Tbl:5.4.3}, otherwise $\dim T_1=0$. The dimension of
$T_2$ equals 2 for the second row, 1 for rows NN 1,3,4,6-9, and 0
otherwise.

We embed $\widetilde{H}$ into $\widetilde{G}$ via
$(h,t_1,t_2)\mapsto (ht_1,t_1,t_2)$. Equip $V$ with a natural
structure of an $\widetilde{H}$-module: $H$ acts on $V$ as before,
$T_1$ trivially, and $T_2$ via the identification $T_2\cong
Z(\GL(V)^H)$. The $\widetilde{G}$-varieties $X,
\widetilde{G}*_{\widetilde{H}}V$ are isomorphic.

Since $\a_{G,X}=\t$, the projection
$\a_{\widetilde{G},X}\rightarrow\t$  is surjective.  Thus
$\a_{\widetilde{G},X}, \l_{0\,\widetilde{G},X}$ are mutually
orthogonal and $\a_{\widetilde{G},X}\oplus
\l_{0\,\widetilde{G},X}=\t\oplus\t_1\oplus \t_2$. To prove the claim
it is enough to compute $\a_{\widetilde{G},X}\subset \t\oplus
\t_1\oplus\t_2$. In virtue of the isomorphism $X\cong
\widetilde{G}*_{\widetilde{H}}V$, the computation can be done by
using the algorithm from Subsection~\ref{SUBSECTION_algorithm1}.
However, in some cases one can simplify the computation. Namely, for
rows NN1-3,6,7 there is an antistandard parabolic $Q\subset
\widetilde{G}$ with the standard Levi subgroup $M\subset Q$ such
that $H=(M,M)$. Applying Proposition~\ref{Prop:1.5.4}, we reduce the
computation of $\a_{\widetilde{G},X}$ to computing the spaces
$\a_{\bullet,\bullet}$ for certain linear actions.  For rows 2,3,6,7
these linear actions are spherical, so the Cartan spaces can be
extracted, for example, from the second table in~\cite{Leahy},
Section 2.

Now we proceed to the claim  on the Weyl groups. The projection
$\a_0$ of $\l_{0\,\widetilde{G},X}$ to $\g$ lies in $\t\cap
(\a_{\widetilde{G},X}\cap\g)^\perp$. It follows from Proposition
\ref{Prop:3.4.3} that $\a_0\subset \t^{W_{G,X}}$. Therefore
$W_{G,X}\subset W(\g)_{\xi}$ for $\xi\in \a_0\setminus \{0\}$. By
Proposition \ref{Prop:1.2.11}, $W_{G,X}$ is one of the groups listed
in Table \ref{Tbl:5.2.8}. Now the equalities $W_{G,X}=W(\g)_{\xi}$
are easily checked case by case.
\end{proof}

It remains to compute the Weyl groups for  triples NN 9,11 (with
$k=2$), 18, 23-25,28,31 from Table \ref{Tbl:5.0.3}.

Note that $s_{\alpha}\in W(\g,\h,V)$ for a short root $\alpha\in
\Delta(\g) $. Indeed, Proposition~\ref{Prop:5.3.1} implies that
$S^{(\alpha)}\rightsquigarrow_\g T^*(G*_HV)$ iff $\alpha\in
\Delta(\g) $ is a long root. Our claim follows from
Corollary~\ref{Cor:5.2.3}.

{\it N9.} Here $\g=\so_{2n+1},\h=\sl_{n}, V=\{0\}$. Put
$\m:=\t+\g^{(\alpha_1,\ldots,\alpha_{n-1})}$ (so that $\h=[\m,\m]$),
$\q=\m+\b^-$ and let $M,Q$ be the corresponding Levi and parabolic
subgroups. The representation of $H\cong \SL_{n}$ in $\Rad_u(\q)$
coincides with $\bigwedge^2 \tau^*+\tau^*$. From
Corollary~\ref{Cor:3.4.3} it follows that $W(\g,\h,V)\cap
M/T=W(\sl_n,\sl_n,\bigwedge^2\C^{n*}\oplus \C^{n*})$. By above, the
last Weyl group is generated by
$s_{\varepsilon_i-\varepsilon_{i+2}}, i=\overline{1,n-2}$. Thence
$W(\g,\h,V)$ coincides with the group indicated in
Table~\ref{Tbl:5.0.3}.

{\it N11.} Here $\g=\so_7,\h=\so_6$ and $V$ is the direct sum of two
copies of the semispinor $\h$-module $V(\pi_3)$. We assume that $\h$
is embedded into $\so_7$ as the annihilator of $e_4$. It follows
from assertion  3 of Proposition~\ref{Prop:1.2.1} that
$W(\g,\h,V)=W(\g,\h_0)$, where  $\h_0$ denotes the s.s.g.p. of the
$\h$-module $V$. 
%
For $\h_0$ one may take the following subalgebra
$$\h_0=\{\begin{pmatrix}0&0&0&0&0&0&0\\x&a&b&0&0&0&0\\
y&c&-a&0&0&0&0\\ 0&0&0&0&0&0&0\\ -t&0&0&0&a&-b&0
\\-z&0&0&0&-c&-a&0\\ 0&z&t&0&-y&-x&0\end{pmatrix}, x,y,z,t,a,b,c\in \C\}.$$

Put $\m:=\t+\g^{(\varepsilon_2-\varepsilon_3)},\q=\b^-+\m$.  It is
seen directly that $\s:=\h_0\cap \m\cong \sl_2$ is a Levi subalgebra
in $\h_0$ and  $\Rad_u(\h_0)\subset \Rad_u(\q)$. The nontrivial part
of the $\s$-module $\q/\h_0$ is two-dimensional. Applying Corollary
\ref{Cor:3.4.3} to $Q$ and $M$, we see that
$s_{\varepsilon_2-\varepsilon_3}\not\in W(\g,\h_0)$.

It remains to show that $s_{\varepsilon_1-\varepsilon_2}\in
W(\g,\h_0)$. Put
$$\widetilde{\h}_0=\{\begin{pmatrix}u&0&0&0&0&0&0\\x&a&b&0&0&0&0\\y&c&d&0&0&0&0\\0&0&0&0&0&0&0\\
-t&v&0&0&-d&-b&0\\-z&0&-v&0&-c&-a&0\\
0&z&t&0&-y&-x&-u\end{pmatrix}, a,b,c,d,u,v,x,y,z,t\in \C\}.$$

By \cite{Wasserman}, Table B, row 7,
$\a(\g,\widetilde{\h}_0)=\langle \varepsilon_1,\varepsilon_2\rangle$
and $W(\g,\widetilde{\h}_0)$ is generated by
$s_{\varepsilon_1-\varepsilon_2},s_{\varepsilon_2}$. From assertion
1 of Proposition~\ref{Prop:1.2.1} it follows that
$W(\g,\widetilde{\h}_0)$ is a subquotient in $W(\g,\h_0)$. Thus
$s_{\varepsilon_1-\varepsilon_2}\in W(\g,\h_0)$.

{\it N18.} Here $\g=\so_9,\h=G_2, V=\{0\}$.

We may assume that $H\subset M:=TG^{(\alpha_2,\alpha_3,\alpha_4)}$.
Put $\q:=\b^{-}+\m$. Applying Corollary~\ref{Cor:3.4.3} to $Q,M$, we
see that $W(\g,\h)\cap M/T=W(\so_7,\h,V(\pi_3))$. The last group has
been already computed. It is generated by
$s_{\varepsilon_2-\varepsilon_4},s_{\varepsilon_3},s_{\varepsilon_4}$.
Therefore $W(\g,\h)$ is generated either by
$s_{\varepsilon_1-\varepsilon_2},s_{\varepsilon_2-\varepsilon_4},s_{\varepsilon_3},s_{\varepsilon_4}$,
or by $s_{\varepsilon_1-\varepsilon_3},s_{\varepsilon_3},$
$s_{\varepsilon_2-\varepsilon_4},s_{\varepsilon_4}$. According to
Corollary~\ref{Cor:5.2.3}, it remains to check
$S^{(A)}\not\rightsquigarrow_\g T^*(G/H)$, where
$A=\{\varepsilon_1-\varepsilon_2,\varepsilon_3-\varepsilon_4\}$.
Indeed, there are no $G^{(A)}$-fixed point in $T^*(G/H)$ because
$G^{(A)}\subset G$ is not conjugate to a subgroup in $H=G_2$.

{\it N23.} Here $\g=\so_{11},\h=\mathfrak{spin}_7, V=\{0\}$.

As above, we may assume that $H\subset
M:=TG^{(\alpha_2,\ldots,\alpha_5)}$. Analogously to the previous
case, we obtain $W(\g,\h)\cap M/T=W(\so_9,\h,V(\pi_3))$. The last
Weyl group is generated
$s_{\varepsilon_2-\varepsilon_3},s_{\varepsilon_3-\varepsilon_5},
s_{\varepsilon_4},s_{\varepsilon_5}$. As above, one should check
that $S^{(A)}\not\rightsquigarrow_\g T^*(G/H)$, where
$A=\{\varepsilon_1-\varepsilon_2,\varepsilon_3-\varepsilon_4\}$.
Assume the converse. If  $\g^{(A)}\sim\s$ for some subalgebra
$\s\subset\h$, then  both simple ideals of $\s$ are of index 1 in
$\h$. Such a subalgebra $\s\subset\h\cong\so_7$ is unique, it
coincides with $\so_4\subset\so_7$. By Remark~\ref{Rem:5.3.2}, $(\s,
V_0)\rightsquigarrow_\h U:=\g/\h$, where $\dim V_0=8$. Therefore
$\codim_U (\h v+U^{\s})=8$ for some $v\in U^{\s}$. From this
observation we deduce that $\dim U^{\s}+\dim \h/\n_\h(\s)\geqslant
\dim U-8$. There is an isomorphism of $\h$-modules $U\cong
V(\pi_1)\oplus V(\pi_3)^{\oplus 3}$. It can be easily seen that
$\dim V(\pi_1)^{\s}=3, \dim V(\pi_3)^{\s}=0, \dim
\h/\n_\h(\s)=21-9=12$. Contradiction.

{\it N24.} Here $\g=\so_{13},\h=\so_{10}, V=V(\pi_4)$.

Again, we may assume that $H\subset (M,M)$, where
$M:=TG^{\alpha_2,\ldots,\alpha_6}$. Using Corollary~\ref{Cor:3.4.3},
we have $W(\g,\h,V)\cap M/T=W([\m,\m],\h,V(\pi_1)\oplus V(\pi_4))$.
As above, it remains to show that $S^{(A)}\not\rightsquigarrow_\g
G*_H(\h^\perp\oplus V\oplus V^*)$, where
$A=\{\varepsilon_1-\varepsilon_2,\varepsilon_4-\varepsilon_5\}$.
There is a unique (up to $H$-conjugacy) subalgebra $\s\subset \h$
that is $G$-conjugate with $\g^{(A)}$, namely, $\s=\h^{(A')}$, where
$(A')=\{\varepsilon_1-\varepsilon_2,\varepsilon_3-\varepsilon_4\}$
(the roots are taken in the root system $\Delta(\so_{10})$).

Analogously to the previous case, one should prove that $\dim
\widetilde{V}^{\s}+\dim \h/\n_\h(\s)<\dim \widetilde{V}-8$, where
$\widetilde{V}:=\h^\perp\oplus V\oplus V^*\cong V(\pi_1)^{\oplus
3}\oplus V(\pi_4)\oplus V(\pi_5)$. It can be easily seen that $\dim
V(\pi_1)^\s=2, \dim \h/\n_{\h}(\s)=32$. To compute $\dim
(V(\pi_4)\oplus V(\pi_5))^\s$ we note that the weight system of the
$\h$-module $V(\pi_4)\oplus V(\pi_5)$ consists of elements
$\frac{1}{2}\sum_{i=1}^5\pm\varepsilon_i$ without multiplicities.
The subspace $(V(\pi_4)\oplus V(\pi_5))^\s$ is the direct sum of
weight subspaces. A weight vector $v_\lambda$ of weight  $\lambda$
is annihilated by $\s$ iff
$(\lambda,\varepsilon_1-\varepsilon_2)=(\lambda,\varepsilon_3-\varepsilon_4)=0$.
Therefore $\dim (V(\pi_4)\oplus V(\pi_5))^\s=8$ whence the required
inequality.

{\it N25.} Here $\g=\sp_{2n},\h=\sp_{2m},
m=\left[\frac{n+1}{2}\right]$, $V=V(\pi_1)$ for odd $n$ and 0
otherwise.

It is enough to check $W(\g,\h,V)\neq W(\g)$. We do this by
induction on $n$. For  $n=1$ we get $\g=\h=\sl_2,V=V(\pi_1)$ and
$W(\g,\h,V)=\{1\}$. Now let $n>1$ and set
 $M:=TG^{\alpha_2,\ldots,\alpha_n}$. One may assume that $H\subset
(M,M)$. Using Proposition~\ref{Cor:3.4.3}, we see that
$W(\g,\h,V)\cap M/T=W([\m,\m],\h, V(\pi_1)\oplus V)$. The reduced
triple for $([\m,\m],\h,V(\pi_1)\oplus V)$ coincides with
$(\sp_{2(n-1)},\sp_{2[n/2]},V(\pi_1)^{\oplus (2\{n/2\})})$. We are
done by the inductive assumption.

{\it N28.} Here  $\g=G_2, \h=A_2, V=V(\pi_1)$.

Again, it is enough to show that  $W(\g,\h,V)\neq W(\g)$. Note that
$\alpha_2,\alpha_2+3\alpha_1$ is a simple root system in $\h$. The
subalgebra $\h_0=\g^{(\alpha_2)}\oplus
\g^{-\alpha_2-3\alpha_1}\oplus\g^{-2\alpha_2-3\alpha_1}$ is the
s.s.g.p. for the $\h$-module $V$. Note that   $\h_0$ is tamely
embedded into $\q:=\g^{\alpha_2}\oplus\b^{-}$. Set
$M:=TG^{(\alpha_2)}$. Note that the restriction of the $\h_0$-module
$\Rad_u(\q)/\Rad_u(\h_0)$ to $\g^{(\alpha_2)}$ is the irreducible
two-dimensional module. Applying Corollary \ref{Cor:3.4.3} to the
homogeneous space $G/H_0$ and the pair $(Q,M)$, we see that
$W(\g,\h_0)\cap M/T=\{1\}$. But $W(\g,\h,V)=W(\g,\g/\h_0)$ by
assertion 3 of Proposition~\ref{Prop:1.2.1}.

{\it N31.} Here $\g=F_4, \h=B_3, V=\{0\}$.

Set $M:=TG^{\alpha_2,\alpha_3,\alpha_4}$. We may assume that
$\h=[\m,\m]$. Using Corollary~\ref{Cor:3.4.3}, we see that
$W(\g,\h,V)\cap M/T=W([\m,\m],\h,V(\pi_3)\oplus V(\pi_1))$. As we
have shown above, the latter is generated by
$s_{\varepsilon_2-\varepsilon_4},s_{\varepsilon_3},s_{\varepsilon_4}$.
In particular, $W(\g,\h,V)\neq W(\g)$. Recall that $s_{\alpha}\in
W(\g,\h,V)$ for all simple roots   $\alpha$. The reflections
$s_{\alpha}$, where $\alpha$ is a short root, and
$s_{\varepsilon_2-\varepsilon_4}$ generate the subgroup of $W(\g)$
indicated in Table~\ref{Tbl:5.0.3}. This subgroup is maximal.

\subsection{Completing the proof of Theorem~\ref{Thm:5.0.2}}\label{SUBSECTION_Weyl_aff4}
To prove Theorem~\ref{Thm:5.0.2} it remains
\begin{enumerate} \item To check that any $W$-quasiessential triple is
$W$-essential and vice-versa. \item  To show how to determine a
$W$-essential part of a given admissible triple.
\end{enumerate}
The following lemma solves the both problems.
\begin{Lem}\label{Lem:5.5.1}
\begin{enumerate}
\item
Any  $W$-quasiessential triple is $W$-essential.
\item Let $(\g,\h,V)$ be an admissible triple and $\h_0$ a minimal ideal of $\h$ such that
the conditions $S^{(\alpha)}\rightsquigarrow_\g T^*(G*_HV)$ and
$S^{(\alpha)}\rightsquigarrow_\g T^*(G*_{H_0}V)$ are equivalent for
any root $\alpha$. Then $(\g,\h_0,V/V^{\h_0})$ is a
$W$-quasiessential triple and $W(\g,\h,V)=W(\g,\h_0,V/V^{\h_0})$.
\end{enumerate}
\end{Lem}

Assertion 2 of Lemma~\ref{Lem:5.5.1} implies that any $W$-\!\!
essential triple is $W$-quasiessential. Moreover, a triple
$(\g,\h_0,V/V^{\h_0})$ from assertion 2 is either the trivial triple
$(\g,0,0)$ or one of the triples from Table~\ref{Tbl:5.0.3}.
Inspecting Table~\ref{Tbl:5.0.3}, we note that an ideal
$\h_0\subset\h$ is determined uniquely. Therefore  the triple
$(\g,\h_0,V/V^{\h_0})$  is a $W$-essential part of $(\g,\h,V)$.

\begin{proof}
Suppose an admissible triple $(\g,\h,V)$ is  $W$-quasiessential but
not  $W$-essential. Then, by the definition of a $W$-essential
triple, there exists a proper ideal $\h_0\subset\h$ such that
$W(\g,\h,V)=W(\g,\h_0,V)$. On the other hand, there is $\alpha\in
\Delta(\g) $ such that $S^{(\alpha)}\rightsquigarrow_\g T^*(G*_HV)$,
$S^{(\alpha)}\not\rightsquigarrow_\g T^*(G*_{H_0}V)$.
Corollary~\ref{Cor:5.2.3} implies that $s_{w\alpha}\in W(\g,\h_0,V)$
for all $w\in W$. By the computation of the previous subsection,
there is $w\in W$ such that $s_{w\alpha}\not\in W(\g,\h,V)$.
Contradiction.

Proceed to assertion 2. The triple  $(\g,\h_0,V/V^{\h_0})$ is
$W$-quasiessential, thanks to the minimality condition for $\h_0$.
If $\h_0=\{0\}$, then, by Corollary~\ref{Cor:5.2.3},
$W(\g,\h,V)=W(\g)$. If $(\g,\h_0,V/V^{\h_0})$ is one of the triples
 NN28,32,33 from Table~\ref{Tbl:5.0.3}, then $\h=\h_0$ whence the equality of
 the Weyl groups. So we may assume that $\g\neq G_2$ and that $S^{(\alpha)}\rightsquigarrow_\g
T^*(G*_HV)$ iff $\alpha$ is a long root. By
Corollary~\ref{Cor:5.2.3}, $W(\g,\h,V)$ contains all reflections
corresponding to short roots. Assume, at first, that $\g$ is a
classical Lie algebra. By Proposition \ref{Prop:1.2.11},
$W(\g,\h,V)$ is one of the subgroups listed in
Table~\ref{Tbl:5.2.8}. But any of those groups containing all
reflections corresponding to short roots is maximal among all proper
subgroups of $W(\g)$ generated by reflections. Since
$W(\g,\h,V)\subset W(\g,\h_0,V/V^{\h_0})$, the previous subsection
yields $W(\g,\h_0,V/V^{\h_0})=W(\g,\h,V)$. It remains to consider
the case $\g=F_4$. Analogously to  the classical case, $W(\g,\h,V)$
contains a reflection corresponding to a long root. Otherwise,
$W(\g,\h,V)$ is not large in $W(\g)$ (in Definition~\ref{Def:1.2.14}
take
$\alpha=\varepsilon_1-\varepsilon_2,\beta=\varepsilon_2-\varepsilon_3$).
Any subgroup in $W(\g)$ containing all reflections corresponding to
short roots and some reflection corresponding to a long root is
maximal.
\end{proof}

\subsection{Weyl groups of affine homogeneous spaces}\label{SUBSECTION_Weyl_aff5}
In this subsection $\g$ is a reductive Lie algebra, $\h$ its
reductive subalgebra.

Before stating the main result we make the following remark. The
Weyl group $W(\g,\h)$ depends only on the commutant $[\h,\h]$, see
Corollary~\ref{Cor:3.4.6}. Therefore we may assume that  $\h$ is
semisimple. The strategy of the computation is similar to that
above: we compute the Weyl groups only for so called {\it
$W$-essential} subalgebras $\h\subset\g$ and then show how to reduce
the general case to this one.

\begin{defi}\label{Def:5.5.1}
A semisimple subalgebra  $\h\subset\g$ is called $W$-essential if
any ideal $\h_0\subset \h$ such that $\a(\g,\h)=\a(\g,\h_0),
W(\g,\h)=W(\g,\h_0)$ coincides with $\h$.
\end{defi}

\begin{Prop}\label{Prop:5.5.2} Let $\g$ be a semisimple Lie algebra and $\h$ its semisimple subalgebra.
Then the following claims hold.
\begin{enumerate}
\item Let $\h$ be a  $W$-essential subalgebra of $\g$. Then the pair $(\g,\h)$
is  contained either in Table~\ref{Tbl:5.0.3} or in
Table~\ref{Tbl:1.7.4}. In the latter case $W(\g,\h)$ coincides with
the group indicated in Table~\ref{Tbl:6.3.2}.
\item There is a unique minimal ideal $\h^{W-ess}\subset\h$ such that $\a(\g,\h^{W-ess})=\a(\g,\h),W(\g,\h^{W-ess})=W(\g,\h)$
It coincides with a maximal ideal of $\h$ that is a  $W$-essential
subalgebra of $\g$.
\end{enumerate}
\end{Prop}
\begin{proof}
The case  $\a(\g,\h)=\t$ is clear. By definition, any $\a$-essential
semisimple subalgebra of $\g$ is $W$-essential. Now let $\h$ be a
$W$-essential subalgebra in $\g$ with $\a(\g,\h)\varsubsetneqq \t$.
Let $L_0=L_{0\,G,G/H}^\circ$, $\underline{X}$ be the distinguished
component in $(G/H)^{L_0}$,
$\underline{G}=N_G(\l_0,\underline{X})/L_0,
\underline{H}=N_H(\l_0)/L_0$, $\underline{\t}=\underline{\g}\cap
\t$. Denote by $F$ the subgroup of $\underline{G}$ consisting of all
elements leaving invariant the distinguished Borel and Cartan
subalgebras in $\underline{\g}$. It follows from assertion 1 of
Proposition~\ref{Prop:6.0.2} that
$\underline{X}=\underline{G}/\underline{H}$. By
Theorem~\ref{Thm:3.0.3},
$W(\g,\h)=W(\underline{\g},\underline{\h})\leftthreetimes
F/\underline{T}$. From assertion 4 of Proposition~\ref{Prop:6.0.1}
it follows that $N_G(L_0)=N_G(L_0)^\circ N_{H^{ess}}(L_0)$.
Therefore $N_G(L_0)=N_G(L_0)^\circ N_H(L_0)$ or, equivalently,
$N_G(L_0,\underline{X})=N_G(L_0)$. Hence $N_G(L_0)^\circ
F=N_G(L_0)$. Inspecting Table~\ref{Tbl:6.2.2} and using
Theorem~\ref{Thm:5.0.2}, we see that for any subalgebra
$\underline{\h}_1\subset\underline{\g}$ such that $\h^{ess}$ is an
ideal $\underline{\h}_1$ and
$\a(\underline{\g},\underline{\h}_1)=\underline{\t}$ the equality
$W(\underline{\g},\underline{\h})=W(\underline{\g})$ holds. Thence
$W(\g,\h)=W(\g,\h^{ess})=W(\underline{\g})\leftthreetimes
F/\underline{T}$. In particular, if $\h^{ess}\neq \{0\}$, then
$\h^{W-ess}=\h^{ess}$.
\end{proof}

\section{Computation of the Weyl groups for homogeneous spaces}\label{SECTION_Weyl_homogen}
\subsection{Introduction}\label{SUBSECTION_Weyl_homogen_intro}
In this section we complete the computation of the groups
$W(\g,\h)$.

At first, note that Corollary \ref{Cor:3.4.6} allows one to reduce
the computation to the case when a maximal reductive subalgebra of
$\h$ is semisimple. In this case $G/H$ is quasiaffine, thanks to
Sukhanov's criterium, see \cite{Sukhanov}.

%

Now recall that we have reduced the computation of $W(\g,\h)$ to the
case when $\a(\g,\h)=\t$ (Theorem~\ref{Thm:3.0.3} and results of
Section \ref{SECTION_distinguished}) and $\g$ is simple
(Proposition~\ref{Prop:3.4.5}).

The computation of $W(\g,\h)$ for $\g=\so_5,G_2$ is carried out in
Subsection \ref{SUBSECTION_Weyl_homogen1}. Basically, it is built
upon Proposition \ref{Prop:1.2.9}. In the beginning of Subsection
\ref{SUBSECTION_Weyl_homogen3} we show how to compute the space
$\t^{W(\g,\h)}$. It turns out that for $\g$ of type $A,D,E$ the
group $W(\g,\h)$ is determined uniquely by $\t^{W(\g,\h)}$. This is
the main result of Subsection \ref{SUBSECTION_Weyl_homogen3}, its
proof is based on the restrictions on $W(\g,\h)$ obtained in
Proposition \ref{Prop:8.2.1}. The computation for algebras  of type
$C,B,F$ and rank bigger than 2 is carried out in Subsections
\ref{SUBSECTION_Weyl_homogen4}-\ref{SUBSECTION_Weyl_homogen6}.
Remark \ref{Rem:8.3.7} allows us to make an   additional restriction
on $\h$: we require $\t^{W(\g,\h)}=0$.

\begin{Rem}\label{Rem:8.0.1}
The computation of $W(\g,\h)$   for all groups of rank 2 essentially
does not use Proposition \ref{Prop:8.2.1} and results of
\cite{Wasserman} on the classification of all wonderful varieties of
rank 2. Let us note that the classification of all spherical
varieties of rank 2 is not very difficult. The reductions described
above allow to reduce the computation of the Weyl groups of
varieties of rank 2 to the case when $G$ itself has rank 2.
\end{Rem}

\subsection{Types $B_2,G_2$}\label{SUBSECTION_Weyl_homogen1}
In this subsection $G=\SO_5,G_2$ and  $\h$ denotes a subalgebra of
$\g$ such that $\a(\g,\h)=\t$ and a  maximal reductive subalgebra
$\h_0$ of $\h$ is semisimple.

Suppose that $\Rad_u(\q)\subset \h$ for some parabolic subalgebra
$\q\subset\g$. We may assume that $\q$ is antistandard, let $\m$ be
its standard Levi subalgebra. Then, by Corollary~\ref{Cor:3.4.7},
$W(\g,\h)=W(\m,\m\cap\h)$.  Below in this subsection we suppose that
$\h$ does not contain the unipotent radical of a parabolic.

\begin{Prop}\label{Prop:8.1.3}
Let $\g=\so_5$ and a subalgebra $\h\subset\g$ satisfy the above
assumptions. Then the following conditions are equivalent.
\begin{enumerate}
\item $W(\g,\h)\neq W(\g)$.
\item $\h$ is conjugate to $\sl_2^{diag}$.
\end{enumerate}
Under conditions (1),(2),  $W(\g,\h)$ is generated by
$s_{\varepsilon_i},i=1,2$.
\end{Prop}
\begin{proof}
We have already computed the groups $W(\g,\h)$ for reductive
subalgebras $\h$. Suppose now that $W(\g,\h)\neq W(\g)$. Put
$\q_1=\b^{-}+\g^{\varepsilon_2},\q_2=\b^{-}+\g^{\varepsilon_1-\varepsilon_2}$.

Firstly, we consider the case when $H$ is unipotent. Then
$\dim\h\leqslant 4$. Since $\h$ does not contain  the unipotent
radical of a parabolic, we see that $\dim\h\leqslant 3$ and if
$\dim\h=3$, then $\h\sim_G
\widetilde{\h}:=\Span_\C(e^{-\varepsilon_1-\varepsilon_2},
e^{-\varepsilon_1},e^{\varepsilon_2-\varepsilon_1}+e^{-\varepsilon_2})$.
The closure of $G\widetilde{\h}$ in $\Gr(\g,3)$ contains
$G\Rad_u(\q_i),i=1,2$. Applying Corollary~\ref{Cor:3.4.7} and
Proposition~\ref{Prop:1.2.9}, we get $s_{\alpha_1},s_{\alpha_2}\in
W(\g,\h)$.

If $\dim\h=1$, then there exists a unipotent subalgebra
$\widehat{\h}$ of dimension 2 containing $\h$. Indeed
$\n_{\u_1}(\h)\neq\h$, where $\u_1$ denotes the maximal unipotent
subalgebra of $\g$ containing $\h$. Let us check that
$W(\g,\h)=W(\g)$ whenever $\h$ is a unipotent subalgebra of
dimension 2. By Proposition~\ref{Prop:1.2.9}, we may assume that
$G\h\subset \Gr(\g,2)$ is closed or, equivalently, $\n_\g(\h)$ is a
parabolic subalgebra of $\g$. From this one easily deduces that
$\h\sim_G
\h_0:=\Span_\C(e^{-\varepsilon_1-\varepsilon_2},e^{-\varepsilon_1})$.
But  $\h_0$ is conjugate to a subalgebra in $\widetilde{\h}$ whence,
thanks to Corollary~\ref{Cor:1.2.2}, $W(\g,\h)=W(\g)$.

It remains to consider the case $\h_0\cong\sl_2$. It is easily
deduced from $\a(\g,\h)=\t$ that $\h$ is semisimple.
\end{proof}

\begin{Prop}\label{Prop:8.1.4}
Suppose $\g=G_2$ and $\h$ satisfies the above assumptions. Then the
following conditions are equivalent:
\begin{enumerate}
\item $W(\g,\h)\neq W(\g)$.
\item $\h$ is conjugate to one of the following subalgebras
$\h_1:=\g^{(\alpha_2)}\oplus \g^{-3\alpha_1-\alpha_2}\oplus
\g^{-3\alpha_1-2\alpha_2}$, $\h_2:=\h_1+ \g^{-2\alpha_1-\alpha_2}$.
\end{enumerate}
  $W(\g,\h_1)$ is generated by
  $s_{\alpha_1},s_{\alpha_1+\alpha_2}$, and
 $W(\g,\h_2)$ is generated by $s_{\alpha_1+\alpha_2}$.
\end{Prop}
\begin{proof}
Again, put $\q_1=\b^{-}+\g^{\alpha_1},\q_2=\b^{-}+\g^{\alpha_2}$
($\alpha_2$ is a long root).

The subalgebra $\h_1$ coincides with the s.s.g.p. for the action
$G:G*_{A_2}\C^3$ and the equality for $W(\g,\h_1)$ is deduced from
Theorem \ref{Thm:5.0.2}.  Now put
$\widetilde{\h}_2:=\h_2+\C\alpha_1^\vee$. Clearly,
$\widetilde{\h}_2\subset \n_\g(\h_2)$. The subalgebra
$\widetilde{\h}_2$ is tamely contained in $\q_2$. Applying results
of Subsection~\ref{SUBSECTION_prelim1}, we obtain
$\a(\g,\widetilde{\h}_2)=\C(\alpha_1+\alpha_2)$. Since
$\widetilde{\h}_2$ does not contain a maximal unipotent subalgebra,
we have $W(\g,\widetilde{\h}_2)\neq \{1\}$. The required equality
for $W(\g,\h_2)$ stems from Corollary~\ref{Cor:3.4.6}.

Now suppose (1) holds.   Firstly, we show that
$\h_0\not\sim_G\g^{(\alpha_2)}$ implies $\h_0=\{0\}$. Assume the
converse. Since $\a(\g,\h)=\t$, we have $\h_0\neq A_2$. So if
$\h_0\neq 0$, then $\h_0$ contains an ideal conjugate to
$\g^{(\alpha_1)}$. We easily check that $\h$ is conjugate to a
subalgebra in $\g^{(\alpha_1)}\oplus\g^{(3\alpha_1+2\alpha_2)}$. The
Weyl group of the latter coincides with $W(\g)$. So $\h_0=\{0\}$, in
other words, $\h$ is unipotent.

We have $\dim\h\leqslant 5$. If $\dim\h=5$, then, thanks to
$\h\not\sim_G\Rad_u(\q_1),\Rad_u(\q_2)$, we get
$\h\sim_G\widetilde{\h}:=[\b^{-},\b^{-}]+
\C(e^{-\alpha_1}+e^{-\alpha_2})$. As in the proof of Proposition
\ref{Prop:8.1.4}, we get $W(\g,\h)=W(\g)$.

Now let $\dim\h\leqslant 4$. Let us check that $W(\g,\h)=W(\g)$. By
Proposition \ref{Prop:1.2.9}, it is enough to check the claim only
when $\n_\g(\h)$ is a parabolic subalgebra of $\g$. But in this case
$\h\sim_G [\b^-,\b^-]\subset\widetilde{\h}$. From
Proposition~\ref{Prop:1.2.1} it follows that $W(\g,\h)=W(\g)$.

So it remains to consider the case $\h_0=\g^{(\alpha_2)}$. By
Theorem~\ref{Thm:5.0.2},  $\h$ is not reductive. We may assume that
$\h$ is tamely contained in $\q_2$. Since $\h\not\subset
\g^{(\alpha_2)}\oplus \g^{(\alpha_2+2\alpha_1)}$ and
$\Rad_u(\q_2)\not\subset\h$, we get the required list of subalgebras
$\h$.
\end{proof}

\subsection{Types $A,D,E$}\label{SUBSECTION_Weyl_homogen3}
In this subsection $\g$ is a simple Lie algebra and $\h$ its
subalgebra such that $\a(\g,\h)=\t$ and a maximal reductive
subalgebra of $\h$ is simple.

 Firstly, we show how to compute the space $\t^{W(\g,\h)}$.

\begin{Prop}\label{Prop:8.3.6}
Let $\g,\h$ be such as above and $\z$ a commutative reductive
subalgebra of $\n_\g(\h)/\h$ containing  all semisimple elements of
$\z(\n_\g(\h)/\h)$. Let $\widetilde{\h}$ denote the inverse image of
$\z$ in $\n_\g(\h)$ under the natural epimorphism
$\n_\g(\h)\twoheadrightarrow \n_\g(\h)/\h$. Then
$\a(\g,\widetilde{\h})$ is the orthogonal complement to
$\t^{W(\g,\h)}$ in $\t$.
\end{Prop}
\begin{proof}
By Corollary \ref{Cor:3.4.6}, $\t^{W(\g,\h)\perp}\subset
\a(\g,\widetilde{\h})$. Let $\widehat{H}$ denote the inverse image
of the subgroup $\A_{G,G/H}^\circ\subset N_G(H)/H$ (see Definition
\ref{Def:1.2.13}) in $N_G(H)/H$. By Proposition \ref{Prop:1.2.10},
$\a(\g,\widehat{\h})=\t^{W(\g,\h)\perp}$.  By
Lemma~\ref{Lem:1.2.11}, $\widehat{H}/H\subset Z(N_G(H)/H)$. By
Lemma~\ref{Lem:1.2.12}, $\widehat{H}/H$ is a torus. These two
observations yield $\widehat{H}\subset \widetilde{H}$ whence
$\t^{W(\g,\h)\perp}\supset \a(\g,\widetilde{\h})$.
\end{proof}

\begin{Rem}\label{Rem:8.3.7}
We use the notation of Proposition \ref{Prop:8.3.6} and for $\z$
take the ideal of $\z(\n_\g(\h)/\h)$ consisting of all semisimple
elements. Suppose that $\a(\g,\widetilde{\h})\neq \t$. Then we can
reduce the computation of $W(\g,\h)=W(\g,\widetilde{\h})$ to the
computation of $W(\underline{\g},\underline{\h})$, where $\rank
[\underline{\g},\underline{\g}]<\rank\g$, $\rank
\underline{G}/\underline{H}=\rank\underline{G}$, as follows.

Put $\widetilde{G}=G\times \C^\times$. There is a quasiaffine
homogeneous space $\widetilde{X}$ of $\widetilde{G}$ such that there
is a $G$-equivariant principal $\C^\times$-bundle
$\widetilde{X}\rightarrow G/\widetilde{H}$, see Subsection
\ref{SUBSECTION_algor2} for details. The stable subalgebra
$\widetilde{\h}_1$ of $\widetilde{X}$ is naturally identified with
$\widetilde{\h}$. By Proposition \ref{Prop:1.2.5},
$W(\g,\widetilde{\h})\cong W(\widetilde{\g},\widetilde{\h}_1)$. Then
we apply Theorem \ref{Thm:3.0.3} to $\widetilde{X}$.

 However, for some algebras $\g$ the group
$W(\g,\h)$ is uniquely determined by $\t^{W(\g,\h)}$, see below, so
the reduction described above is unnecessary.
\end{Rem}

\begin{Prop}\label{Prop:8.3.1}
Let $\g=\sl_n,n\geqslant 3,\so_{2n},n\geqslant 4$ or $E_l,l=6,7,8$,
and $\h$ be as above. Then $W(\g,\h)=W(\g)\cap
Z_G(\t^{W(\g,\h)})/T$.
\end{Prop}
\begin{proof}
Set $\Delta_0:=\Delta(\g)\cap\Span_\C(\widehat{\Pi}(\g,\h))$. Our
goal is to prove that $\widehat{\Pi}(\g,\h)$ is a system of simple
roots in $\Delta_0$. Let $\Pi_0$ be the system of simple roots in
$\Delta_0$ positive on the dominant Weyl chamber of $\g$. Then for
any $\alpha\in \widehat{\Pi}(\g,\h)$ there are positive integers
$n_\gamma,\gamma\in \Pi_0,$ such that $\alpha=\sum_{\gamma\in
\Pi_0}n_\gamma\gamma$ and for any $\gamma\in \Pi_0$ there are
rationals $m_\alpha,\alpha\in \widehat{\Pi}(\g,\h)$ such that
$\gamma=\sum_{\alpha\in \widehat{\Pi}(\g,\h)}m_\alpha\alpha$. So for
any $\gamma\in \Pi_0$ there is $\alpha\in \widehat{\Pi}(\g,\h)$ such
that $\Supp(\gamma)\subset \Supp(\alpha)$.

Set
$\Pi_1:=\widehat{\Pi}(\g,\h)\cap\Pi(\g),\widehat{\Pi}_2:=\widehat{\Pi}(\g,\h)\setminus\Pi_1,
\Pi_2:=\cup_{\beta\in \widehat{\Pi}_2}\Supp(\beta)$. Recall
(Proposition \ref{Prop:8.2.1}) that for any
$\beta\in\widehat{\Pi}_2$ there are adjacent simple roots
$\beta_1,\beta_2$ such that $\beta=\beta_1+\beta_2$. It follows that
$\Pi_0\subset \widehat{\Pi}(\g,\h)\cup\Pi_2$. It is enough to prove
that $\Pi_0\cap\Pi_2=\varnothing.$

Since $(\alpha,\beta)\leqslant 0$ for any $\alpha,\beta\in
\widehat{\Pi}(\g,\h)$, we see that $\alpha\not\in \Supp(\beta)$ for
any $\alpha\in \Pi_1, \beta\in \widehat{\Pi}_2$. Our claim will
follow if we check $\Span_\C(\widehat{\Pi}_2)\cap
\Pi(\g)=\varnothing$. Assume the converse, there is a subset
$\Sigma\subset \widehat{\Pi}_2$ such that $\alpha:=\sum_{\beta\in
\Sigma}n_\beta \beta\in \Pi(\g)$ with all $n_\beta$ nonzero. There
are at least two simple roots in $\Pi_2$ such that no more than one
of their neighbors  lies in $\cup_{\beta\in\Sigma}\Supp(\beta)$. Let
$\beta_1,\beta'_1$ be such roots. Then there are $\beta,\beta'\in
\Sigma$ such that $\beta_1\not\in \Supp(\gamma)$ (resp.
$\beta_1'\not\in \Supp(\gamma)$) for any $\gamma\in
\Pi_2,\gamma\neq\beta$ (resp., $\gamma\neq\beta'$). It follows that
$\beta_1,\beta_1'\in \Supp(\sum_{\beta\in \Sigma}n_\beta \beta)$,
contradiction.
\end{proof}

\subsection{Type $C_l,l>2$}\label{SUBSECTION_Weyl_homogen4}
In this subsection $\g$ is a simple Lie algebra and $\h$ its
subalgebra such that $\a(\g,\h)=\t$ and a maximal reductive
subalgebra of $\h$ is semisimple.

\begin{Prop}\label{Prop:8.4.1}
Suppose $\g\cong\so_{2n+1},\sp_{2n},n>2,F_4$. Choose roots
$\alpha,\beta\in\Delta(\g)$ as follows:
$(\alpha,\beta)=(\varepsilon_{n-1},\varepsilon_{n})$ for $\g=
\so_{2n+1}$,
$(\alpha,\beta)=(\varepsilon_{n-1}+\varepsilon_n,\varepsilon_{n-1}-\varepsilon_{n})$
for $\g=\sp_{2n}$, $(\alpha,\beta)=(\varepsilon_3,\varepsilon_4)$
for $\g=F_4$. Suppose $s_{\alpha},s_{\beta}\in W(\g,\h),
s_{\alpha+\beta}\not\in W(\g,\h)$. Then there is a subalgebra
$\s\subset\g$ satisfying the following two conditions:
\begin{enumerate}
\item $\s\sim_G \g^{(\alpha+\beta)}$.
\item $\g/(\g^\s\oplus\s)\cong(\h/\h^\s)^{\oplus 2}\oplus (\C^2)^{\oplus 2}$ (an isomorphism of $\h$-modules).
\end{enumerate}
\end{Prop}
\begin{proof}
For two $\s$-modules $V_1,V_2$ we write $V_1\sim V_2$ if
$V_1/V_1^\s\cong V_2/V_2^\s$. Set
$M:=TG^{(\beta,\alpha-\beta)},Q:=B^-M,X:=G/H$. Let $Z_0$ denote a
rational quotient for the action $\Rad_u(Q):X$. Modifying $Z_0$ if
necessary, one may assume that $M$ acts regularly on $Z_0$. By
Proposition 8.2 from \cite{comb_ham}, $W_{M,Z_0}$ is generated by
$s_{\alpha},s_{\beta}$. Thanks to assertion 4 of Proposition
\ref{Prop:1.2.1}, $W_{M,Z_0}=W(\m,\m_z)$ for a general point $z\in
Z_0$. From Proposition~\ref{Prop:8.1.3} it follows that $\m_z$ is
reductive and $\s:=[\m_z,\m_z]\sim_M \m^{(\alpha+\beta)}$. By
\cite{comb_ham}, Lemma 7.12,  the action $\Rad_u(Q):G/H$ is locally
free. Since $Z_0$ is a rational quotient for the action
$\Rad_u(Q):G/H$, we see that $\q_x\sim_Q \g^{(\alpha+\beta)}$ for
$x\in X$ in general position. Thence there is $x\in X$ such that
$\s\subset\g_x$ and $T_xX\sim \q/\s$. We may assume that $x=eH$. It
follows that $\g/\h\sim \q/\s\sim\Rad_u(\q)+(\C^2)^{\oplus 2}$. Note
that $\g\sim \Rad_u(\q)^{\oplus 2}\oplus (\C^2)^{\oplus 2}\oplus\s$
whence (2).
\end{proof}

\begin{Lem}\label{Lem:8.4.2}
Let $\g,\h$ satisfy the assumptions of the previous proposition and
 $\s\subset \h$ possesses properties 1,2. Further, let $Q$ be a
 parabolic subgroup of $G$ such that $\h$ is tamely contained in
 $\q$ and $M$ a Levi subgroup of $Q$ such that $\h_0:=\m\cap\h$ is a maximal reductive
  subalgebra of $\h$ containing $\s$. Then there are simple ideals $\h_1\subset\h_0,\m_1\subset\m_0$ such that
  $\s\subset\h_1\subset\m_1$ and
$S^{\s}\rightsquigarrow_{\h_1}\m_1/\h_1\oplus
\Rad_u(\q)/\Rad_u(\h)\oplus (\Rad_u(\q)/\Rad_u(\h))^*$.
\end{Lem}
\begin{proof}
We preserve the notation of Proposition \ref{Prop:8.4.1}. Note that
$\alpha+\beta$ is a long root. Therefore $\s$ is contained in a
simple ideal $\m_1\subset\m$ and $i(\s,\m_1)=1$, by
(\ref{eq:5.3:1}). Moreover, $\s\sim_M \m^{(\gamma)}$ for some root
$\gamma\in \Delta(\m)$ (see Lemma~\ref{Lem:5.3.3}). Analogously,
$\s$ is contained in some simple ideal $\h_{1}\subset \h_0$. Since
$\s\subset \m_1$, we get $\h_1\subset\m_1$. Property 2 of
Proposition \ref{Prop:8.4.1} implies that the nontrivial parts of
the $\s$-modules $(\m/\h_0)\oplus V\oplus V^*$ and $(\h_1/\s)\oplus
(\C^2)^{\oplus 2}$ coincide, where $V:=\Rad_u(\q)/\Rad_u(\h)$.
Thence, see the proof of Lemma~\ref{Lem:5.3.6},
$l_{\h_1}(\m/\h_0\oplus V\oplus V^*)=1-\frac{4}{k_{\h_1}}$. Besides,
we remark that the $\h_1$-module $(\m/\h_0)\oplus V\oplus V^*$ is
orthogonal and its s.s.g.p is trivial. The latter follows easily
from the observation that the s.s.g.p. for the action
$M:M*_H(\m/\h\oplus V\oplus V^*)$ is trivial. One can list all such
modules by using table from~\cite{AEV} and results from
\cite{Elash1}. It turns out that the nontrivial part of such a
module is presented in Table~\ref{Tbl:5.3.10}.
\end{proof}

Our algorithm for computing $W(\g,\h)$ for $\g\cong\sp_{2n}$ is
based on the following proposition.

\begin{Prop}\label{Prop:8.4.3}
Suppose  $\g\cong\sp_{2n},n>2$. Assume, in addition, that
$\t^{W(\g,\h)}=0$. Further, suppose $\h$ is tamely contained in an
antistandard parabolic subalgebra
 $\q\subset \g$ and $\m\cap\h$ is a maximal reductive subalgebra of $\h$, where $\m$ is the standard Levi
 subalgebra of $\q$. Then the following conditions are equivalent.
\begin{enumerate}
\item $W(\g,\h)\neq W(\g)$.
\item
 $\g^{(\alpha_{n-1},\alpha_n)}\subset\q$, there are simple ideals $\h_1\subset\h_0, \m_1\subset\m$
 such that
$\m_1=\g^{(\alpha_k,\ldots,\alpha_n)},k\leqslant n$, $
\h_1=\sp_{2l}\subset\m_1\cong\sp_{2(n-k)}$, and the $\h_1$-modules
$\Rad_u(\q)/(\Rad_u(\h)+\Rad_u(\q)^{\h_1})$ and $(\C^{2l})^{\oplus
4l-2(n-k)}$ are isomorphic.
\end{enumerate}
Under  conditions (1),(2), $W(\g,\h)$ is generated by all
reflections corresponding to short roots.
\end{Prop}
\begin{proof}
Let us introduce some notation. Let $I$ be a finite set. By $V$ we
denote the vector space with the basis $\varepsilon_i,i\in I$. Let
$A_I$ (resp., $B_I,C_I,D_I$) be a linear group acting on $V$ as the
Weyl
group of type $A_{\#I}$ (resp. $B_{\#I}, C_{\#I},D_{\#I}$).%

Since $W(\g,\h)$ is generated by reflections and $\t^{W(\g,\h)}=0$,
we see that there exists a partition $I_1,\ldots, I_k$ of
$\{1,2,\ldots, n\}$ such that $W(\g,\h)=\prod_{i=1}^k \Gamma_i$,
where $\Gamma_j=D_{I_j}$ or $C_{I_j}$.  From
Proposition~\ref{Prop:8.2.1} it follows that
$\varepsilon_i+\varepsilon_j\not\in \widehat{\Pi}(\g,\h)$ unless
$(i,j)=(n,n),(n,n-1)$ or $(n-1,n)$. Thus $k=1$. (2)$\Rightarrow$(1)
stems from Corollary~\ref{Cor:3.4.3} applied to $\q$ and results of
the previous section. Now assume that $W(\g,\h)\neq W(\g)$ whence
$W(\g,\h)=D_n$. The assumptions of Proposition~\ref{Prop:8.4.1} are
satisfied. By that proposition, $\h_0$ contains a subalgebra
conjugate to $\g^{(\alpha_n)}$. Since $\alpha_n$ is a long root, any
subalgebra of $\m$ that is $G$-conjugate to $\g^{(\alpha_n)}$ is
contained in an ideal $\m_1$ of the form
$\g^{(\alpha_k,\ldots,\alpha_n)}$. Now (1)$\Rightarrow$(2) follows
from Lemma~\ref{Lem:8.4.2} and assertion 1 of
Proposition~\ref{Prop:5.3.1}.
\end{proof}

\subsection{Type $B_l,l>2$}\label{SUBSECTION_Weyl_homogen5}
 In this subsection, if otherwise is not indicated, $\g=\so_{2n+1},n\geqslant 3$, and $\h\subset\g$
 is such that $\a(\g,\h)=\t,\t^{W(\g,\h)}=\{0\}$ and the maximal reductive subalgebra of $\h$ is semisimple.
\begin{Lem}\label{Lem:8.5.1}
If $W(\g,\h)\neq W(\g)$, then $W(\g,\h)=B_I\times B_J$, where $I$ is
one of the following subsets of  $\{1,2\ldots,n\}$ and
$J=\{1,\ldots,n\}\setminus I$:
\begin{enumerate}
\item $I=\{n,n-2,\ldots,n-2k\}, 0\leqslant k< \frac{n-1}{2}$.
\item $I=\{n-1,\ldots,n-2k+1\},0<k<\frac{n}{2}$.
\end{enumerate}
\end{Lem}
\begin{proof}
As in the proof of Proposition~\ref{Prop:8.4.3}, there is a
partition $I_1,\ldots,I_l$ of $\{1,2\ldots,n\}$ such that
$W(\g,\h)=\prod_{j=1}^l\Gamma_j$, where $\Gamma_j=B_{I_j}$ or
$D_{I_j}$. Since
$\varepsilon_{i}+\varepsilon_{j},\varepsilon_k\not\in\widehat{\Pi}(\g,\h)$
for arbitrary $i,j$ and $k<n-1$ (Proposition~\ref{Prop:8.2.1}), we
get $l=2,\Gamma_j=B_{I_j}$.  It follows that $n\in I_1, n-1\in I_2$.
By Proposition \ref{Prop:8.2.1}, if $k,k+1\in I_j,j=1,2$ for
$1<k<n$, then  $k-1$ or $k+2$ lies in $I_j$. This implies the claim
of the present proposition.
\end{proof}

We remark that, by definition, $1\not\in I$.

\begin{Lem}\label{Lem:8.5.3}
Let $\g \cong \so_{2n+1},n\geqslant 3$ or $F_4$. Suppose $W(\g,\h)$
is a maximal proper subgroup generated by reflections in $W(\g)$
(for $\g=\so_{2n+1}$ we have checked this in Lemma \ref{Lem:8.5.1}).
Let $\q,\m,\h_1,\m_1$ be such as in Lemma \ref{Lem:8.4.2}. Then:
\begin{enumerate}
\item Subalgebras $\h_1\subset\h$, $\m_1\subset\m$ are determined uniquely.
\item Let $\underline{\h}$ denote the subalgebra in $\h$ generated by $\h_1$ and
 $[\h_1,\Rad_u(\h)]$. Then
$W(\g,\underline{\h})=W(\g,\h)$.
\end{enumerate}
\end{Lem}
\begin{proof}
It follows from assertion 1 of Proposition \ref{Prop:5.3.1} that
$\m_1\neq \m_2$ for two different pairs $(\m_1,\h_1), (\m_2,\h_2)$
satisfying the assumptions of Lemma~\ref{Lem:8.4.2}. Thanks to
Corollary \ref{Cor:1.2.2}, $W(\g,\h)\subset W(\g,\underline{\h})$.
The subalgebra $\underline{\h}$ is tamely contatined in $\q$.
Applying Corollary~\ref{Cor:3.4.3} to $\q$, we see that
$W(\m_2)\not\subset W(\g,\h)$, $W(\m_2)\subset
W(\g,\underline{\h})$. Since $W(\g,\h)$ is maximal, this proves
assertion 1. Analogously, $W(\m_1)\not\subset W(\g,\underline{\h})$,
since the nontrivial parts of the  $\h_1$-modules
$\Rad_u(\q)/\Rad_u(\h), \Rad_u(\q)/\Rad_u(\underline{\h})$ are the
same. This proves assertions 2.
\end{proof}

For indices $i_1,i_2$ we write $i_1\sim i_2$ whenever $i_1,i_2\in I$
or $i_1,i_2\in J$.

\begin{Lem}\label{Lem:8.5.2}
Suppose $\h$ is tamely contained in an antistandard parabolic
subalgebra $\q\subset\g$ and some maximal reductive subalgebra
$\h_0\subset\h$ is contained in the standard Levi subalgebra
$\m\subset\g$. Then
\begin{enumerate}
\item $\m_1=\g^{(\alpha_{k+1},\ldots,\alpha_l)}$, where $n-2\leqslant l\leqslant
n$, $i\sim j$ for $i,j\leqslant k$. Further, for $i,j\in
\{k+1,k+2,\ldots, \min(l+1,n)\}$ we have
$s_{\varepsilon_i-\varepsilon_j}\in
W(\m_1,\h_1,\Rad_u(\q)/\Rad_u(\h))$ iff $i\sim j$.
\item Suppose $k+1\sim k+2$.  Let $m\in \{k+1,k+2,\ldots, \min(l,n-1)\}$ be the maximal integer such that
$m\sim m+1$. Then $I=\{m+2k, k=\overline{1,[(n-m)/2]}\}$.
\end{enumerate}
\end{Lem}
\begin{proof}
Thanks to Lemma \ref{Lem:8.5.1}, $W(\g,\h)$  is maximal among all
proper subgroups of $W(\g)$ generated by reflections. There are
$k,l$ such that $\m_1=\g^{(\alpha_{k+1},\ldots,\alpha_l)}$.
 Let
$\underline{\q}\subset\g$ denote the antistandard parabolic
subalgebra of $\g$ corresponding to simple roots
$\alpha_1,\ldots,\alpha_{k-1},$
$\alpha_{k+1},\ldots,\alpha_l,\alpha_{l+2},$ $\ldots,\alpha_n$.
Since $[\h_1,\Rad_u(\q)]\subset\Rad_u(\underline{\q})$, we see that
$\underline{\h}$ is tamely contained in $\underline{\q}$. Now the
both assertion of this lemma stem from Corollary~\ref{Cor:3.4.3}
applied to $\underline{\q}$ and the explicit form of $W(\g,\h)$
indicated in Lemma~\ref{Lem:8.5.1}.
\end{proof}

Summarizing results obtained above in this subsection, we see that
it remains to compute $W(\g,\h)$ under the following assumptions on
$(\g,\h)$:
\begin{itemize}
\item[(a)] $\a(\g,\h)=\t,\t^{W(\g,\h)}=0$.
\item[(b)] A maximal reductive subalgebra $\h_1\subset\h$ is simple
and the Lie algebra $\Rad_u(\h)$ is generated by
$[\h_1,\Rad_u(\h)]$. \item[(c)] There are integers $
l\in\{n-2,n-1,n\},k>l$ such that $\h$ is tamely contained in the
antistandard parabolic subalgebra $\q\subset\g$ corresponding to
$\Pi(\g)\setminus \{\alpha_k,\alpha_{l+1}\},\h_1\subset
\m_1:=\g^{(\alpha_{k+1},\ldots,\alpha_l)}$. \item[(d)] $\h_1$ is not
contained in a proper Levi subalgebra of $\m_1$.
\item[(e)] The group $W(\m_1,\h_1,\Rad_u(\q)/\Rad_u(\h))$ does not
contain a reflection of the form
$s_{\varepsilon_i-\varepsilon_{i+1}},$ $k<i<l$.
\end{itemize}

\begin{Prop}\label{Prop:8.5.4}
Let  $\h,\h_1,\q,\m_1$ be such as in (a)-(e). The inequality
$W(\g,\h)\neq W(\g)$ holds iff $\h$ is $G$-conjugate to a subalgebra
from the following list:

1) $\h$ consists of all block matrices of the following form

\begin{equation*}
\begin{pmatrix}
0&0&0&0&0\\
0&0&0&0&0\\
x_{31}&0&x_{33}&0&0\\
0&0&0&0&0\\
x_{51}&0&x_{31}'&0&0
\end{pmatrix}
\end{equation*}
The sizes of blocks (from left to right and from top to bottom):
$k-1,1,7,1,k-1$; $x_{ij}'$  denotes the matrix $-I_px_{ij}^TI_q$,
where
$I_p=(\delta_{i+j,p+1})_{i,j=1}^p,I_q=(\delta_{i+j,q+1})_{i,j=1}^{q}$,
$x_{33},x_{51}$ are arbitrary elements of $G_2\hookrightarrow
\so_7,\so_k$. In this case  $I=\{n-1\}$.

2) $\h$ consists of all block matrices of the following form
\begin{equation*}
\begin{pmatrix}
0&0&0&0&0\\
x_{21}&x_{22}&0&0&0\\
0&0&0&0&0\\
x_{41}&0&0&x_{22}'&0\\
x_{51}&x_{41}'&0&x_{21}'&0
\end{pmatrix}
\end{equation*}

Here the sizes of blocks are $k,n-k,1,n-k,k$; $x_{22},x_{51}$ are
arbitrary elements of $\sl_{n-k},\so_k$. We have
$I=\{k+2,k+4,\ldots,k+2i,\ldots\}$.

3) $\h$ consists of all block matrices of the following form:
\begin{equation}\label{eq:8.5:2}
\begin{pmatrix}0&0&0&0&0&0&0\\
0&0&0&0&0&0&0\\
x_{31}&x_{32}&x_{33}&0&0&0&0\\
0&0&0&0&0&0&0\\
x_{51}&x_{52}&\iota(x_{32})&0&x_{33}'&0&0\\
x_{61}&0&x_{52}'&0&x_{32}'&0&0&\\
x_{71}&x_{61}'&x_{51}'&0&x_{31}'&0&0\\
\end{pmatrix}
\end{equation}
The sizes of blocks (from left to right and from top to bottom):
$k-1,1,3,1,3,1,k-1$;  $\iota$ denotes an isomorphism of
$\sl_3$-modules $\C^3$ and $\bigwedge^2\C^{3*}$; $x_{33},x_{71}$ are
arbitrary elements of $\sl_3$ and $\so_{k-1}$. Here $I=\{n-1\}$.

4) $\h$ consists of all block matrices of the following form
\begin{equation*}
\begin{pmatrix}
0&0&0&0&0\\
x_{21}&x_{22}&0&0&0\\
0&0&0&0&0\\
x_{41}&x_{32}&0&x_{22}'&0\\
x_{51}&x_{41}'&0&x_{21}'&0
\end{pmatrix}
\end{equation*}

Here the sizes of blocks are $k,n-k,1,n-k,k$ and $n-k$ is even;
$x_{22},x_{51}$ are arbitrary elements of $\sp_{n-k},\so_k$, while
$x_{32}$ lies in the nontrivial part of the $\sp_{n-k}$-module
$\bigwedge^2\C^{n-k}$. We have $I=\{k+2,k+4,\ldots,k+2i,\ldots\}$.
Here $I=\{k+2,k+4,\ldots,n\}$.

5) $\h$ consists of all block matrix of the following form
\begin{equation}\label{eq:8.5:11}
\begin{pmatrix}
0&0&0&0&0&0&0\\
0&0&0&0&0&0&0\\
x_{31}&x_{32}&x_{33}&0&0&0&0\\
x_{41}&x_{42}&\iota(x_{32})&0&0&0&0\\
x_{51}&x_{52}&x_{53}+x_{42}\omega&\iota(x_{32})'&x_{33}'&0&0\\
x_{61}&0&x_{52}'&x_{42}'&x_{32}'&0&0\\
x_{71}&x_{61}'&x_{51}'&x_{41}'&x_{31}'&0&0
\end{pmatrix}
\end{equation}

The sizes of blocks are $k-1,1,r,1,r,1,k-1$; $x_{ij}'$ have the same
meaning as in (\ref{eq:8.5:2}); $\iota$ denotes an isomorphism of
$\sp_r$-modules $\C^r,\C^{*r}$; $x_{71}\in
\so_{k-1},x_{33}\in\sp_{r}$, $x_{53}$ lies in the nontrivial part of
the $\sp_r$-module $\bigwedge^2\C^{r*}$, and $\omega$ is an
appropriate nonzero element in $(\bigwedge^2\C^{r*})^{\h_1}$. We
have $I=\{k+1,k+3,\ldots,n-1\}$.
\end{Prop}
Note that in all five cases the equality $\t^{W(\g,\h)}=0$ can be
checked by using Proposition \ref{Prop:8.3.6} (compare with the
proof of  Proposition~\ref{Prop:8.6.2} below).

\begin{proof}
Note that $\Rad_u(\q)/\Rad_u(\h)$ is an $\h_1$-submodule in
$\Rad_u(\q)$. Using Theorem \ref{Thm:5.0.2}, we obtain the following
list of possible triples
$(\m_1,\h_1,\Rad_u(\q)/(\Rad_u(\h)+\Rad_u(\q)^{\h_1}))$:

\begin{itemize}
\item $(\sl_r,\sl_r,\bigwedge^2\C^{r*}\oplus \C^{r*}),r\geqslant 3$.
\item $(\sl_r,\sl_r,\bigwedge^2\C^{r*}\oplus \C^r)$,  $r\geqslant 4$ is even.
\item $(\sl_r,\sp_r,\C^{r})$,  $r\geqslant 2$ is even.
\item $(\so_7,G_2,V(\pi_1))$.
\end{itemize}

{\it The case  $(\sl_r,\sl_r,\bigwedge^2\C^{r*}\oplus \C^r)$.} Let
us check that this triple cannot occur. Since $r>3$, we see that
$\Rad_u(\q)$ contains a unique  $\m_1$-submodule isomorphic to
$\bigwedge^2\C^{r*}$. It is spanned by
$e_{-\varepsilon_i-\varepsilon_j}$, $k+1\leqslant i\leqslant k+r$.
On the other hand, any $\m_1$-submodule of $\Rad_u(\q)$ isomorphic
to $\C^{*r}$ is contained in $\h$. Vectors of the form $[\xi,\eta]$,
where $\xi,\eta$ lie in the $\C^{r*}$-isotypical component of the
$\m_1$-module $\Rad_u(\q)$, generate a submodule of $\Rad_u(\q)$
isomorphic to $\bigwedge^2\C^{r*}$.

{\it The case  $(\so_7,G_2,V(\pi_1))$.} Here
$\m_1=\g^{(\alpha_{n-2},\alpha_{n-1},\alpha_n)}$. Conjugating $\h$
by an appropriate element of $G^{(\alpha_1,\ldots,\alpha_{k-1})}$,
we may assume that $\h$ coincides with the subalgebra indicated in
1). In particular, $\h$ is tamely contained in the antistandard
parabolic corresponding to simple roots
$\alpha_{n-3},\ldots,\alpha_n$. Applying Corollary~\ref{Cor:3.4.3}
and then using Theorem~\ref{Thm:5.0.2}, we obtain that
$s_{\alpha_{n-3}}\in W(\g,\h)\cap
W(\g^{(\alpha_{n-3},\ldots,\alpha_n)})$ whence $I=\{n-1\}$.

{\it The case
$(\sl_r,\sl_r,\bigwedge^2\C^{r*}\oplus\C^{r*}),r\geqslant 3$.} Let
$V_0$ denote the $\m_1$-submodule of $\Rad_u(\q)$ spanned by
$e_{-\varepsilon_i-\varepsilon_j}$, $k+1\leqslant i<j\leqslant k+r$.
The isotypical component of type $\C^{r*}$ in the $\m_1$-module
$\Rad_u(\q)$ is the direct sum of $\Span_\C(e_{-\varepsilon_i})$,
$\Span_\C(e_{-\varepsilon_i-\varepsilon_j}), j\leqslant k$,
$\Span_\C(e_{-\varepsilon_i\pm\varepsilon_j}), j>k+r$, where $i$
ranges from $k+1$ to $k+r$. If $k+r=n-1$, then, conjugating  $\h$ by
an element from $G^{(\alpha_n)}$, one may assume that
$\Span_\C(\varepsilon_{-i},i=\overline{k+1,k+r})\subset \h$. Since
$[\h,\h]\subset\h$, it follows that $V_0\subset \h$. But if $r>3$,
then $V_0$ is a unique submodule of $\Rad_u(\q)$ isomorphic to
$\bigwedge^2\C^{r*}$. So $k+r=n$ provided $r>3$. Further, if $k+r=n$
and $V_0\not\subset\h$, then
$\Span_\C(e_{-\varepsilon_i-\varepsilon_j}, i=\overline{k+1,k+r})$
is the isotypical component of type $\C^{r*}$ in the $\h_1$-module
$\Rad_u(\h)$. Otherwise, the commutators of elements from the
isotypical $\C^{r*}$-component in $\Rad_u(\h)$ generate $V_0$.

Suppose   $\Rad_u(\h)$ is such as indicated in 2). By the previous
paragraph, this is always the case when  $r>3$. Put
$\widetilde{\h}=\t+\g^{(\alpha_1,\ldots,\alpha_{k-1})}+\h+V_0$. By
\cite{Wasserman},
$\a(\g,\widetilde{\h})=\Span_\C(\varepsilon_1,\varepsilon_{k+1})$
and $W(\g,\widetilde{\h})=B_{\{1,k+1\}}$. Using
Corollary~\ref{Cor:1.2.2}, we see that
$s_{\varepsilon_1-\varepsilon_{k+1}}\subset W(\g,\h)$ whence
$I=\{k+2,k+4,\ldots\}$.

Now suppose $r=3$ and $\h$ is not conjugate to the subalgebra from
2). Let us check that
 $h:=\sum_{i=k+1}^{k+3}\varepsilon_i\not\in\n_\g(\h)$.  Assume the converse. Note that $V_0\subset\Rad_u(\h)$.
 The subalgebra
$\m_1+\C h$ acts trivially on $\Rad_u(\q)^{\m_1}$ and as $\gl_3$ on
$\Rad_u(\q)/(\Rad_u(\h)+\Rad_u(\q)^{\m_1})$. Applying
Proposition~\ref{Prop:1.5.4} to $\q$, one easily gets that
$\a(\g,\h+\C h)\subset \varepsilon_{k+2}^\perp$ whence
$\varepsilon_{k+2}\in\t^{W(\g,\h)\perp}$.

If $V_0\subset \h$ (in particular, if $k+r=n-1$), then
$h\in\n_\g(\h)$. So we may assume that $k=n-3, r=3$ and $h\not\in
\n_\g(\h)$. Here $\h$ is
$G^{(\alpha_1,\ldots,\alpha_{k-1})}$-conjugate to the subalgebra
indicated in 3) so we assume that $\h$ coincides with that
subalgebra. Put
$\tau(t):=diag(t,\ldots,t,1,\ldots,1,t^{-1},\ldots,t^{-1})$ ($t$
occurs $k$ times). The limit $\lim_{t\rightarrow
0}\tau(t)\Rad_u(\h)$ exists and is the subalgebra from 2). Using
Proposition~\ref{Prop:1.2.9}, we get $I=\{n-1\}$.

{\it The case $(\sl_r,\sp_r,\C^r)$, $r$ is even.} Let us introduce
some notation. Put
$V_{21}:=\Span_\C(e_{\varepsilon_i-\varepsilon_j},i=\overline{k+1,k+r},
j=\overline{1,k})$,
$V_{31}:=\Span_\C(e_{-\varepsilon_j},e_{-\varepsilon_j\pm\varepsilon_i},j=\overline{1,k},i>k+r)$,
$V_{41}:=\Span_\C(e_{-\varepsilon_i-\varepsilon_j},i=\overline{k+1,k+r},
j=\overline{1,k})$,
$V_{51}:=\Span_\C(e_{-\varepsilon_i-\varepsilon_j}, 1\leqslant
i<j\leqslant k)$,
$V_{32}:=\Span_\C(e_{-\varepsilon_j},e_{-\varepsilon_j\pm\varepsilon_i},j=\overline{k+1,k+r},i>k+r)$,
$V_{42}:=\Span_\C(e_{-\varepsilon_i-\varepsilon_j}, k+1\leqslant
i<j\leqslant k+r)$. Further, let $V_{42}^+,V_{42}^0$ denote the
nontrivial and the trivial parts of the $\h_1$-module $V_{42}$,
respectively. Clearly,
\begin{equation}\label{eq:8.5:6}
\Rad_\u(\q)=V_{21}\oplus V_{31}\oplus V_{41}\oplus V_{51}\oplus
V_{32}\oplus V_{42}.\end{equation} The isotypical component of type
$\C^r$ (resp., of type $V(\pi_2)$, the trivial one) of the
$\h_1$-module $\Rad_u(\q)$ coincides with $V_{21}\oplus V_{41}\oplus
V_{32}$ (resp., $V_{42}^+$, $V_{31}\oplus V_{51}\oplus V_{42}^0$).
So $V_{42}^+\subset \Rad_u(\h)$.

Note that
\begin{equation}\label{eq:8.5:4}
\begin{split}
&[V_{32}',V_{32}']=V_{42},\\
&[V_{31}',V_{31}']=V_{51},\\
&[V_{21},V_{32}]=V_{31},\\
&[V_{21},V_{42}^0]=V_{41}\\
&[V_{21},V_{42}^+]=V_{41}, r>2,\\
&[V_{21},V_{41}]=V_{51},\\
&[V_{32},V_{31}]=V_{41},\\
&[V_{ij},V_{pq}]=0 \text{ for the other pairs }(i,j),(p,q).
\end{split}
\end{equation}

In (\ref{eq:8.5:4}) we denote by $V_{32}'$, resp. $V_{31}'$, either
the whole module $V_{32}$, resp. $V_{31}$, (for $k+r=n$), or the
direct sum of any two $\h_1$-, resp.,
$\g^{(\alpha_1,\ldots,\alpha_{k-1})}$-submodules in $V_{32}$,
resp.., in $V_{31}$, (for $k+r=n-1$).

Suppose $h:=\sum_{i=k+1}^{k+r}\varepsilon_k\in\n_\g(\h)$. In this
case either $V_{21}\subset\Rad_u(\h)$ or $V_{32}\oplus V_{41}\subset
\Rad_u(\h)$.

If $V_{42}^0\subset\h$,  then the nontrivial part of the $\h_1+\C
h$-module $\Rad_u(\q)/\Rad_u(\h)$ does not contain a trivial
$\h_1$-submodule. In this case  $\a(\g,\h+\C h)\neq \t$. The latter
is deduced from Proposition~\ref{Prop:1.5.4} applied to $\q$. So
$V_{42}^0\not\subset \h$. Since $[\h,\h]\subset\h$, we get
\begin{equation}\label{eq:8.5:7}
\Rad_u(\h)\subset V_{21}\oplus V_{31}\oplus V_{41}\oplus
V_{51}\oplus V_{32}\oplus V_{42}^+.
\end{equation}

Now let us suppose that (\ref{eq:8.5:7}) holds. By (\ref{eq:8.5:4}),
the projection of $\Rad_u(\h)$ to $V_{32}$ (with respect to
decomposition (\ref{eq:8.5:6})) is trivial whence $\h$ is the
subalgebra described in 4). Let us check that  $1\sim k+1$ for this
$\h$. By Proposition \ref{Prop:1.2.9} and computations above it is
enough to show that $\h\in \overline{G\h^1}$, where $\h^1$ is the
subalgebra 2). Let $\xi\in V_{42}^0$ be nonzero. Then
$\lim_{t\rightarrow \infty}\exp(t\xi)\h^1=\h$.

So it remains to consider $\h$ such that (\ref{eq:8.5:7}) does not
hold and
\begin{equation}\label{eq:8.5:12}\sum_{i=k+1}^{k+r}\varepsilon_i\not\in
\n_\g(\h).\end{equation} One can deduce from
(\ref{eq:8.5:4}),(\ref{eq:8.5:12}) that
\begin{equation}\label{eq:8.5:13} V_{51}\oplus V_{41}=[(V_{31}\oplus V^0_{42})\cap\Rad_u(\h),
(V_{31}\oplus V^0_{42}\oplus V_{21}\oplus
V_{32})\cap\Rad_u(\h)].\end{equation} From  (\ref{eq:8.5:13}) it
follows that
\begin{equation}\label{eq:8.5:8}V_{51}\oplus V_{41}\oplus V_{42}^+\subset \Rad_u(\h).\end{equation}
Further, if $k+r=n-1$, then (\ref{eq:8.5:4}) implies
\begin{equation}\label{eq:8.5:10} V_{42}^0\subset\Rad_u(\h).\end{equation}
It follows from (\ref{eq:8.5:13})-(\ref{eq:8.5:10}) that
\begin{equation}\label{eq:8.5:9}\Rad_u(\q)/\Rad_u(\h)\cong\C^r.\end{equation}

Suppose  (\ref{eq:8.5:8})-(\ref{eq:8.5:9}) hold. In this case
$\n_\g(\h)$ contains
$h_1:=2\sum_{i=1}^k\varepsilon_i+\sum_{i=k+1}^r\varepsilon_i$. As
above, one may check that $\a(\g,\h+\C h_1)\neq \t$ whence
$\t^{W(\g,\h)}\neq \{0\}$.

So $k+r=n$ and $\Rad_u(\h)$ projects surjectively onto
$V_{32},V_{21},V_{42}^0,V_{31}$. From this observation and
(\ref{eq:8.5:8}) it follows that $\h$ is
$G^{(\alpha_1,\ldots,\alpha_{k-1})}$ conjugate to the subalgebra 5).

Put
$\tau(t)=diag(1,\ldots,1,t,\ldots,t,1,t^{-1},\ldots,t^{-1},1,\ldots,1)$
($t$ occurs $r$ times). Let $\h^1$ denote $\lim_{t\rightarrow
0}\tau(t)\h$. Note that $\h^1$ is tamely contained in the
antistandard parabolic subalgebra corresponding to $\Pi(\g)\setminus
\{\alpha_n,\alpha_{n-r-1}\}$. Applying Corollary~\ref{Cor:3.4.3} to
this subalgebra and using Theorem~\ref{Thm:5.0.2}, we see that
$s_{\varepsilon_k-\varepsilon_{k+2}}\in W(\g,\h^1)$. By
Proposition~\ref{Prop:1.2.9}, $k\sim k+2$ whence
$I={k+1,k+3,\ldots,n-1}$.
\end{proof}

\subsection{Type $F_4$}\label{SUBSECTION_Weyl_homogen6}
In this subsection $\g=F_4$ and $\h$ denotes a subalgebra of $\g$
such that $\a(\g,\h)=\t, \t^{W(\g,\h)}=\{0\}$ and a maximal
reductive subalgebra of $\h$ is semisimple.

Firstly, we describe all possible sets $\widehat{\Pi}(\g,\h)$.

\begin{Lem}\label{Lem:8.6.1}
Suppose $W(\g,\h)\neq W(\g)$. Then
$\widehat{\Pi}(\g,\h)=\{\alpha_1,\alpha_2,\alpha_2+\alpha_3,
\alpha_4\}$ or $\{\alpha_1,\alpha_2,\alpha_2+\alpha_3,
\alpha_3+\alpha_4\}$.
\end{Lem}
\begin{proof}
If follows from Proposition \ref{Prop:8.2.1} that
$\widehat{\Pi}(\g,\h)\subset
\{\alpha_1,\alpha_2,\alpha_3,\alpha_4,\alpha_1+\alpha_2,\alpha_2+\alpha_3,\alpha_3+\alpha_4\}$.
Let us check that $\alpha_1+\alpha_2\not\in \widehat{\Pi}(\g,\h)$.
Assume the converse. Then $\alpha_1,\alpha_2\not\in
\widehat{\Pi}(\g,\h)$ because
$(\alpha_1,\alpha_1+\alpha_2)=(\alpha_2,\alpha_1+\alpha_2)>0$. Since
$\Span_\C(\widehat{\Pi}(\g,\h))=\t$, we obtain $\alpha_2+\alpha_3\in
\widehat{\Pi}(\g,\h)$. As above, $\alpha_3\not\in
\widehat{\Pi}(\g,\h)$ whence $\alpha_3+\alpha_4,\alpha_4\in
\widehat{\Pi}(\g,\h)$. This is not possible, for
$(\alpha_3+\alpha_4,\alpha_4)>0$. So $\alpha_1+\alpha_2\not\in
\widehat{\Pi}(\g,\h)$ is proved. Thus $\alpha_1\in
\widehat{\Pi}(\g,\h)$.

Let us check that $\alpha_2\in \widehat{\Pi}(\g,\h)$. Assume the
converse. Then $\alpha_2+\alpha_3\in
\widehat{\Pi}(\g,\h),\alpha_3+\alpha_4\not\in\widehat{\Pi}(\g,\h)$
(the latter stems from $\#\widehat{\Pi}(\g,\h)=4,
(\alpha_3+\alpha_4,\alpha_3)=(\alpha_3+\alpha_4,\alpha_4)>0$).
Therefore
$\widehat{\Pi}(\g,\h)=\{\alpha_1,\alpha_2+\alpha_3,\alpha_3,\alpha_4\}$.
But the Weyl group corresponding to this simple root system
coincides with $W(\g)$. Contradiction.

Let us check that $\alpha_3\not\in \widehat{\Pi}(\g,\h)$. Indeed,
otherwise $\alpha_4\in \widehat{\Pi}(\g,\h)$, which contradicts
$W(\g,\h)\neq W(\g)$. Summarizing, we get
$\{\alpha_1,\alpha_2\}\subset \widehat{\Pi}(\g,\h)\subset
\{\alpha_1,\alpha_2,\alpha_2+\alpha_3,\alpha_3+\alpha_4,\alpha_4\}$.
Now the claim of the proposition is easy.

\end{proof}

The both possible subgroups $W(\g,\h)\subsetneq W(\g)$ are maximal
among subgroups of $W(\g)$ generated by reflections. So the
assertion of Lemma~\ref{Lem:8.5.3} holds. Let
$\underline{\h},\h_1,\m_1$ have the same meaning as in Lemma
\ref{Lem:8.5.3}. Recall that $\alpha_3\in \Delta(\m_1)$. If
$\alpha_4\in \Delta(\m_1)$, then $s_{\alpha_4}\in W(\g,\h)$ iff
$s_{\alpha_4}\in W(\m_1,\h_1,\Rad_u(\q)/\Rad_u(\h))$, thanks to
Corollary \ref{Cor:3.4.3}. Replacing $\h$ with $\underline{\h}$,  we
may assume that
\begin{itemize}
\item[(a)] $\a(\g,\h)=\t,\t^{W(\g,\h)}=0$.
\item[(b)] A maximal reductive subalgebra $\h_1\subset\h$ is simple
and the Lie algebra $\Rad_u(\h)$ is generated by
$[\h_1,\Rad_u(\h)]$.
\item[(c)] $\m_1=\g^{(\alpha_i,\alpha_{i+1},\ldots,\alpha_3)},i=1,2,3$.
\item[(d)] $\h_1$ is not contained
in a proper Levi subalgebra of $\m_1$.
\end{itemize}

\begin{Prop}\label{Prop:8.6.2}
Let $\h,\h_1,\q,\m_1$ be such as above. Then $W(\g,\h)\neq W(\g)$
precisely in the following cases:

1) $\h_1=\g^{(\alpha_2,\alpha_3)}$,  $\Rad_u(\h)$ is spanned by
$e_{\alpha},
\alpha=-\varepsilon_1-\varepsilon_2,-\varepsilon_i,-\varepsilon_i\pm\varepsilon_j,i=1,2,j=3,4$.
In this case
$\widehat{\Pi}(\g,\h)=\{\alpha_1,\alpha_2,\alpha_2+\alpha_3,\alpha_3+\alpha_4\}$.

2) $\h_1=\g^{(\alpha_3)}$,  $[\h_1,\Rad_u(\h)]$ is spanned by
$e_{\alpha}$ with
$\alpha=-\varepsilon_i\pm\varepsilon_j,i=1,2,j=3,4,
(-\varepsilon_1-\varepsilon_2\pm\varepsilon_3\mp\varepsilon_4)/2$
and by $v_1:=x_1e_{-\varepsilon_3}+y_1
e_{(-\varepsilon_1+\varepsilon_2-\varepsilon_3+\varepsilon_4)/2},
v_2:=x_2e_{-\varepsilon_4}+y_2
e_{(-\varepsilon_1+\varepsilon_2+\varepsilon_3-\varepsilon_4)/2}$,
where $(x_1,y_1),(x_2,y_2)$ are nonzero pairs of complex numbers
such that $\Span_\C(v_1,v_2)$ is an $\h_1$-submodule in $\g$. Here
we have three conjugacy classes of subalgebras corresponding to the
collections $x_i,y_i$ with $x_1=x_2=0$, $y_1=y_2=0$,
$x_1x_2y_1y_2\neq 0$. The equality
$\widehat{\Pi}(\g,\h)=\{\alpha_1,\alpha_2,\alpha_2+\alpha_3,\alpha_4\}$
holds.

3) $\h_1=\g^{(\alpha_3)}$, $[\h_1,\Rad_u(\h)]$ is spanned by
$e_{\alpha}$ with $\alpha=-\varepsilon_1\pm\varepsilon_i,
-\varepsilon_2-\varepsilon_i,i=3,4,
(-\varepsilon_1-\varepsilon_2\pm\varepsilon_3\mp\varepsilon_4)/2$
and by $v_1:=e_{\varepsilon_4-\varepsilon_2}+x_1
e_{-\varepsilon_3}+y_1e_{(-\varepsilon_1+\varepsilon_2-\varepsilon_3+\varepsilon_4)/2}$,
$v_2:=e_{\varepsilon_3-\varepsilon_2}+x_2 e_{-\varepsilon_4}+y_2
e_{(-\varepsilon_1+\varepsilon_2-\varepsilon_4+\varepsilon_3)/2}$,
$u_1:=x_1'e_{-\varepsilon_3}+y_1'e_{(-\varepsilon_1+\varepsilon_2-\varepsilon_3+\varepsilon_4)/2}$,
$u_2:=x_2'e_{-\varepsilon_4}+y_2'e_{(-\varepsilon_1+\varepsilon_2-\varepsilon_4+\varepsilon_3)/2}$.
Here $x_i,y_i,x_i',y_i', i=1,2$ are complex numbers such that
$(x_1,y_1)$ and $(x_1',y_1')$ are linearly independent and
$\Span_\C(v_1,v_2),\Span_\C(u_1,u_2)$ are $\h_1$-submodules. There
are also three classes of conjugacy corresponding to the collections
$x_i,y_i,x_i',y_i'$ with $x_1'=x_2'=0$, $y_1'=y_2'=0$,
$x_1'x_2'y_1'y_2'\neq 0$. Finally,
$\widehat{\Pi}(\g,\h)=\{\alpha_1,\alpha_2,\alpha_2+\alpha_3,\alpha_3+\alpha_4\}$
\end{Prop}
\begin{proof}
If $\rank\m_1>1$, then
$W(\m_1,\h_1,\Rad_u(\q)/\Rad_u(\h))=W(\g,\h)\cap M_1/T$ is generated
by all reflections corresponding to short roots in $\Delta(\m_1)$.
By our conventions, $\h_1$ is not contained in a proper Levi
subalgebra of $\m_1$. Using Theorem~\ref{Thm:5.0.2}, we see that
$\h_1=\m_1$ and the nontrivial part of the $\m_1$-module
$\Rad_u(\q)/\Rad_u(\h)$ is isomorphic to the direct sum of $r$
copies of the tautological $\m_1\cong\sp_{2r}$-module.

If $\m_1=\g^{(\alpha_1,\alpha_2,\alpha_3)}\cong \sp_6$, then the
nontrivial part of the $\m_1$-module $\Rad_u(\q)/\Rad_u(\h)$
coincides with $V(\pi_3)$. Thus in this case $W(\g,\h)=W(\g)$.

Now let $\m_1=\g^{(\alpha_2,\alpha_3)}\cong\sp_4$. The multiplicity
of $\C^4$ in $\Rad_u(\q)$ equals 2. It follows that $\h$ is the
subalgebra indicated in 1). The ideal of $\z(\n_\g(\h)/\h)$
consisting of all semisimple elements is one-dimensional, its
preimage $\widetilde{\h}$ in $\n_\g(\h)$ is spanned by $\h$ and
$\varepsilon_1-\varepsilon_2$. Applying
Proposition~\ref{Prop:1.5.4}, we see that
$\a(\g,\widetilde{\h})=\t$. By Proposition \ref{Prop:8.3.6}, we get
$\t^{W(\g,\h)}=\{0\}$.

Now we check
$\widehat{\Pi}(\g,\h)=\{\alpha_1,\alpha_2,\alpha_2+\alpha_3,\alpha_3+\alpha_4\}$.
Note that $\h\in \overline{G\g^{(\alpha_2,\alpha_3,\alpha_4)}}$.
Indeed, $\h=\lim_{t\rightarrow
\infty}\exp(te_{-\varepsilon_1-\varepsilon_2})\g^{(\alpha_2,\alpha_3,\alpha_4)}$.
Proposition \ref{Prop:1.2.9} and Theorem \ref{Thm:5.0.2} yield the
claim.

Till the end of the proof $\m_1=\g^{(\alpha_3)}$. Put
$V_1:=\Span_\C(e_{-\varepsilon_3},e_{-\varepsilon_4})$,
$V_{2i}^{\pm}:=\Span_{\C}(e_{-\varepsilon_i\pm\varepsilon_3},
e_{-\varepsilon_i\pm\varepsilon_4}),$ $i=1,2$,
$V_3^{\pm}:=\Span_\C(e_{(-\varepsilon_1\pm\varepsilon_2+\varepsilon_3-\varepsilon_4)/2},
e_{(-\varepsilon_1\pm\varepsilon_2-\varepsilon_3+\varepsilon_4)/2})$,
$V_{21}:=V_{21}^+\oplus V_{21}^-, V_2^+:=V_{21}^+\oplus V_{22}^+,
V_{22}:=V_{22}^+\oplus V_{22}^-,V_2^-=:V_{21}^-\oplus V_{22}^-,
V_2:=V_{21}\oplus V_{22}, V_3:=V_3^+\oplus V_3^-$. The isotypical
component of type $\C^2$ in $\Rad_u(\q)$ coincides with $V_1\oplus
V_2\oplus V_3$.


One easily verifies that

\begin{equation}\label{eq:8.6:1}
\begin{split}
&[V_1,V_1]=\C e_{-\varepsilon_3-\varepsilon_4},\\
&[V_1,V_{2i}^+]=\C e_{-\varepsilon_i},\\
&[V_1,V_3^{\pm}]=\C e_{(-\varepsilon_1\pm\varepsilon_2-\varepsilon_3-\varepsilon_4)/2},\\
&[V_{21}^\pm,V_{22}^\mp]=[V_3^{-},V_3^{-}]=\C e_{-\varepsilon_1-\varepsilon_2},\\
&[V_{22}^+,V_3^+]=\C e_{(-\varepsilon_1-\varepsilon_2+\varepsilon_3+\varepsilon_4)/2},\\
&[V_{22}^-,V_3^+]=\C e_{(-\varepsilon_1-\varepsilon_2-\varepsilon_3-\varepsilon_4)/2},\\
&[V_3^+,V_3^+]=\C e_{\varepsilon_2-\varepsilon_1},\\
&[V_3^+,V_3^-]=\C e_{-\varepsilon_1}.
\end{split}
\end{equation}

The other commutators of $V_1,V_{2i}^{\pm},V_3^\pm$ vanish.

Let us check that $e_{-\varepsilon_1-\varepsilon_2}\in \Rad_u(\h)$.
Indeed, according to (\ref{eq:8.6:1}), one may consider the
commutator on $V_{2}$ as a nondegenerate skew-symmetric form. The
intersection $\Rad_u(\h)\cap V_2$ is 6-dimensional whence not
isotropic.

Let us show that
$e_{-(\varepsilon_1+\varepsilon_2+\varepsilon_3+\varepsilon_4)/2}\in
\Rad_u(\h)$. The dimension of $\Rad_u(\h)\cap (V_1\oplus
V_{22}^-\oplus V_3)$ is not less than 6. By (\ref{eq:8.6:1}), it is
enough to consider the situation when this dimension equals 6 and
the projection of $\Rad_u(\h)\cap (V_1\oplus V_{22}^-\oplus V_3)$ to
any of four submodules $V_1,V_{22}^-, V_3^{\pm}$ is nonzero. If
$V_{22}^-\subset \Rad_u(\h)$, then
$e_{-(\varepsilon_1+\varepsilon_2+\varepsilon_3+\varepsilon_4)/2}$
lies in the commutator of $V_{22}^-$ with an $\m_1$-submodule of
$\Rad_u(\h)$ not contained in $V_1\oplus V_{22}^-\oplus V_3^-$. So
we may assume that $\Rad_u(\h)\cap (V_1\oplus V_{22}^-\oplus V_3^-)$
contains a submodule projecting nontrivially to $V_1$ and
$\Rad_u(\h)\cap (V_{22}^-\oplus V_3^-)$ contains a submodule
projecting nontrivially to $V_3^-$. The commutator of two such
submodules contains an element of the form
$e_{-(\varepsilon_1+\varepsilon_2+\varepsilon_3+\varepsilon_4)/2}+ x
e_{-\varepsilon_1-\varepsilon_2}$. It remains to recall that
$e_{-\varepsilon_1-\varepsilon_2}\in \Rad_u(\h)$.

Analogously (by considering $V_1\oplus V_{21}^+\oplus V_3$ instead
of $V_1\oplus V_{22}^-\oplus V_3$), we get
$e_{-\varepsilon_1}\in\Rad_u(\h)$.

Now we show that  $V_2^-\oplus V_{21}^+\oplus
V_3^-\subset\Rad_u(\h)$. Assume the converse. Then the projection of
$\Rad_u(\h)$ to $V_1\oplus V_{22}^+\oplus V_3^+$ is surjective. Note
that the commutator of $V_1\oplus V_2\oplus V_3$ and $V_2^-\oplus
V_{21}^+\oplus V_3^-$ lies in
$\Span_\C(e_{-\varepsilon_1-\varepsilon_2},e_{(-\varepsilon_1-
\varepsilon_2-\varepsilon_3-\varepsilon_4)/2}, e_{-\varepsilon_1})$.
It follows from (\ref{eq:8.6:1}) that $\Rad_u(\h)^{\m_1}$ is spanned
by the r.h.s.'s of (\ref{eq:8.6:1}). From this one can deduce that
$V_2^-\oplus V_{21}^+\oplus V_3^-\subset
[\Rad_u(\h)^{\m_1},\Rad_u(\h)]$. Contradiction.

So $\Rad_u(\h)=V_{21}\oplus V_{22}^-\oplus V_3^-\oplus U$, where
$U\subset V_1\oplus V_{22}^+\oplus V_3^+$ is a fourdimensional
$\m_1$-submodule. (\ref{eq:8.6:1}) implies
$$[U,U]\subset V':=\Span_\C(e_{-\varepsilon_3-\varepsilon_4},
e_{-\varepsilon_2},e_{\varepsilon_2-\varepsilon_1},e_{(-\varepsilon_1-\varepsilon_2+\varepsilon_3+\varepsilon_4)/2},
e_{(-\varepsilon_1\pm\varepsilon_2-\varepsilon_3-\varepsilon_4)/2}).$$

Note that $\dim [U,U]=2$ if $V_{22}^+\subset U$ and 3 otherwise. The
subalgebra $\h$  such that $U\subset\Rad_u(\h)$ is denoted by
$\h_U$. The torus $T_0$ with the Lie algebra
$\t_0:=(\varepsilon_3-\varepsilon_4)^\perp\subset \t$ acts on
$V_1\oplus V_{22}^+\oplus V_3^+$ and on $V'$. If $U,U_0$ are
fourdimensional $\m_1$-submodules in $V_1\oplus V_{22}^+\oplus
V_3^+$ such that $\dim [U,U]=\dim [U_0,U_0]$ and $U_0\in
\overline{T_0U}$, then $\h_{U_0}\in \overline{T_0\h_U}$.

Let us suppose for a moment that $\t\subset\n_\g(\h_U)$. Let
$\Sigma_U$ denote the set of all roots $\beta$ with
$\g^{\beta}\subset [V_1\oplus V_2\oplus V_3,V_1\oplus V_2\oplus
V_3]\setminus \Rad_u(\h)$. For any $\Sigma\subset \Sigma_U$ the
subspace $\widetilde{\h}:=\h_U\oplus
\bigoplus_{\beta\in\Sigma}\g^{\beta}\subset\g$ is a subalgebra of
$\g$ normalized by $\t$. The equality $\t^{W(\g,\h)}=\{0\}$ holds
iff $\a(\g,\widetilde{\h}+\t_0)=\t$ (Proposition~\ref{Prop:8.3.6}).
Applying Proposition~\ref{Prop:1.5.4} to $\q$, one can show that
$\a(\g,\widetilde{\h}+\t_0)=\t$ iff the subspace spanned by
$\Sigma_U\setminus \Sigma$ and
$(-\varepsilon_1+\varepsilon_2+\varepsilon_3+\varepsilon_4)/2$ is
threedimensional.

There are  seven  $T_0$-orbits of submodules $U\subset V_1\oplus
V_{22}^+\oplus V_3^+$. To the $T_0$-orbit of $U$ we associate the
triple $(x_1,x_2,x_3)$ consisting of 0 and 1 by the following rule:
 $V_1\subset U$ (resp. $V_{22}^+\subset U, V_3^+\subset U$) iff $x_1=0$ (resp., $x_2=0,x_3=0$).
Abusing the notation, we write  $\h_{(x_1,x_2,x_3)}$ instead of
$\h_U$. The $T_0$-orbit corresponding to $(x_1,x_2,x_3)$ lies in the
closure of that corresponding to $(y_1,y_2,y_3)$ iff $x_i\leqslant
y_i$. Consider the seven possible variants case by case.

{\it The case $(1,0,0)$.} Here $\t\subset\n_\g(\h)$ and
$\Sigma_U=\{-\varepsilon_3-\varepsilon_4,-\varepsilon_2,(-\varepsilon_1+\varepsilon_2
-\varepsilon_3-\varepsilon_4)/2\}$. By above, $\t^{W(\g,\h)}=\{0\}$.
Let us check that
$\widehat{\Pi}(\g,\h)=\{\alpha_1,\alpha_2,\alpha_2+\alpha_3,\alpha_4\}$.
Let $\widetilde{\h}$ be the subspace of $\g$ spanned by $\h$ and
$e_{(-\varepsilon_1+\varepsilon_2\pm(\varepsilon_3+\varepsilon_4))/2}$,
it is a subalgebra. Applying Proposition~\ref{Prop:1.5.4}, we see
that $\a(\g,\widetilde{\h})=\t$. Let $\widetilde{\q}$ denote the
antistandard parabolic subalgebra of $\g$ corresponding to the
simple roots $\alpha_2,\alpha_3,\alpha_4$. Note that
$\Rad_u(\widetilde{\q})\subset\widetilde{\h}$. Therefore
$W(\g,\widetilde{\h})=W(\g^{(\alpha_2,\alpha_3,\alpha_4)},\widetilde{\h}/\Rad_u(\widetilde{\q}))$
(Corollary~\ref{Cor:3.4.7}). The last group was computed in
Subsection~\ref{SUBSECTION_Weyl_aff3}, it is generated by
$s_{\alpha_2},s_{\alpha_2+\alpha_3},s_{\alpha_4}$. It remains to
recall that $W(\g,\widetilde{\h})\subset W(\g,\h)$.

{\it The case $(1,0,1)$.} As we have remarked above,
$\h_{(1,0,0)}\in \overline{T_0\h}$ whence
$W(\g,\h)=W(\g,\h_{(1,0,0)})$.

{\it The case $(0,0,1)$.} Here $\h\in \overline{T_0\h_{(1,0,1)}}$.
By above, $\t^{W(\g,\h)}=\{0\}$
($\Sigma_U=\{\varepsilon_2-\varepsilon_1,
-\varepsilon_1+\varepsilon_2-\varepsilon_3-\varepsilon_4)/2,
(-\varepsilon_1-\varepsilon_2 +\varepsilon_3+\varepsilon_4)/2\}$).

{\it The case $(0,1,0)$.} We have
$\Sigma_U=\{-\varepsilon_2,(-\varepsilon_1-\varepsilon_2
+\varepsilon_3+\varepsilon_4)/2\}\}$. Thus $\t^{W(\g,\h)}\neq
\{0\}$.

{\it The case $(1,1,0)$.} The subalgebra
$\t_1=\Span_\C(\varepsilon_1,2\varepsilon_2-\varepsilon_3-\varepsilon_4)$
is contained $\n_\g(\h)$. On the other hand,
$\n_\g(\h+\t_1)/(\h+\t_1)$ does not contain a semisimple element,
since $\n_\g(\m_1+\t_1)=\m_1+\t$. To verify $\t^{W(\g,\h)}=\{0\}$ it
is enough to check $\a(\g,\h+\t_1)=\t$. This is done by using
Proposition~\ref{Prop:1.5.4}. Let us show that
$\widehat{\Pi}(\g,\h)=\{\alpha_1,\alpha_2,\alpha_2+\alpha_3,\alpha_3+\alpha_4\}$.
Let $\widetilde{\q}$ denote the antistandard parabolic subgroup of
$\g$ corresponding to $\alpha_2,\alpha_3,\alpha_4$. It is enough to
check that $s_{\alpha_4}\not\in W(\g,\h)\cap
W(\g^{(\alpha_2,\alpha_3,\alpha_4)})$. It follows from
\cite{comb_ham}, Proposition 8.2, that $W(\g,\h)\cap
W(\g^{(\alpha_2,\alpha_3,\alpha_4)})=W_{G^{(\alpha_2,\alpha_3,\alpha_4)},
\widetilde{Q}/H}$. By Proposition \ref{Prop:8.5.4},
$s_{\alpha_4}\not\in W(\g^{(\alpha_2,\alpha_3,\alpha_4)},
\g^{(\alpha_2,\alpha_3,\alpha_4)}\cap\h)$ and we are done.

{\it The case $(1,1,1)$.} Since $\h_{(1,1,0)}\in \overline{T_0\h}$,
we have
$\widehat{\Pi}(\g,\h)=\{\alpha_1,\alpha_2,\alpha_2+\alpha_3,\alpha_3+\alpha_4\}$.

{\it The case $(0,1,1)$.} Analogously to the case $(1,1,0)$, one can
show that  $\t^{W(\g,\h)}=\{0\}$. Since $\h\in
\overline{T_0\h_{(1,1,1)}}$, we see that
$\widehat{\Pi}(\g,\h)=\{\alpha_1,\alpha_2,\alpha_2+\alpha_3,\alpha_3+\alpha_4\}$.
\end{proof}

\section{Algorithms}\label{SECTION_algorithm}
\subsection{Algorithms for computing Cartan
spaces}\label{SUBSECTION_algorithm1} {\it Case 1.} Suppose $X=G/H$,
where $H$ is a reductive subgroup of $G$. Using Theorem 1.3,
\cite{ranks}, we compute  $\a(\g,\h)$. Then, applying
Proposition~\ref{Prop:6.0.1}, we compute a point in the
distinguished component of $(G/H)^{L_{0\,G,G/H}^\circ}$. Finally,
using Proposition~\ref{Prop:6.0.2}, we find the whole distinguished
component.

{\it Case 2.} Here $X=G*_HV$ is an affine homogeneous vector bundle
and $\pi:G*_HV\rightarrow G/H$ is the natural projection. Applying
the algorithm of case 1 to $G/H$, we compute the space $\a(\g,\h)$
and find a point $x$ in the distinguished component of
$(G/H)^{L_{0\,G,G/H}}$. Applying the following algorithm to the
group $L_0:=L_{0\,G,G/H}^\circ$ and the module $V=\pi^{-1}(x)$, we
compute $\a_{L_0,V}$.

\begin{Alg}\label{Alg:9.1.1} Let $G$ be a connected
reductive algebraic group and $V$ a $G$-module. Put $G_0=G, V_0=V$.
Assume that we have already constructed a pair $(G_i,V_i)$, where
$G_i$ is a connected subgroup in $G_0$, $B_i:=B\cap G_i$. Choose a
$B_i$-semiinvariant vector $\alpha\in V_i^*$. Put
$V_{i+1}=(\u_i^-\alpha)^0$, where $\u_i^-$ is a maximal unipotent
subalgebra of $\g_i$ normalized by $T$ and opposite to $\b_i$ and
the superscript $^0$ means the annihilator. Put
$G_{i+1}=Z_{G_i}(\alpha)$. The group $G_{i+1}$ is  connected and
$L_{0\,G_i,V_i}=L_{0\,G_{i+1},V_{i+1}}$. Note that $\rank
[\g_{i+1},\g_{i+1}]\leqslant \rank [\g_i,\g_i]$ with the equality
iff $\alpha\in V^{[\g_i,\g_i]}$. Thus if $[\g_i,\g_i]$ acts
non-trivially on $V$, then we may assume that $\rank
[\g_{i+1},\g_{i+1}]< \rank [\g_i,\g_i]$. So $V=V^{[\g_k,\g_k]}$ for
some $k$. Here $L_{0\,G,V}=L_{0\,G_k,V_k}$ coincides with the unit
component of the inefficiency kernel for the action $G_k:V_k$.
\end{Alg}

Thanks to Proposition~\ref{Prop:1.5.8},
$\a_{G,X}=\a(\g,\h)+\a_{L_0,V}$. Finally, using
Proposition~\ref{Prop:6.0.2}, we determine the distinguished
component of  $X^{L_{0\,G,X}^\circ}$.

{\it Case 3.} Suppose $X=G/H$, where $H$ is a nonreductive subgroup
of $G$. We find a parabolic subgroup $Q\subset G$ tamely containing
$H$ by using the following algorithm.

\begin{Alg}\label{Alg:9.1.2}
Put $\n_0=\Rad_u(\h)$. If an algebra $\n_i\subset \g, i\geqslant 0,$
consisting of nilpotent elements is already constructed, we put
$\n_{i+1}=\Rad_u(\n_\g(\n_i))$. Clearly, $\n_i\subset\n_{i+1}$. It
is known that if $\n_{i+1}=\n_i$, then $\n_i$ is the unipotent
radical of a parabolic $\q$ (see \cite{Bourbaki}, ch.8, $\S$ 10,
Theorem 2). Clearly, $\Rad_u(H)\subset \Rad_u(Q)$. It is clear from
construction that $N_G(\n_0)\subset Q$. In particular, $H\subset Q$.
\end{Alg}

Further, we choose a Levi subgroup $M\subset Q$ such that $M\cap H$
is a maximal reductive subgroup of $H$ and an element $g\in G$ such
that $gQg^{-1}$ is an antistandard parabolic subgroup and $gMg^{-1}$
is its standard Levi subgroup. Replace  $(Q,M,H)$ with
$(gQg^{-1},gMg^{-1},gHg^{-1})$. Put $X=Q^-/H$. Using
Remark~\ref{Rem:1.5.3}, we construct an $M$-isomorphism of $X$ and
an affine homogeneous vector bundle. By
Proposition~\ref{Prop:1.5.4}, $\a_{G,G/H}=\a_{M,X}$. If  $G/H$ is
quasiaffine, we use Proposition~\ref{Prop:6.0.3} and obtain a point
in the distinguished component of $(G/H)^{L_{0\,G,G/H}^\circ}$.
Applying Proposition~\ref{Prop:6.0.2}, we determine the whole
distinguished component.

\subsection{Computation of Weyl groups}\label{SUBSECTION_algor2}
{\it Case 1.} Let $X=G*_HV$, where $H$ is a reductive subgroup of
$G$,  $V$ is an $H$-module. Suppose  $\rank_G(X)=\rank(G)$. Let
$G=Z(G)^\circ G_1\ldots G_k$ be the decomposition into the locally
direct product of the center and simple normal subgroups. Put
$H_i=G_i\cap H$. By Proposition \ref{Prop:3.4.5}, the equality
$W(\g,\h,V)=\prod_{i=1}^k W(\g_i,\h_i,V)$ holds. The computation of
$W(\g_i,\h_i,V)$ is carried out by using Theorem~\ref{Thm:5.0.2}.

{\it Case 2.} Let $X=G*_HV$, where $H$ is a reductive subgroup of
$G$,  $V$ is an $H$-module. We find  $\a_{G,X}$ and a point in the
distinguished component of $\underline{X}\subset X^{L_0}$, where
$L_0=L_{0\,G,X}^\circ$, lying in $G/H$, by using the algorithm of
case 2 of Subsection~\ref{SUBSECTION_algorithm1}. We may assume that
$eH$ lies in the distinguished component. Put
$\underline{G}=(N_G(L_0)^\circ N_H(L_0))/L_0,
\underline{H}=N_H(L_0)/L_0, \underline{V}=V^{L_0}$. We have a
$\underline{G}$-isomorphism $\underline{X}\cong
\underline{G}*_{\underline{H}}\underline{V}$
(Proposition~\ref{Prop:6.0.2}). Put
$\Gamma=N_{\underline{G}}(\underline{B},\underline{T})/\underline{T}$,
where $\underline{B},\underline{T}$ are the distinguished Borel
subgroup and the maximal torus of $\underline{G}$. By
Theorem~\ref{Thm:3.0.3},
$\a_{G,X}=\underline{\t}=\a_{\underline{G}^\circ,\underline{X}},
W_{G,X}=W_{\underline{G}^\circ,\underline{X}}\leftthreetimes
\Gamma$. This reduces the computation of $W_{G,X}$ to the previous
case.

{\it Case 3.} Let $X=G/H$ be the quasiaffine homogeneous space with
$\rank_G(G/H)=\rank(G)$. Let $G=Z(G)^\circ G_1 G_2\ldots G_k$ be the
decomposition into the locally direct product of the center and
simple normal subgroups. Put $H_i:=G_i\cap H$. Then, thanks to
Proposition  \ref{Prop:3.4.5}, $W(\g,\h)=\prod W(\g_i,\h_i)$. So we
reduce the computation to the case when $G$ is simple.

If $\g\cong\sl_2$, then $W(\g,\h)$ is trivial iff $\h$ is a
one-dimensional unipotent subalgebra.

Suppose  $\g\cong\so_5,G_2$. If $\h$ does not contain the unipotent
radical of a parabolic, then $W(\g,\h)$ can be computed by using
Propositions~\ref{Prop:8.1.3}, \ref{Prop:8.1.4}. Otherwise, applying
Corollary~\ref{Cor:3.4.7}, we reduce the computation of $W(\g,\h)$
to the case $\rank [\g,\g]\leqslant 1$.

Below we suppose  $\g\not\cong \sl_2,\so_5,G_2$.

Let $\widetilde{\h}$ be the inverse image of the reductive part of
$\z(\n_\g(\h)/\h)$ under the canonical epimorphism
$\n_\g(\h)\twoheadrightarrow \n_\g(\h)/\h$. Compute
$\a(\g,\widetilde{\h})$. By Proposition~\ref{Prop:8.3.6},
$\a(\g,\widetilde{\h})=\t^{W(\g,\h)}$.

When $\g$ is of type $A,D,E$ we recover $W(\g,\h)$ from
$\t^{W(\g,\h)}$ using Proposition~\ref{Prop:8.3.1}.

Now let $\g\cong \so_{2l+1},\sp_{2l},l\geqslant 3, F_4$. Recall, see
Corollary~\ref{Cor:3.4.6} that $W(\g,\h)$ is identified with
$W(\g,\widetilde{\h})$. If $\t^{W(\g,\h)}\neq\{0\}$, then
$\a(\g,\widetilde{\h})\neq \t$ and we proceed to  case 4. Note that
under the reduction of  case 4 to  case 3 $\rank[\g,\g]$ decreases.
So we may assume that $\t^{W(\g,\h)}=\{0\}$.

Let us find a parabolic subalgebra  $\q\subset\g$ tamely containing
$\h$ (see Algorithm~\ref{Alg:9.1.2}). Conjugating $(\q,\h)$ by an
element of $G$, we may assume that $\q$ is an antistandard parabolic
subalgebra of $\g$ and that $\h_0:=\m\cap\h$ is a Levi subalgebra of
$\h$, where $\m\subset \q$ is the standard Levi subalgebra of $\q$.

If $\g\cong \sp_{2l}$, then the computation is carried out by using
Proposition~\ref{Prop:8.4.3}. In the remaining cases we inspect
Table~\ref{Tbl:5.3.10} and find simple ideals $\h_1\subset \h_0,
\m_1\subset\m$ satisfying the assumptions of Lemma~\ref{Lem:8.4.2}.
Denote by $\underline{\h}$ the subalgebra of $\g$ generated by
$\h_1$ and $[\h_1,\Rad_u(\h)]$. According to Lemmas~\ref{Lem:8.5.2},
\ref{Lem:8.5.1}, \ref{Lem:8.6.1}, $W(\g,\underline{\h})=W(\g,\h)$.
The groups $W(\g,\underline{\h})$ for $\g=\so_{2l+1}$ are computed
in Subsection~\ref{SUBSECTION_Weyl_homogen5}, see especially
Proposition~\ref{Prop:8.5.4} and the preceding discussion. The case
$\g=F_4$ is considered in Subsection~\ref{SUBSECTION_Weyl_homogen6}
(Proposition~\ref{Prop:8.6.2}).

{\it Case 4.} Here $X=G/H$ is an arbitrary homogeneous space. At
first we find a parabolic subalgebra $\q\subset \g$ tamely
containing $\h$. Let us choose a subalgebra $\q_0\subset \q$
annihilating a  $\q$-semiinvariant vector  $v\in V$, where $V$ is a
$G$-module (if $\q$ is standard, one can take a highest vector of an
appropriate irreducible module for $v$). Let $\chi\in \X(Q)$ we the
character of $\q$ such that $\xi v=\chi(\xi)v$ for all $\xi\in\q$.
If $\h\subset\q_0$, then  $\h$ is observable (i.e., the homogeneous
space $G/H$ is quasiaffine), thanks to Sukhanov's theorem,
\cite{Sukhanov}. If not, put $\widetilde{G}=G\times\C^\times$, and
embed $\q$ into $\widetilde{\g}$ via $\iota:\xi\mapsto
(\xi,-\chi(\xi))$. There is the natural representation of
$\widetilde{G}$ in $V$ such that $\widetilde{G}_v=\iota(Q)$. So
$\iota(\h)$ is observable in $\widetilde{\g}$. Further, we have a
$G$-equivariant principal $\C^\times$-bundle
$\widetilde{G}/\iota(H^\circ)\rightarrow G/H^\circ$. From
Proposition~\ref{Prop:1.2.5} it follows that
$W(\g,\h)=W(\widetilde{\g},\iota(\h))$. So  replacing $(\g,\h)$ with
$(\widetilde{\g},\iota(\h))$, if necessary, we may (and will) assume
that $\h$ observable.

We compute $\a(\g,\h)$ and determine the distinguished component
$\underline{X}\subset X^{L_0}$, where $L_{0}:=L_{0\,G,G/H}^\circ$.
Put $\underline{G}:=(N_G(L_0)^\circ N_H(L_0))/L_0,
\underline{H}:=N_H(L_0)/L_0$. Let
$\underline{B},\underline{T},\Gamma$ have the same meaning as in
case 2. According to Proposition~\ref{Prop:6.0.2},
$\underline{X}\cong \underline{G}/\underline{H}$. Using Theorem
\ref{Thm:3.0.3}, we get
$\a(\g,\h)=\underline{\t}=\a(\underline{\g},\underline{\h}),
W(\g,\h)=W(\underline{\g},\underline{\h})\leftthreetimes \Gamma$.

\bigskip
{\Small Chair of Higher Algebra, Department of Mechanics and
Mathematics, Moscow State University.

\noindent E-mail address: ivanlosev@yandex.ru}
\end{document}